\documentclass[12pt]{article}
\usepackage{amssymb,amsmath,latexsym}
\usepackage{color}

\allowdisplaybreaks

\begin{document}

\begin{doublespace}

\def\1{{\bf 1}}
\def\ind{{\bf 1}}
\def\nn{\nonumber}
\newcommand{\I}{\mathbf{1}}

\def\sA {{\cal A}} \def\sB {{\cal B}} \def\sC {{\cal C}}
\def\sD {{\cal D}} \def\sE {{\cal E}} \def\sF {{\cal F}}
\def\sG {{\cal G}} \def\sH {{\cal H}} \def\sI {{\cal I}}
\def\sJ {{\cal J}} \def\sK {{\cal K}} \def\sL {{\cal L}}
\def\sM {{\cal M}} \def\sN {{\cal N}} \def\sO {{\cal O}}
\def\sP {{\cal P}} \def\sQ {{\cal Q}} \def\sR {{\cal R}}
\def\sS {{\cal S}} \def\sT {{\cal T}} \def\sU {{\cal U}}
\def\sV {{\cal V}} \def\sW {{\cal W}} \def\sX {{\cal X}}
\def\sY {{\cal Y}} \def\sZ {{\cal Z}}

\def\bA {{\mathbb A}} \def\bB {{\mathbb B}} \def\bC {{\mathbb C}}
\def\bD {{\mathbb D}} \def\bE {{\mathbb E}} \def\bF {{\mathbb F}}
\def\bG {{\mathbb G}} \def\bH {{\mathbb H}} \def\bI {{\mathbb I}}
\def\bJ {{\mathbb J}} \def\bK {{\mathbb K}} \def\bL {{\mathbb L}}
\def\bM {{\mathbb M}} \def\bN {{\mathbb N}} \def\bO {{\mathbb O}}
\def\bP {{\mathbb P}} \def\bQ {{\mathbb Q}} \def\bR {{\mathbb R}}
\def\bS {{\mathbb S}} \def\bT {{\mathbb T}} \def\bU {{\mathbb U}}
\def\bV {{\mathbb V}} \def\bW {{\mathbb W}} \def\bX {{\mathbb X}}
\def\bY {{\mathbb Y}} \def\bZ {{\mathbb Z}}
\def\R {{\mathbb R}} \def\RR {{\mathbb R}} \def\H {{\mathbb H}}
\def\n{{\bf n}} \def\Z {{\mathbb Z}}

\newcommand{\expr}[1]{\left( #1 \right)}
\newcommand{\cl}[1]{\overline{#1}}
\newtheorem{thm}{Theorem}[section]
\newtheorem{lemma}[thm]{Lemma}
\newtheorem{defn}[thm]{Definition}
\newtheorem{prop}[thm]{Proposition}
\newtheorem{corollary}[thm]{Corollary}
\newtheorem{remark}[thm]{Remark}
\newtheorem{example}[thm]{Example}
\numberwithin{equation}{section}
\def\ee{\varepsilon}
\def\qed{{\hfill $\Box$ \bigskip}}
\def\NN{{\mathcal N}}
\def\AA{{\mathcal A}}
\def\MM{{\mathcal M}}
\def\BB{{\mathcal B}}
\def\CC{{\mathcal C}}
\def\LL{{\mathcal L}}
\def\DD{{\mathcal D}}
\def\FF{{\mathcal F}}
\def\EE{{\mathcal E}}
\def\QQ{{\mathcal Q}}
\def\SS{{\mathcal S}}
\def\RR{{\mathbb R}}
\def\R{{\mathbb R}}
\def\L{{\bf L}}
\def\K{{\bf K}}
\def\S{{\bf S}}
\def\A{{\bf A}}
\def\E{{\mathbb E}}
\def\F{{\bf F}}
\def\P{{\mathbb P}}
\def\N{{\mathbb N}}
\def\eps{\varepsilon}
\def\wh{\widehat}
\def\wt{\widetilde}
\def\pf{\noindent{\bf Proof.} }
\def\pff{\noindent{\bf Proof} }
\def\cp{\mathrm{Cap}}

\title{\Large \bf Heat kernels of non-symmetric jump processes: beyond the stable case}

\author{{\bf Panki Kim}\thanks{This work was supported by the National Research Foundation of Korea(NRF) grant funded by the Korea government(MSIP) (No. 2016R1E1A1A01941893)
}
\quad {\bf Renming Song\thanks{Research supported in part by a grant from
the Simons Foundation (208236)}} \quad and
\quad {\bf Zoran Vondra\v{c}ek}
\thanks{Research supported in part by the Croatian Science Foundation under the project 3526}
}

\date{}

\maketitle

\begin{abstract}
Let $J$ be the L\'evy density of a symmetric L\'evy process in
$\R^d$ with its L\'evy exponent satisfying a weak lower scaling condition at infinity.
Consider the non-symmetric and non-local operator
$$
\LL^{\kappa}f(x):= \lim_{\eps \downarrow 0} \int_{\{z \in \R^d: |z|>\eps\}}
(f(x+z)-f(x))\kappa(x,z)J(z)\, dz\, ,
$$
where $\kappa(x,z)$ is a 
Borel function on $\R^d\times \R^d$ satisfying
$0<\kappa_0\le \kappa(x,z)\le \kappa_1$,  $\kappa(x,z)=\kappa(x,-z)$ and
$|\kappa(x,z)-\kappa(y,z)|\le \kappa_2|x-y|^{\beta}$ for some 
$\beta\in (0, 1]$.
We construct the heat kernel $p^\kappa(t, x, y)$ of $\LL^\kappa$,  establish its upper bound as well as its fractional derivative and gradient estimates. Under an additional weak upper scaling condition at infinity, we also establish a lower bound for the heat kernel $p^\kappa$.
\end{abstract}



\section{Introduction}

Suppose that $d\ge 1$, $\alpha\in (0,2)$ and $\kappa(x,z)$ is a Borel 
function on $\R^d\times \R^d$ such that
\begin{equation}\label{e:intro-kappa}
0<\kappa_0\le \kappa(x,z)\le \kappa_1\, , \quad \kappa(x,z)=\kappa(x,-z)\, ,
\end{equation}
and 
for some $\beta\in (0,1]$,
\begin{equation}\label{e:intro-kappa-holder}
|\kappa(x,z)-\kappa(y,z)|\le \kappa_2|x-y|^{\beta}\, .
\end{equation}
The operator
\begin{equation}\label{e:intro-stable-operator}
\LL_{\alpha}^{\kappa}f(x)= \lim_{\eps \downarrow 0} 
\int_{\{z \in \R^d: |z|>\eps\}}
(f(x+z)-f(x))\frac{\kappa(x,z)}{|z|^{d+\alpha}}\, dz
\end{equation}
is a  non-symmetric and non-local stable-like operator.  
In the recent paper \cite{CZ}, Chen and Zhang studied the heat kernel of 
$\LL_{\alpha}^{\kappa}$ and its sharp two-sided estimates.
As the main result of the paper,  they proved the existence and uniqueness of a non-negative jointly continuous function $p_{\alpha}^{\kappa}(t,x,y)$ in $(t,x,y)\in (0,1]\times \R^d\times \R^d$ solving the equation
$$
\partial_t p_{\alpha}^{\kappa}(t,x,y)=\LL_{\alpha}^{\kappa} p_{\alpha}^{\kappa}(t,\cdot,y)(x)\, , \quad x\neq y\, ,
$$
and satisfying four properties - an upper bound, H\"older's estimate, fractional derivative estimate and continuity, cf.~\cite[Theorem 1.1]{CZ} for details. They also proved some other properties of the heat kernel  $p_{\alpha}^{\kappa}(t,x,y)$ such as conservativeness, Chapman-Kolmogorov equation, lower bound, gradient estimate and studied the corresponding semigroup. Their paper is the 
first one to address these questions for not necessarily symmetric non-local stable-like operators. 
These operators can be regarded as the non-local
counterpart of elliptic operators in non-divergence form. In this context the H\"older continuity of $\kappa(\cdot, z)$ in \eqref{e:intro-kappa-holder} is a natural assumption. 

The goal of this paper is to extend the results of \cite{CZ} to more general operators than the ones defined in \eqref{e:intro-stable-operator}. These operators will be non-symmetric and not necessarily stable-like. We will replace the kernel $\kappa(x,z) |z|^{-d-\alpha}$ with a kernel  
$\kappa(x,z)J(z)$ 
where $\kappa$ still satisfies \eqref{e:intro-kappa} and \eqref{e:intro-kappa-holder}, but 
$J(z)$ is the L\'evy density
of a rather general symmetric L\'evy process. 
Here are the precise assumptions that we make.

Let $\phi:(0,\infty)\to (0,\infty)$ be a Bernstein function without drift and killing. Then
$$
\phi(\lambda)=\int_{(0,\infty)}\left(1-e^{-\lambda t}\right)\mu(dt),
$$
where $\mu$ is a measure on $(0,\infty)$ satisfying 
$\int_{(0,\infty)} (t\wedge 1)\mu(dt)<\infty$. 
Here and throughout this paper, we use the notation $a \wedge b := \min \{ a, b\}$
and $a\vee b:=\max \{ a, b\}$.
Without loss of generality we assume that $\phi(1)=1$. Define $\Phi:(0,\infty)\to (0,\infty)$ by $\Phi(r)=\phi(r^2)$ and let $\Phi^{-1}$ be its inverse. The function $x\mapsto \Phi(|x|)=:\Phi(x)$, $x\in \R^d$, $d\ge 1$, is negative definite and hence it is the characteristic exponent of an isotropic L\'evy process on $\R^d$. This process can be obtained by subordinating a $d$-dimensional Brownian motion by an independent subordinator with Laplace exponent $\phi$. 
The L\'evy measure of this process has a density $j(|y|)$ where $j:(0,\infty)\to (0,\infty)$ is the function given by
$$
j(r)=\int_{(0,\infty)}(4\pi t)^{-d/2}e^{-\frac{r^2}{4t}}\mu(dt)\, .
$$
Thus we have
$$
\Phi(x)=\int_{\R^d\setminus \{0\}}\left(1-\cos (x\cdot y)\right) j(|y|)\, dy\, .
$$

Note that when $\phi(\lambda)=\lambda^{\alpha/2}$, $0<\alpha<2$, we have $\Phi(r)=r^{\alpha}$, the corresponding subordinate Brownian motion is an isotropic $\alpha$-stable process and 
$j(r)=c(d, \alpha)  \, r^{-d-\alpha}$.

Our main assumption is the following \emph{weak lower scaling condition at infinity}: There exist $\delta_1\in (0,2]$ and $a_1 \in (0,1)$ such that
\begin{equation}\label{e:intro-wsc}
a_1 \lambda^{\delta_1}\Phi(r)\le \Phi(\lambda r)\, ,\quad \lambda\ge 1, r\ge 1\, .
\end{equation}
This condition implies that $\lim_{\lambda\to \infty}\Phi(\lambda)=\infty$ and hence $\int_{\R^d\setminus \{0\}}j(|y|)dy=\infty$ (i.e.,~the subordinate Brownian motion is not a compound Poisson process). 
The weak lower scaling condition at infinity
governs the short-time small-space behavior of the subordinate Brownian motion. We also need a weak condition on the behavior of $\Phi$ near zero. We assume that
\begin{equation}\label{e:intro-psibound}
\int_0^1 \frac{\Phi(r)}{r}\, dr=C_* <\infty\, .
\end{equation}

The following function will play a prominent role in the paper. For $t>0$ 
and $x\in \R^d$ we define
\begin{equation}\label{e:intro-rho-def}
\rho(t,x)=\rho^{(d)}(t,x):={\Phi\left(\left(\frac{1}{\Phi^{-1}(t^{-1})}+|x|\right)^{-1}\right)}\left(\frac{1}{\Phi^{-1}(t^{-1})}+|x|\right)^{-d}\, .
\end{equation}
In case when $\Phi(r)=r^{\alpha}$ we see that $\rho(t,x)=(t^{1/ \alpha}+|x|)^{-d-\alpha}$. It is well known that $t(t^{1/ \alpha}+|x|)^{-d-\alpha}$ is comparable to the heat kernel $p(t,x)$ of the isotropic $\alpha$-stable process in $\R^d$. 
We will prove later in this paper  (see Proposition \ref{p:upperestonp}) that $t\rho(t,x)$ is an upper bound of the heat kernel of the subordinate Brownian motion with characteristic exponent $\Phi$.

We assume that $J: \R^d \to (0, \infty)$ is symmetric 
in the sense that $J(x)=J(-x)$ for all $x\in \R^d$
and there exists $\gamma_0>0$ such that 
\begin{equation}\label{e:psi1}
\gamma^{-1}_0 j(|y|)\le J(y) \le \gamma_0 j(|y|), \quad \mbox{for all } y\in \R^d\, .
\end{equation}
Following \eqref{e:intro-stable-operator}, we define a non-symmetric and non-local operator
\begin{equation}\label{e:intro-operator}
\LL^{\kappa}f(x)=\LL^{\kappa, 0}f(x):= \mathrm{p.v.}\int_{\R^d}
(f(x+z)-f(x))\kappa(x,z)J(z)\, dz
:=\lim_{\eps \downarrow 0} \LL^{\kappa, \eps}f(x)\, ,
\end{equation}
where $$
\LL^{\kappa, \eps}f(x):=\int_{|z| > \eps}
(f(x+z)-f(x))\kappa(x,z)J(z)\, dz, \quad \eps>0.
$$

The following theorem is the main result of this paper.

\begin{thm}\label{t:intro-main}
Assume that $\Phi$ satisfies \eqref{e:intro-wsc} and \eqref{e:intro-psibound}, 
that $J$ satisfies \eqref{e:psi1}, 
and that $\kappa$  satisfies \eqref{e:intro-kappa} and \eqref{e:intro-kappa-holder}.
Suppose there exists a function $g:\R^d \to (0, \infty)$ such that
\begin{align}
\label{e:nshf}
\lim_{x \to \infty}g(x)=\infty \quad \text{ and } \quad \LL^{\kappa}g(x)/g(x) \text{ is bounded from above.}
\end{align}
Then there exists a unique non-negative jointly continuous function $p^{\kappa}(t,x,y)$ 
on $(0,\infty)\times \R^d\times \R^d$ solving
\begin{equation}\label{e:intro-main-1}
\partial_t p^{\kappa}(t,x,y)=\LL^{\kappa}p^{\kappa}(t,\cdot, y)(x)\, , \quad x\neq y\, ,
\end{equation}
and satisfying the following properties:

\noindent
(i) (Upper bound) For every $T\ge 1$, there is a constant $c_1>0$ so that for all $t\in (0,T]$ and $x,y\in \R^d$,
\begin{equation}\label{e:intro-main-2}
p^{\kappa}(t,x,y)\le c_1 t \rho(t,x-y)\, .
\end{equation}

\noindent (ii) (Fractional derivative estimate) For any $x,y\in \R^d$, $x\neq y$, the map $t\mapsto \LL^{\kappa}p^{\kappa}(t,\cdot,y)(x)$ is continuous in $(0, \infty)$, 
and, for each $T\ge 1$ 
 there is a constant $c_2>0$ so that for all $t\in (0,T]$, $\eps \in [0,1]$ and $x,y\in \R^d$,
\begin{equation}\label{e:intro-main-4}
| \LL^{\kappa, \eps}p^{\kappa}(t,\cdot,y)(x)|\le c_2 \rho(t,x-y)\, .
\end{equation}

\noindent
(iii) (Continuity) For any bounded and uniformly continuous function $f:\R^d\to \R$,
\begin{equation}\label{e:intro-main-5}
\lim_{t\downarrow 0}\sup_{x\in \R^d}\left| \int_{\R^d}p^{\kappa}(t,x,y)f(y)\, dy-f(x)\right|=0\, .
\end{equation}
Moreover, the constants $c_1$ and $c_2$ can be chosen so that they depend only on 
$T$,  $\Phi^{-1}(T^{-1})$,  $d$, $a_1$, $\delta_1$, $C_*$, $\beta$, $\gamma_0$,  $\kappa_0$, $\kappa_1$ and $\kappa_2$. 
\end{thm}

The assumption \eqref{e:nshf} is a quite mild one. For example, if $\int_{|z|>1} |z|^\eps j(|z|)dz < \infty$ for some $\eps>0$, then 
\eqref{e:nshf} holds, see Remark \ref{r:sufjeps} below.

Some further properties of the heat kernel $p^{\kappa}(t,x,y)$ are listed in the following result.
\begin{thm}\label{t:intro-further-properties}
Suppose that the assumptions of Theorem \ref{t:intro-main} are satisfied. 

\noindent
(1) (Conservativeness) For all 
$(t,x)\in (0, \infty) \times \R^d$, 
\begin{equation}\label{e:intro-main-6}
\int_{\R^d}p^{\kappa}(t,x,y)\, dy =1\, .
\end{equation}

\noindent
(2) (Chapman-Kolmogorov equation) For all $s,t > 0$ and all $x,y\in \R^d$,
\begin{equation}\label{e:intro-main-7}
\int_{\R^d}p^{\kappa}(t,x,z)p^{\kappa}(s,z,y)\, dz =p^{\kappa}(t+s,x,y)\, .
\end{equation}

\noindent
(3) 
(Joint H\"{o}lder continuity) For every $T\ge 1$ and $\gamma\in (0,\delta_1) \cap (0,1]$, there is a constant 
$c_3=c_3(T, d,\delta_1, a_1, \beta, C_*, \Phi^{-1}(T^{-1}), \gamma_0, \kappa_0, \kappa_1, \kappa_2)>0$
 such that for all $0<s\le t\le T$ and $x,x',y\in \R^d$,
\begin{equation}
\label{e:intro-main-8}
|p^{\kappa}(s,x,y)-p^{\kappa}(t,x',y)|
\le c_3 \left( |t-s|+|x-x'|^\gamma t\, \Phi^{-1}(t^{-1}) \right) 
(\rho(s, x-y)\vee\rho(s,  x'-y))\, .
\end{equation}
Furthermore, if the constant $\delta_1$ in \eqref{e:intro-wsc} belongs to $(2/3,2)$ and
the constant $\beta$ in \eqref{e:intro-kappa-holder} satisfies 
 $\beta+\delta_1 >1$  then \eqref{e:intro-main-8} holds with 
$\gamma=1$.

\noindent
(4)
(Gradient  estimate) 
If $\delta_1 \in (2/3,2)$,
and $\beta+\delta_1 >1$, 
then for every $T\ge 1$, there exists
$c_4=c_4(T, d,\delta_1, a_1, \beta, C_*, \Phi^{-1}(T^{-1}), \gamma_0, \kappa_0, \kappa_1, \kappa_2)>0$
so that for all $x,y\in \R^d$, $x\neq y$, and $t\in (0,T]$,
\begin{equation}\label{e:intro-main-9}
|\nabla_x p^{\kappa}(t,x,y)| \le 
c_4 \Phi^{-1}(t^{-1}) t \rho(t, |x-y|)\, .
\end{equation}
\end{thm}

Note that the gradient estimate \eqref{e:intro-main-9} is an improvement of the corresponding estimate \cite[(4.19)]{CZ} 
in the sense that
the parameter $\delta_1$ could be smaller than one as long as it is still larger than $2/3$ and $\beta+\delta_1>1$.

For $t>0$, define the operator $P_t^{\kappa}$ by
\begin{equation}\label{e:intro-semigroup}
P_t^{\kappa}f(x)=\int_{\R^d} p^{\kappa}(t,x,y)f(y)\, dy\, ,\quad x\in \R^d\, ,
\end{equation}
where $f$ is a non-negative (or bounded) Borel function on $\R^d$, and let $P_0^{\kappa}=\mathrm{Id}$. Then
by Theorems \ref{t:intro-main} and  \ref{t:intro-further-properties},  $(P_t^{\kappa})_{t\ge 0}$ is a Feller semigroup with the strong Feller property.  
Let $C_b^{2, \eps} (\R^d)$ be  the space of bounded twice differentiable functions in $\R^d$ whose second derivatives are uniformly H\"older continuous.
We further have  

\begin{thm}\label{t:intro-semigroup}
Suppose that the assumptions of Theorem \ref{t:intro-main} are satisfied.

\noindent
(1) (Generator) Let $\eps >0$. For any $f\in C_b^{2, \eps} (\R^d)$,  we have
\begin{equation}\label{e:intro-main-10}
\lim_{t\downarrow 0}\frac{1}{t}\left(P_t^{\kappa}f(x)-f(x)\right)=\LL^{\kappa}f(x)\, ,
\end{equation}
and the convergence is uniform.

\noindent
(2) (Analyticity) The 
semigroup $(P_t^{\kappa})_{t\ge 0}$ of $\LL^{\kappa}$ is analytic in $L^p(\R^d)$ for every $p\in [1,\infty)$.
\end{thm}

Finally, 
under an additional assumption, we prove by probabilistic methods
a lower bound for the heat kernel $p^{\kappa}(t,x,y)$. The \emph{weak upper scaling condition} 
 means that there 
 exist $\delta_2\in (0,2)$ and $a_2>0$ such that
\begin{equation}\label{e:intro-wusc}
\Phi(\lambda r)\le a_2 \lambda^{\delta_2}\Phi(r)\, ,\quad \lambda\ge 1, r\ge 1\, .
\end{equation}

\begin{thm} \label{t:intro-lower-bound}
Suppose  that $\Phi$ satisfies \eqref{e:intro-wsc}, \eqref{e:intro-wusc} and \eqref{e:intro-psibound}, 
that $J$ satisfies \eqref{e:psi1}, 
and that $\kappa$ satisfies \eqref{e:intro-kappa} and \eqref{e:intro-kappa-holder}.
Suppose also that there exists a function $g:\R^d \to (0, \infty)$ such that
\eqref{e:nshf} holds.
For every $T\ge 1$, there exists 
$c_5=c_5(T, d,\delta_1,\delta_2, \gamma_0, C_*, \Phi^{-1}(T^{-1}), a_1, a_2$, $ \beta, \kappa_0, \kappa_1, \kappa_2)>0$ 
such that for all $t\in (0,T]$,
\begin{equation}\label{e:intro-main-11}
p^{\kappa}(t,x,y)\ge c_5
\begin{cases}
\Phi^{-1}(t^{-1})^d & \text{ if } |x-y| \le 3\Phi^{-1}(t^{-1})^{-1},\\
 t j\left( |x-y|\right)  & \text{ if } |x-y|  >  3\Phi^{-1}(t^{-1})^{-1}.
\end{cases}
\end{equation}
In particular, for all $T, M \ge 1$, there exists 
$c_6=c_6(T, d,\delta_1,\delta_2, \gamma_0, C_*, \Phi^{-1}(T^{-1}), a_1, a_2, \beta, \kappa_0,$ $\kappa_1, \kappa_2)>0$
 for all $t\in (0,T]$ and $x,y\in \R^d$ with $|x-y| \le M$,
\begin{equation}\label{e:intro-main-111}
p^{\kappa}(t,x,y)\ge c_6 t \rho(t,x-y)\, .
\end{equation}
\end{thm}

Theorems \ref{t:intro-main}-\ref{t:intro-lower-bound} generalize \cite[Theorem 1.1]{CZ}. Note that the lower bound \eqref{e:intro-main-111} of $p^{\kappa}(t,x,y)$ is stated only for 
$|x-y|\le M$. 
This is natural in view of the fact that \eqref{e:intro-wsc} and \eqref{e:intro-wusc} only give information about short-time small-space behavior of the underlying subordinate Brownian motion.
We remark in passing that, the upper bound \eqref{e:intro-main-2} may not be sharp
under the assumptions \eqref{e:intro-wsc} and \eqref{e:intro-psibound}. When $\Phi$
satisfies scaling conditions both near infinity and near the origin, see \cite[(H1) and (H2)]{KSV}, the upper bound \eqref{e:intro-main-2} is sharp in the sense that the lower bound \eqref{e:intro-main-111}
is valid for all $x,y\in \R^d$.

The assumptions  \eqref{e:intro-wsc}, \eqref{e:intro-psibound}, \eqref{e:nshf} and
\eqref{e:intro-wusc} are very weak conditions and they are satisfied by many subordinate Brownian motions. For the reader's convenience, we list some examples 
of $\phi$, besides the Laplace exponent of the stable subordinator, such that 
$\Phi(r)=\phi(r^2)$ satisfies these assumptions.

\begin{description}
\item{(1)} $\phi(\lambda)=\lambda^{\alpha_1} + \lambda^{\alpha_2}$, $0<\alpha_1<\alpha_2<1$;
\item{(2)} $\phi(\lambda)=(\lambda+\lambda^{\alpha_1})^{\alpha_2}$, $\alpha_1, \alpha_2\in (0, 1)$;
\item{(3)} $\phi(\lambda)=(\lambda+m^{1/\alpha})^{\alpha}-m$, $\alpha\in (0, 1)$,
$m>0$;
\item{(4)} $\phi(\lambda)=\lambda^{\alpha_1}(\log(1+\lambda))^{\alpha_2}$, $\alpha_1\in (0, 1)$,
$\alpha_2\in (0, 1-\alpha_1]$;
\item{(5)} $\phi(\lambda)=\lambda^{\alpha_1}(\log(1+\lambda))^{-\alpha_2}$, $\alpha_1\in (0, 1)$,
$\alpha_2\in (0, \alpha_1)$;
\item{(6)} $\phi(\lambda)=\lambda/\log(1+\lambda^{\alpha})$, $\alpha\in (0, 1)$.
\end{description}

The functions in (1)--(5) satisfy \eqref{e:intro-wsc}, \eqref{e:intro-psibound},
\eqref{e:intro-wusc}  and \eqref{e:nshf} (see \eqref{e:jupper} and Remark \ref{r:sufjeps}); while the function in (6) satisfies \eqref{e:intro-wsc}, \eqref{e:intro-psibound} and \eqref{e:nshf}, but does not satisfy \eqref{e:intro-wusc}. The function
$\phi(\lambda)=\lambda/\log(1+\lambda)$ satisfies \eqref{e:intro-wsc}, but does not satisfy the other two conditions.

In order to prove our main results, 
we follow the ideas and the road-map from \cite{CZ}. At many stages we encounter substantial technical difficulties due to the fact that in the stable-like case one deals with power functions while in the present situation the power functions are replaced with a quite general $\Phi$ and its variants. We also strive to simplify the proofs and streamline the presentation. In some places we provide full proofs where in \cite{CZ} only an indication is given. 
On the other hand, we skip some proofs which would be almost identical to the corresponding ones in \cite{CZ}.
Below is a detailed outline of the paper with emphasis on 
the
main differences from \cite{CZ}.

In Section 2 we start by introducing the basic setup, state again the assumptions, 
and
derive some of the consequences.
In Subsection 2.1 we discuss convolution inequalities, cf.~Lemma \ref{l:convolution}. While in \cite{CZ} these involve power functions, the most challenging task in the present setting was to find appropriate versions of these inequalities. 
The main new technical result here is Lemma \ref{l:convolution}.

In Section 3 we first study the heat kernel $p(t,x)$ of a 
symmetric L\'evy process $Z$ with 
L\'evy density $j_Z$ comparable to the L\'evy density $j$ 
of the subordinate Brownian motion with characteristic exponent $\Phi$. We  prove 
the joint Lipschitz continuity  of $p(t,x)$
and then, based on a result from \cite{KS},  that $t\rho(t,x)$ is the upper bound of 
$p(t,x)$ for all $x\in \R^d$ and small $t$, 
cf.~Proposition \ref{p:upperestonp}. In Subsection 3.1, 
we provide  
some useful estimates on functions of $p(t,x)$.
In Subsection 3.2, we specify $j_Z$ by assuming $j_Z(z)=\mathfrak{K}(z)J(z)$, with
$\mathfrak{K}$ being symmetric and bounded between two positive constants.
Let $\LL^{\mathfrak{K}}$ be the infinitesimal generator of 
the corresponding process
and let $p^{\mathfrak{K}}$ be its heat kernel.
We look at the continuous dependence of $p^{\mathfrak{K}}$
with respect to $\mathfrak{K}$. This subsection follows the ideas and proofs from \cite{CZ} with additional technical difficulties.

Given a function $\kappa$ satisfying  \eqref{e:intro-kappa} and \eqref{e:intro-kappa-holder}, 
we define, for a fixed $y\in \R^d$, 
$\mathfrak{K}_y=\kappa(y,\cdot )$ and denote by $p_y(t,x)$ the heat kernel of the freezing operator  $\LL^{\mathfrak{K}_y}$. Various estimates and joint continuity of $p_y(t,x)$ are shown in Subsection 4.1. The rest of Section 4 is devoted to constructing the heat kernel $p^{\kappa}(t,x,y)$ of the operator $\LL^{\kappa}$. The heat kernel should have the form
\begin{equation}\label{e:intro-p-kappa}
p^{\kappa}(t,x,y)=p_y(t,x-y)+\int_0^t \int_{\R^d} p_z(t-s,x-z)q(s,z,y)\, dz\, ds\, ,
\end{equation}
where according to Levi's method the function $q(t,x,y)$ solves the integral equation
\begin{equation}\label{e:intro-q}
q(t,x,y)=q_0(t,x,y)+\int_0^t \int_{\R^d} q_0(t-s,x-z)q(s,z,y)\, dz\, ds\, ,
\end{equation}
with $q_0(t,x,y)=(\LL^{\mathfrak{K}_x}-\LL^{\mathfrak{K}_y})p_y(t,x-y)$. The main result is Theorem \ref{t:definition-of-q} showing existence and joint continuity of $q(t,x,y)$ satisfying \eqref{e:intro-q}. We follow \cite[Theorem 3.1]{CZ}, and give a full proof.
Joint continuity and various estimates of $p^{\kappa}(t,x,y)$ defined by $\eqref{e:intro-p-kappa}$ are given in Subsection 4.3.

Section 5 contains proofs of Theorems \ref{t:intro-main}--\ref{t:intro-lower-bound}. 
We start with a version of a non-local maximum principle 
in Theorem \ref{t:nonlocal-max-principle} which is somewhat different from the one in \cite[Theorem 4.1]{CZ}, 
continue with two results about the semigroup $(P_t^{\kappa})_{t\ge 0}$ and then complete the proofs.

In this paper, we use the following notations.
We will use ``$:=$" to denote a definition, which is  read as ``is defined to be".
For any two positive functions $f$ and $g$,
$f\asymp g$ means that there is a positive constant $c\geq 1$
so that $c^{-1}\, g \leq f \leq c\, g$ on their common domain of
definition.
 For a set $W$ in $\R^d$, $|W|$ denotes the Lebesgue measure of $W$ in $\R^d$. 
 For  a function space $\H(U)$
on an open set $U$ in $\R^d$, we let   
$\H_c(U):=\{f\in\H(U): f \mbox{ has  compact support}\},$
$\H_0(U):=\{f\in\H(U): f \mbox{ vanishes at infinity}\}$ and $\H_b(U):=\{f\in\H(U): f \mbox{ is bounded}\}$.

 Throughout the rest of this paper, the positive constants
 $\delta_1,\delta_2, \gamma_0, a_1, a_2, \beta, \kappa_0, \kappa_1, \kappa_2, C_i$,
$i=0,1,2,\dots $, can be regarded as fixed.
In the statements of results and the proofs, the constants $c_i=c_i(a,b,c,\ldots)$, $i=0,1,2,  \dots$, denote generic constants depending on $a, b, c, \ldots$, whose exact values are unimportant.
They start
anew in each statement and each proof.
The dependence of the constants on the dimension $d \ge 1$, $C_*$, $\Phi^{-1}((2T)^{-1})$, $\Phi^{-1}(T^{-1})$   and $\gamma_0$
may not be mentioned explicitly.

\section{Preliminaries} 

It is well known that the Laplace exponent $\phi$  of a subordinator is a Bernstein function and 
\begin{equation}\label{e:Berall}
\phi(\lambda t)\le \lambda\phi(t) \qquad \text{ for all } \lambda \ge 1, t >0\, .
\end{equation}
For notational convenience, in this paper, we denote  $\Phi(r)=\phi(r^2)$ and without loss of generality we assume that 
$\Phi(1)=1$.

Throughout this paper $\phi$ is the Laplace exponent of a subordinator and 
$\Phi(r)=\phi(r^2)$ 
 satisfies the weak lower scaling 
 condition \eqref{e:intro-wsc} at infinity. 
 This can be reformulated as follows: 
There exist $\delta_1 \in (0,2]$ and a positive constant 
$a_1\in (0, 1]$
 such that  for any $r_0\in (0, 1]$, 
\begin{equation}\label{e:lower-scaling}
a_1\lambda^{\delta_1}r_0  ^{\delta_1}\Phi(r)\le \Phi(\lambda r)\, ,\quad \lambda \ge 1, r\ge r_0\, .
\end{equation}
In fact, suppose $r_0 \le  r <1$ and $\lambda \ge 1$. Then, 
  $ \Phi(\lambda r)   \ge a_1\lambda^{\delta_1}r_0  ^{\delta_1}    \Phi(1) \ge a_1\lambda^{\delta_1}r_0  ^{\delta_1}\Phi(r)$    if $\lambda r > 1$, 
  and  $\Phi(\lambda r)   \ge \Phi(r)  \ge a_1\lambda^{\delta_1}r_0  ^{\delta_1}\Phi(r)$  if $\lambda r \le 1$. 
  
Since $\phi$ is a Bernstein function and we assume \eqref{e:lower-scaling}, it follows that 
$\Phi$ is strictly increasing and $\lim_{\lambda\to \infty}\Phi(\lambda)=\infty$.  We denote by $\Phi^{-1} : (0, \infty) \to (0, \infty)$ the inverse function of $\Phi$.

From \eqref{e:Berall} we have 
\begin{equation}\label{e:psi-inverse-sbm}
\Phi^{-1}(\lambda r)\ge \lambda^{1/2}\Phi^{-1}(r)\, ,\quad \lambda \ge 1, r>0\, .
\end{equation}
Moreover, by \eqref{e:lower-scaling}, $\Phi^{-1}$ satisfies the following weak upper scaling condition at infinity: 
For any $r_0\in (0, 1]$, 
\begin{equation}\label{e:lsc-inverse}
\Phi^{-1}(\lambda r)\le a_1^{-1/\delta_1} \Phi^{-1}(r_0)^{-1} \lambda^{1/\delta_1}\Phi^{-1}(r)\, ,
\quad \lambda \ge 1, r\ge r_0\, .
\end{equation}
In fact, from  \eqref{e:lower-scaling} we get $\Phi^{-1}(\lambda r)\le a_1^{-1/\delta_1} \lambda^{1/\delta_1}\Phi^{-1}(r)$ for $ \lambda \ge 1$ and $r\ge 1$. 
Suppose $r_0 \le  r <1$. Then,  $\Phi^{-1}(\lambda r)  \le 1 \le a_1^{-1/\delta_1} \Phi^{-1}(r_0)^{-1} \lambda^{1/\delta_1}\Phi^{-1}(r)$  if $\lambda r \le 1$, and
  $\Phi^{-1}(\lambda r)  \le a_1^{-1/\delta_1}  \lambda^{1/\delta_1} r^{1/\delta_1} 
      \le a_1^{-1/\delta_1}  \lambda^{1/\delta_1}\Phi^{-1}(r_0)^{-1} \Phi^{-1}(r)$  if $\lambda r > 1$.

For $t>0$ and $x\in \R^d$,  we define
functions $r(t,x)$ and $\rho(t,x)$ by
\begin{align*}
r(t,x)=\Phi^{-1}(t^{-1})^d \wedge \frac{t\Phi(|x|^{-1})}{|x|^d}
\end{align*}
and
\begin{equation}\label{e:rho-def}
\rho(t,x)=\rho^{(d)}(t,x):={\Phi\left(\left(\frac{1}{\Phi^{-1}(t^{-1})}+|x|\right)^{-1}\right)}\left(\frac{1}{\Phi^{-1}(t^{-1})}+|x|\right)^{-d}\, .
\end{equation}
Note that, by \cite[Lemma 17]{BGR}, 
\begin{equation}\label{e:relphix}
t\Phi(|x|^{-1})|x|^{-d}\ge \Phi^{-1}(t^{-1})^d \quad \mbox{if and only if  } \quad  t\Phi(|x|^{-1})\ge 1.
\end{equation}
\begin{prop}\label{p:p-q}
 For all $t>0$ and $x\in \R^d$, $t \rho(t,x)\le  r(t,x) \le 2^{d+2} t \rho(t,x)$. 
\end{prop}
\pf {\it Case 1:} $t\Phi(|x|^{-1})\ge 1$. In this case, by \eqref{e:relphix} we have that $r(t,x)= \Phi^{-1}(t^{-1})^d $. 
Since $|x|\le \frac{1}{\Phi^{-1}(t^{-1})}$, we have
\begin{equation}\label{e:0}
\frac{1}{\Phi^{-1}(t^{-1})}\le \frac{1}{\Phi^{-1}(t^{-1})}+|x|\le \frac{2}{\Phi^{-1}(t^{-1})}\, .
\end{equation}
This and \eqref{e:Berall} imply that
$$
t^{-1}=\Phi(\Phi^{-1}(t^{-1}))\ge \Phi\left(\left(\frac{1}{\Phi^{-1}(t^{-1})}+|x|\right)^{-1}\right) \ge \Phi(2^{-1}\Phi^{-1}(t^{-1})) \ge \frac14\Phi(\Phi^{-1}(t^{-1})) =\frac14t^{-1}
$$
and
$$
2^{-d}\Phi^{-1}(t^{-1})^{-d}\le \left(\frac{1}{\Phi^{-1}(t^{-1})}+|x|\right)^{-d}\le \Phi^{-1}(t^{-1})^{-d}\, .
$$
The last two displays imply that $2^{-d-2}\Phi^{-1}(t^{-1})^d   \le t \rho(t,x) \le \Phi^{-1}(t^{-1})^d  $.

\noindent
{\it Case 2:} $t\Phi(|x|^{-1})\le 1$. In this case, by \eqref{e:relphix} we have that $r(t,x)= \frac{t\Phi(|x|^{-1})}{|x|^d}$. Since $|x|\ge \frac{1}{\Phi^{-1}(t^{-1})}$, we have
$$
|x|^{-1}\ge \left(\frac{1}{\Phi^{-1}(t^{-1})}+|x|\right)^{-1} \ge 2^{-1}|x|^{-1}\, .
$$
This with \eqref{e:Berall} implies that
$$
\Phi(|x|^{-1})\ge \Phi\left(\left(\frac{1}{\Phi^{-1}(t^{-1})}+|x|\right)^{-1} \right)
\ge \Phi(2^{-1}|x|^{-1})\ge \frac14 \Phi(|x|^{-1})\, .
$$
The last two displays imply the conclusion of the proposition in Case 2. \qed

\begin{lemma}\label{l:psi-and-f}
Let $T\ge 1$ and $c=(2(2/a_1)^{1/\delta_1} /\Phi^{-1}((2T)^{-1}))^{d+2}$.

\noindent
(a)
For all $0<s<t\le T$ and $x,z\in \R^d$,
\begin{equation}\label{e:rho}
\rho(t-s,x-z)\rho(s,z) \le c\big(\rho(t-s,x-z)+\rho(s,z)\big)\rho(t,x)\, .
\end{equation}

\noindent 
(b) For every $x\in \R^d$ and $0<t/2\le s \le t\le T$, $ \rho(t,x) \le \rho(s,x)\le  2c  \rho(t,x)$.
\end{lemma}
\pf 
(a)
By  \eqref{e:lsc-inverse} we have that  for all $0<t,s \le T$,
\begin{equation}\label{e:psi-inv-sub}
\frac{1}{\Phi^{-1}((t+s)^{-1})}\le
\frac{1}{\Phi^{-1}(2^{-1}(t \vee s)^{-1})}
\le
 c_1\left(\frac{1}{\Phi^{-1}(t^{-1})}+\frac{1}{\Phi^{-1}(s^{-1})}\right)\, ,
\end{equation}
where  
$c_1=(2/a_1)^{1/\delta_1} /\Phi^{-1}((2T)^{-1})
\ge 1$.

Define $\varrho:(0,\infty)\to (0,\infty)$ by 
$
\varrho(r):={r^d}/{\Phi(r^{-1})},
$
so that 
$
\rho(t,x)=(\varrho(\frac{1}{\Phi^{-1}(t^{-1})}+|x|))^{-1}.
$
 For all $a,b>0$, $(a+b)^d\le 2^d (a \vee b)^d$ and, by \eqref{e:Berall}, 
$\Phi((a+b)^{-1})\ge \Phi(2^{-1}(a\vee b)^{-1})\ge 4^{-1} \Phi((a \vee b)^{-1})$. Therefore,
for all $a,b>0$,
\begin{equation}\label{e:f-sub}
\varrho(a+b)\le 2^{d+2}\varrho(a \vee b)\le 2^{d+2}(\varrho(a)+\varrho(b))\, .
\end{equation}
Moreover, 
  \eqref{e:Berall} implies that for $r>0$,
\begin{equation}\label{e:f-doubling}
\varrho(c_1 r)=\frac{(c_1 r)^d}{\Phi(c_1^{-1}r^{-1})}\le c_1^{d+2} \frac{r^d}{\Phi(r^{-1})}= c_1^{d+2} \varrho(r)\, .
\end{equation}

By using \eqref{e:psi-inv-sub}--\eqref{e:f-doubling}, we have 
\begin{align}
&\varrho\left(\frac{1}{\Phi^{-1}(t^{-1})}+|x|\right)\le  \varrho\left(c_1\left(\Big(\frac{1}{\Phi^{-1}((t-s)^{-1})}+|x-z|\Big)+\Big(\frac{1}{\Phi^{-1}(s^{-1})}+|z|\Big)\right)\right)\nn\\
\le & c_1^{d+2}  \varrho\left(\Big(\frac{1}{\Phi^{-1}((t-s)^{-1})}+|x-z|\Big)+\Big(\frac{1}{\Phi^{-1}(s^{-1})}+|z|\Big)\right)\nn\\
\le & (2c_1)^{d+2} \left(\varrho\left(\frac{1}{\Phi^{-1}((t-s)^{-1})}+|x-z|\right)+\varrho\left(\frac{1}{\Phi^{-1}(s^{-1})}+|z|\right)\right)\, .\label{e:psi-and-f}
\end{align}
Thus we have that 
for $0<s<t\le T$ and $x,z\in \R^d$,
\begin{eqnarray*}
\lefteqn{\left(\rho(t-s,x-z)+\rho(s,z)\right)\rho(t,x)}\\
&=&\frac{\varrho\left(\frac{1}{\Phi^{-1}((t-s)^{-1})}+|x-z|\right)+\varrho\left(\frac{1}{\Phi^{-1}(s^{-1})}+|z|\right)}{\varrho\left(\frac{1}{\Phi^{-1}((t-s)^{-1})}+|x-z|\right)\varrho\left(\frac{1}{\Phi^{-1}(s^{-1})}+|z|\right)}\frac{1}{\varrho\left(\frac{1}{\Phi^{-1}(t^{-1})}+|x|\right)}\\
&\ge &(2c_1)^{-d-2}\frac{1}{\varrho\left(\frac{1}{\Phi^{-1}((t-s)^{-1})}+|x-z|\right)\varrho\left(\frac{1}{\Phi^{-1}(s^{-1})}+|z|\right)}\\
&=&(2c_1)^{-d-2}\rho(t-s,x-z)\rho(s,z)\, .
\end{eqnarray*}

\noindent
(b) 
This
follows from \eqref{e:psi-and-f} by taking $s=t/2$, $z=0$ and by using that $\varrho$ is increasing.
\qed

\subsection{Convolution inequalities}
Let $B(a,b)$ be the beta function, i.e., $B(a,b)=\int_0^1 s^{a-1} (1-s)^{b-1}ds$, $a,b>0$.

\begin{lemma}\label{l:convoluton-inequality}
Let $\beta, \gamma,  \eta, \theta\in \R$ be such that 
${\bf 1}_{\beta\ge 0}(\beta/2)  +{\bf 1}_{\beta< 0}(\beta/\delta_1)   +1-\theta>0$ and ${\bf 1}_{\gamma\ge 0}(\gamma/2)  +{\bf 1}_{\gamma< 0}(\gamma/\delta_1)  +1-\eta>0$. Then for 
every $t >0$,  we have
\begin{equation}\label{e:convolution-inequality}
\int_0^t u^{-\eta} \Phi^{-1}(u^{-1})^{-\gamma} (t-u)^{-\theta}\Phi^{-1}((t-u)^{-1})^{-\beta}\, du \le C t^{1-\eta-\theta}\Phi^{-1}(t^{-1})^{-\gamma-\beta}\,.
\end{equation}
Moreover, if $\beta\ge 0$ and $\gamma\ge 0$ then \eqref{e:convolution-inequality} holds for all $t>0$ with 
$C=B(\beta/2+1-\theta, \gamma/2+1-\eta)$.
\end{lemma}
\pf Let $I$ denote the integral in \eqref{e:convolution-inequality}. By the change of variables $s=u/t$ we get that 
$$
I=t^{1-\eta-\theta}\int_0^1 s^{-\eta}\Phi^{-1}(t^{-1}s^{-1})^{-\gamma}(1-s)^{-\theta}\Phi^{-1}(t^{-1}(1-s)^{-1})^{-\beta}\, ds\, .
$$
Since $s^{-1}\ge 1$ and $(1-s)^{-1}\ge 1$, we have by \eqref{e:psi-inverse-sbm} that $\Phi^{-1}(t^{-1}s^{-1})\ge  s^{-1/2}\Phi^{-1}(t^{-1})$ and $\Phi^{-1}(t^{-1}(1-s)^{-1})\ge  (1-s)^{-1/2}\Phi^{-1}(t^{-1})$. 
Moreover, when $t \in (0, T]$, by \eqref{e:lsc-inverse} we have 
$$
\Phi^{-1}(t^{-1}s^{-1})\le  a_1^{-1/\delta_1} \Phi^{-1}(T^{-1})^{-1} s^{-1/\delta_1}\Phi^{-1}(t^{-1})
$$ 
and 
$$
\Phi^{-1}(t^{-1}(1-s)^{-1})\le a_1^{-1/\delta_1}
\Phi^{-1}(T^{-1})^{-1} (1-s)^{-1/\delta_1}\Phi^{-1}(t^{-1}).
$$
Hence,
\begin{eqnarray*}
I &\le & c_1 t^{1-\eta-\theta}\Phi^{-1}(t^{-1})^{-\gamma-\beta}\int_0^1  s^{{\bf 1}_{\gamma\ge 0}(\gamma/2)  +{\bf 1}_{\gamma< 0}(\gamma/\delta_1) -\eta}(1-s)^{{\bf 1}_{\beta\ge 0}(\beta/2)  +{\bf 1}_{\beta< 0}(\beta/\delta_1) -\theta}\, ds \\
& = &  C \Phi^{-1}(t^{-1})^{-\gamma-\beta}\, . 
\end{eqnarray*}
When $\beta\ge 0$ and $\gamma\ge 0$ then the above inequality holds for all $t>0$ with $c_1=1$ so 
$C=B(\beta/2+1-\theta, \gamma/2+1-\eta)$.
\qed

\begin{lemma}\label{l:rsconvolution-inequality}
Suppose that $0<t_1\le t_2<\infty$. Under
the assumptions of Lemma \ref{l:convoluton-inequality},
we have
$$
\lim_{h\to0}\sup_{t\in [t_1, t_2]}\left(\int^h_0+\int^t_{t-h} \right) u^{-\eta} \Phi^{-1}(u^{-1})^{-\gamma} (t-u)^{-\theta}\Phi^{-1}((t-u)^{-1})^{-\beta}\, du=0.
$$
\end{lemma}

\pf Under the assumptions of this lemma, by repeating the argument in the proof of Lemma \ref{l:convoluton-inequality}, we have that for all 
$t \in [t_1, t_2]$,
\begin{eqnarray*}
&&\left(\int^h_0+\int^t_{t-h} \right) u^{-\eta} \Phi^{-1}(u^{-1})^{-\gamma} (t-u)^{-\theta}\Phi^{-1}((t-u)^{-1})^{-\beta}\, du\\
&&\le \ \left( t_1^{1-\eta-\theta}\vee  t_2^{1-\eta-\theta}\right)\left(\Phi^{-1}(t_1^{-1})^{-\gamma-\beta}\vee \Phi^{-1}(t_2^{-1})^{-\gamma-\beta}\right)\\
&&\ \ \ \times\left(\int^{h/t_1}_0+\int^1_{1-h/t_1} \right)s^{{\bf 1}_{\gamma\ge 0}(\gamma/2)  +{\bf 1}_{\gamma< 0}(\gamma/\delta_1) -\eta}(1-s)^{{\bf 1}_{\beta\ge 0}(\beta/2)  +{\bf 1}_{\beta< 0}(\beta/\delta_1) -\theta}\, ds.
\end{eqnarray*}
Now the conclusion of the lemma follows immediately.
\qed

For $\gamma, \beta \in \R$, 
we define
$$
\rho_{\gamma}^{\beta}(t,x):=\Phi^{-1}(t^{-1})^{-\gamma}(|x|^{\beta}\wedge 1) \rho(t,x)\, ,\quad t>0, x\in \R^d\, .
$$
Note that 
 $\rho_0^0(t,x)=\rho(t,x)$.

\begin{remark}\label{r:rho-decreasing}
{\rm
Recall that   $\Phi$ is increasing. 
Thus it is straightforward to see that the following inequalities are true: for $T \ge 1$,
\begin{eqnarray}
 \rho_{\gamma_1}^{\beta}(t,x) \le \Phi^{-1}(T^{-1})^{\gamma_2-\gamma_1} \rho_{\gamma_2}^{\beta}(t,x),& & \quad (t,x)\in (0,T]\times \R^d, \quad  \gamma_2\le \gamma_1\, ,\label{e:nonincrease-gamma}\\
\rho_{\gamma}^{\beta_1}(t,x) \le \rho_{\gamma}^{\beta_2}(t,x),\qquad\qquad\qquad& & \quad (t,x)\in (0,\infty)\times \R^d,  \quad 0 \le \beta_2\le \beta_1 \label{e:nonincrease-beta}\, .
\end{eqnarray}
}
\end{remark}

We record the following inequality: for every $T \ge 1$, $t \in (0, T]$  and $\beta<\delta_1$,
\begin{align}\label{e:integral-psi}
\int_{\Phi^{-1}(T^{-1})/\Phi^{-1}(t^{-1})}^{1}r^{\beta-1}\Phi(r^{-1})dr 
&\le \frac{1}{a_1(\delta_1-\beta)}    \left(\frac{\Phi^{-1}(T^{-1})}{\Phi^{-1}(t^{-1})}\right)^{\beta}
\Phi(\frac{\Phi^{-1}(t^{-1})}{\Phi^{-1}(T^{-1})}) \nn\\
&\le \frac{\Phi^{-1}(T^{-1})^{\beta-2}}{a_1(\delta_1-\beta)}   t^{-1}\Phi^{-1}(t^{-1})^{-\beta}\, .
\end{align}
The first inequality follows immediately by using the lower scaling to get that for $1 \ge r\ge \lambda^{-1}$, 
$\Phi(r^{-1})\le a_1^{-1} \lambda^{-\delta_1}r^{-\delta_1}\Phi(\lambda)$. The second inequality follows 
from \eqref{e:Berall}.

For the remainder of this paper we always assume that 
\eqref{e:intro-psibound} holds.
The following result is a generalization of \cite[Lemma 2.1]{CZ}.

\begin{lemma}\label{l:convolution}

 (a) For every $T\ge 1$, there exists $c_1=c_1(d,\delta_1, a_1, C_*, T, \Phi^{-1}(T^{-1}))>0$ such that for
$0< t\le T$, 
 all $\beta\in [0,\delta_1)$ and $\gamma\in \R$,
\begin{equation}
\label{e:convolution-integrability}
\int_{\R^d}\rho_{\gamma}^{\beta}(t,x)\, dx \le \frac{c_1}{\delta_1-\beta} t^{-1}\Phi^{-1}(t^{-1})^{-\gamma-\beta}\, .
\end{equation}

\noindent
(b)
For every $T\ge 1$, there exists $C_0=C_0(T)=C_0(d,\delta_1, a_1, C_*, T, \Phi^{-1}(T^{-1}))> 0$ such that for all $\beta_1,\beta_2  \ge 0$ with $\beta_1+ \beta_2  <\delta_1$,
 $\gamma_1,\gamma_2\in \R$ and  
$0<s<t\le T$, 
\begin{align}
\label{e:convolution-2}
&\int_{\R^d}\rho_{\gamma_1}^{\beta_1}(t-s,x-z) \rho_{\gamma_2}^{\beta_2}(s,z)\, dz \nn\\
& \le \frac{C_0}{ \delta_1 -\beta_1- \beta_2}\Big((t-s)^{-1}\Phi^{-1}((t-s)^{-1})^{-\gamma_1-\beta_1-\beta_2}\Phi^{-1}(s^{-1})^{-\gamma_2}\nn\\
& \  +\Phi^{-1}((t-s)^{-1})^{-\gamma_1}s^{-1}\Phi^{-1}(s^{-1})^{-\gamma_2-\beta_1-\beta_2} \Big) \rho(t,x)\nonumber \\
  & \  +\frac{C_0}{ \delta_1 -\beta_1- \beta_2}(t-s)^{-1}\Phi^{-1}((t-s)^{-1})^{-\gamma_1-\beta_1}\Phi^{-1}(s^{-1})^{-\gamma_2}\rho_0^{\beta_2}(t,x)\nonumber \\
 & \ +\frac{C_0}{ \delta_1 -\beta_1- \beta_2}\Phi^{-1}((t-s)^{-1})^{-\gamma_1}s^{-1}\Phi^{-1}(s^{-1})^{-\gamma_2-\beta_2}\rho_0^{\beta_1}(t,x)\, .
\end{align}

\noindent
(c) 
Let $T\ge 1$.
For  all $\beta_1,\beta_2  \ge 0$ with $\beta_1+ \beta_2  < \delta_1$,
and all $\theta,\eta \in[0,1], \gamma_1, \gamma_2\in \R$ satisfying
   ${\bf 1}_{\gamma_1\ge 0}(\gamma_1/2)  +{\bf 1}_{\gamma_1< 0}(\gamma_1/\delta_1)
   +\beta_1/2+1-\theta>0$ and
   ${\bf 1}_{\gamma_2\ge 0}(\gamma_2/2)  +{\bf 1}_{\gamma_2< 0}(\gamma_2/\delta_1)+\beta_2/2+1-\eta>0$,
 there exists $c_2>0$ such that for all $0<t \le T$ and $x\in \R^d$,
\begin{align}\label{e:convolution-3}
\lefteqn{\int_0^t \int_{\R^d} (t-s)^{1-\theta}\rho_{\gamma_1}^{\beta_1}(t-s, x-z)s^{1-\eta} \rho_{\gamma_2}^{\beta_2}(s,z)\, dz\, ds }\nonumber \\
& \le c_2 t^{2-\theta-\eta}\left(\rho_{\gamma_1+\gamma_2+\beta_1+\beta_2}^0+\rho_{\gamma_1+\gamma_2+\beta_2}^{\beta_1}+\rho_{\gamma_1+\gamma_2+\beta_1}^{\beta_2}\right)(t,x)\, .
\end{align}
Moreover, when we further assume that $\gamma_1, \gamma_2 \ge 0$,  we can take that
\begin{align}\label{e:convolution-4}
c_2=4\frac{C_0(T)}{ \delta_1 -\beta_1- \beta_2} B\left((\gamma_1+\beta_1)/2+1-\theta,{\gamma_2+\beta_2}/{2}+1-\eta\right).
\end{align}

\end{lemma}
\pf
(a) 
Let $c_1=c_1(d)=d|B(0,1)|$ and $T_1=\Phi^{-1}(T^{-1}) \le 1$. 
We have that for all $0<t \le T$,
\begin{align}
&\Phi^{-1}(t^{-1})^{\gamma} 
\int_{\R^d}\rho_{\gamma}^{\beta}(t,x)\, dx  
=\int_{\R^d}\left(|x|^{\beta} \wedge 1 \right)\rho(t,x)\, dx \nn\\
&\le c_1  \int_0^{T_1/\Phi^{-1}(t^{-1})} r^{\beta+d-1}\frac{\Phi\left(\left(
\frac{1}{\Phi^{-1}(t^{-1})}\right)^{-1}\right)}
{\left(\frac{T_1}{\Phi^{-1}(t^{-1})}\right)^{d}}dr+\nn\\
&\qquad +c_1 
\int_{T_1/\Phi^{-1}(t^{-1})}^{1}r^{\beta-1}\Phi(r^{-1})dr +c_1 
\int_{1}^{\infty}\frac{\Phi(r^{-1})}{r}dr \nn\\
&\le  \frac{c_1 T_1^{\beta} }{\beta+d}   t^{-1}\Phi^{-1}(t^{-1})^d \Phi^{-1}(t^{-1})^{-\beta-d} +c_1 \int_{T_1/\Phi^{-1}(t^{-1})}^{1}r^{\beta-1}\Phi(r^{-1})dr+c_1 \int_0^{1}\frac{\Phi(r)}{r}dr
\label{e:secnlast}\\
&\le c_1 d^{-1}  t^{-1}\Phi^{-1}(t^{-1})^{-\beta} + 
\frac{c_1 T_1^{\beta-2}}{a_1(\delta_1-\beta)}t^{-1}\Phi^{-1}(t^{-1})^{-\beta}
 +c_1 C_*  
\nn\\
&\le c_1 (  d^{-1}+ 
T_1^{-2}a_1^{-1}\delta_1^{-1}(\delta_1-\beta)^{-1}
+C_*a_1^{-1/2} T)t^{-1}\Phi^{-1}(t^{-1})^{-\beta},  \nn
\end{align}
where in the second to last line we used \eqref{e:integral-psi} to estimate the second term in \eqref{e:secnlast} and used \eqref{e:intro-psibound} to estimate the last term in \eqref{e:secnlast}, and in the  last line we used the assumption  $\beta\in [0,\delta_1)$ and  the inequality 
 $ t \Phi^{-1}(t^{-1})^{\beta} \le t(a_1^{-1/\delta_1}(T/t)^{1/\delta_1})^{\beta} \le a_1^{-\beta/\delta_1}T\le a_1^{-1}T$ which
follows from (2.4) with $\lambda=T/t$ and $r_0=r=T^{-1}$.

\noindent
(b) 
Let $c_2=(2(2/a_1)^{1/\delta_1} /\Phi^{-1}((2T)^{-1}))^{d+2}$.
As in the display after  \cite[(2.5)]{CZ}, we have that 
$$\left(|x-z|^{\beta_1}\wedge 1\right)\left(|z|^{\beta_2}\wedge 1\right) \le (|x-z|^{\beta_1+\beta_2}\wedge 1)+\left(|x-z|^{\beta_1}\wedge 1\right)\left(|x|^{\beta_2}\wedge 1\right).$$
By using this and \eqref{e:rho}, we have
\begin{align*}
\lefteqn{\rho_{\gamma_1}^{\beta_1}(t-s,x-z) \rho_{\gamma_2}^{\beta_2}(s,z)}\\
=&\Phi^{-1}((t-s)^{-1})^{-\gamma_1}\Phi^{-1}(s^{-1})^{-\gamma_2}\left(|x-z|^{\beta_1}\wedge 1\right)\left(|z|^{\beta_2}\wedge 1\right)\rho(t-s,x-z)\rho(s,z)\\
\le & c_2\Phi^{-1}((t-s)^{-1})^{-\gamma_1}\Phi^{-1}(s^{-1})^{-\gamma_2}\left(|x-z|^{\beta_1}\wedge 1\right)\left(|z|^{\beta_2}\wedge 1\right)\big(\rho(t-s,x-z)+\rho(s,z)\big)\rho(t,x)\\
\le & c_2\Phi^{-1}((t-s)^{-1})^{-\gamma_1}\Phi^{-1}(s^{-1})^{-\gamma_2}\Big\{(|x-z|^{\beta_1+\beta_2}\wedge 1)+\left(|x-z|^{\beta_1}\wedge 1\right)\left(|x|^{\beta_2}\wedge 1\right)\Big\}\\
&  \quad \times\rho(t-s,x-z)\rho(t,x)\\
&  + c_2\Phi^{-1}((t-s)^{-1})^{-\gamma_1}\Phi^{-1}(s^{-1})^{-\gamma_2}\Big\{(|z|^{\beta_1+\beta_2}\wedge 1)+\left(|x|^{\beta_1}\wedge 1\right)\left(|z|^{\beta_2}\wedge 1\right)\Big\} \rho(s,z)\rho(t,x)\\
=& c_2\Phi^{-1}(s^{-1})^{-\gamma_2}\Big\{\rho_{\gamma_1}^{\beta_1+\beta_2}(t-s,x-z)\rho(t,x)+\rho_{\gamma_1}^{\beta_1}(t-s,x-z)\rho_0^{\beta_2}(t,x)\Big\}\\
 & + c_2 \Phi^{-1}((t-s)^{-1})^{-\gamma_1}\Big\{\rho_{\gamma_2}^{\beta_1+\beta_2}(s,z)\rho(t,x)+\rho_{\gamma_2}^{\beta_2}(s,z) 
\rho_0^{\beta_1}(t,x)\Big\}\, . 
\end{align*}
Since
$\beta_1+ \beta_2  < \delta_1$, now \eqref{e:convolution-2} follows by integrating the above and using \eqref{e:convolution-integrability}.

\noindent 
(c)  By integrating \eqref{e:convolution-2} and using Lemma 
\ref{l:convoluton-inequality}, we get 
\eqref{e:convolution-3}. When we further assume that $\gamma_1, \gamma_2\ge0$,  by integrating \eqref{e:convolution-2} and using 
 the last part of Lemma \ref{l:convoluton-inequality}, 
we get \eqref{e:convolution-3} with the constant
\begin{align*}& C_0\left(B\left(\frac{\gamma_1+\beta_1+\beta_2}{2}+1-\theta,\frac{\gamma_2+2}{2}+1-\eta \right)+B\left(\frac{\gamma_2+\beta_1+\beta_2}{2}+1-\eta,\frac{\gamma_1+2}{2}+1-\theta\right)\right.\\
 &\left. +B\left(\frac{\gamma_1+\beta_1}{2}+1-\theta,\frac{\gamma_2+2}{2}+1-\eta\right) +B\left(\frac{\gamma_2+\beta_2}{2}+1-\eta,\frac{\gamma_1+2}{2}+1-\theta\right)\right),\end{align*}
which is,   using that the beta function $B$ is symmetric and non-increasing in each variable,  less than or equal to $4C_0 B\left((\gamma_1+\beta_1)/2+1-\theta,{\gamma_2+\beta_2}/{2}+1-\eta\right)$. 
\qed

\begin{lemma}\label{l:rsnewlemma}
Suppose $0<t_1\le t_2 < \infty$.
For $\beta\in (0, \delta_1/2)$,
$$
\lim_{h\downarrow0}\sup_{x, y\in \R^d, t\in [t_1, t_2]}\left(\int^h_0+\int^t_{t-h}\right)\int_{\R^d}
\rho^\beta_0(t-s, x-z)(\rho^\beta_0(s, z-y)+\rho_\beta^0(s, z-y))dzds=0.
$$
\end{lemma}

\pf 
We first apply Lemma \ref{l:convolution}(b) and then use Remark \ref{r:rho-decreasing},  to get that for $t\in  [t_1, t_2]$,
\begin{eqnarray*}
&&\int_{\R^d}\rho^\beta_0(t-s, x-z)(\rho^\beta_0(s, z-y)+\rho_\beta^0(s, z-y))dz\\
&&\le c_1((t-s)^{-1}\Phi^{-1}((t-s)^{-1})^{-\beta}+s^{-1}\Phi^{-1}(s^{-1})^{-\beta})\rho(t_1, 0).
\end{eqnarray*}
Now the conclusion of the lemma follows immediately from Lemma \ref{l:rsconvolution-inequality}.
\qed


 \section{Analysis of the  heat kernel of $\LL^{\mathfrak K}$}

 Throughout this paper, $Y=(Y_t,\P_x)$ is a subordinate Brownian motion via an independent subordinator with Laplace exponent $\phi$ and L\'evy measure $\mu$. 
The L\'evy density of $Y$, 
denoted by $j$, is given by
$$
j(x)=j(|x|)=\int_0^{\infty}(4\pi s)^{-d/2}e^{-|x|^2/4s}\, \mu(ds)\, .
$$ 
It is well known that  there exists $c=c(d)$ depending only on $d$ such that 
\begin{align}
\label{e:jupper}
j(r) \le  c\frac{\phi(r^{-2})}{r^{d}}, \quad r>0 
\end{align}
(see \cite[(15)]{BGR}).
The function $r\mapsto j(r)$ is non-decreasing.
Recall that we have assumed  that $r \mapsto \Phi(r)(=\phi(r^2))$, 
the radial part of the characteristic exponent $\Phi$ of $Y$, 
satisfies the weak lower scaling condition at infinity in \eqref{e:lower-scaling}.

Suppose that $Z=(Z_t,\P_x)$ is a purely discontinuous symmetric L\'evy process with characteristic exponent $\psi_Z$ 
such that 
its L\'evy measure admits a density $j_Z$ satisfying
\begin{equation}\label{e:assumption-j0}
{\wh \gamma_0}^{-1}j(|x|)\le j_Z(x) \le \wh \gamma_0 j(|x|)\, , \qquad x\in \R^d\, ,
\end{equation}
for some $\wh \gamma_0\ge 1$. Hence, 
 $\int_{\R^d}j_Z(x)dx=\infty$. 
The characteristic exponents of $Z$, respectively $Y$, are given by
$$
\psi_Z(\xi)=\int_{\R^d}(1-\cos(\xi\cdot y))j_Z(y)\, dy\, ,\qquad \Phi(\xi)=\int_{\R^d}(1-\cos(\xi\cdot y))j(|y|)\, dy\, ,
$$
and satisfy
\begin{equation}\label{e:assumption-j}
{\wh \gamma_0}^{-1}\Phi(|\xi|)\le \psi_Z(\xi) \le \wh \gamma_0 \Phi(|\xi|)\, ,\quad \xi \in \R^d\, .
\end{equation}
Let $\psi$ denote the radial nondecreasing majorant of the characteristic exponent of $Z$, i.e., $\psi(r):= \sup_{|z| \le r} \psi_Z(z)$.
Clearly
$$
\wh \gamma_0^{-1}\Phi(r)\le \psi(r) \le \wh \gamma_0 \Phi(r)\,, \quad   r>0,
\quad\text {and}\quad 
\wh \gamma_0^{-2}\psi(|\xi|)\le \psi_Z(\xi) \le  \psi(|\xi|)\, ,\quad \xi \in \R^d,
$$
and thus $\psi$ also satisfies the weak lower scaling condition at infinity in \eqref{e:lower-scaling}.

By \eqref{e:jupper} and 
\eqref{e:assumption-j0},
\begin{align}
\label{e:jZupper}
j_Z(x) \le 
\wh \gamma_0 \frac{\Phi(|x|^{-1})}{|x|^{d}}.  
\end{align}
Moreover,   for every $n\in \Z_+$,
\begin{align}
\int_{\R^d}&\left|\E\left[e^{i\xi\cdot Z_t}\right]\right||\xi|^n\, d\xi\nn=\int_{\R^d}e^{-t\psi_Z(\xi)}|\xi|^n\, d\xi \le \int_{\R^d}e^{-t\wh \gamma_0^{-1}\Phi(|\xi|)}|\xi|^n\, d\xi \\
&\le  c \left(\int_0^{1}r^{d-1+n} dr+\int_1^{\infty}r^{d-1+n}e^{-t\wh \gamma_0^{-1}a_1r^{\delta_1}}\, dr\right) 
<\infty\, . 
\label{e:phifinte}
\end{align}
It follows from \cite[Proposition 2.5(xii) and Proposition 28.1]{Sat} that $Z_t$ has a density 
$$
p(t,x)=(2\pi)^{-d/2}\int_{\R^d} e^{-ix\cdot \xi}e^{-t\psi_Z(\xi)}\, d\xi =(2\pi)^{-d/2}\int_{\R^d} \cos(x\cdot \xi)e^{-t\psi_Z(\xi)}\, d\xi,
$$
which is infinitely differentiable in $x$.
Let $\LL$ be the infinitesimal generator of $Z$.

\begin{lemma}\label{l:partial-time}
(a) For every $x\in \R^d$, the function $t\mapsto p(t,x)$ is differentiable and 
$$
\frac{\partial p(t,x)}{\partial t}=
-(2\pi)^{-d/2}\int_{\R^d} \cos(x\cdot \xi)\psi_Z(\xi)e^{-t\psi_Z(\xi)}\, d\xi =\LL p(t,x)\, .
$$

\noindent
(b) For every $\eps >0$ there exists a constant 
$c=c(d, \delta_1, a_1, \wh \gamma_0, \eps)>0$
 such that for all $s,t \ge \eps$ and all $x,y\in \R^d$,
$$
|p(t,x)-p(s,y)|\le c \left(|t-s|+|x-y|\right)\, .
$$
\end{lemma}
\pf 
(a) Note that for any $t\ge 0$ and any $h\in \R$ such that $t+h\ge 0$,
$$
\frac{p(t+h,x)-p(t,x)}{h}=
(2\pi)^{-d/2}\int_{\R^d} \cos(x\cdot \xi)e^{-t\psi_Z(\xi)}\frac{e^{-h\psi_Z(\xi)}-1}{h}d\xi.
$$
The absolute value of the integrand is bounded by 
$2\wh \gamma_0 \Phi(|\xi|)e^{-\wh \gamma_0^{-1}\Phi(|\xi|)}$ 
which is integrable since $\Phi(|\xi|)\le |\xi|^2$. The claim follows from the dominated convergence theorem by letting $h\to 0$. 
The last equality in the statement of the lemma follows from \cite[Example 4.5.5]{Jac}.

\noindent 
(b) By the triangle inequality we have that
\begin{eqnarray*}
\lefteqn{|p(t,x)-p(s,y)|\le \int_{\R^d}\left| \cos(x\cdot \xi)-\cos(y\cdot \xi)\right| e^{-t\psi_Z(\xi)}\, d\xi}\\
&  +&\int_{\R^d}\left|\cos(y\cdot \xi)\right| \left|e^{-t\psi_Z(\xi)}-e^{-s\psi_Z(\xi)}\right |\, d\xi =:I_1+I_2\, .
\end{eqnarray*}
Clearly, $|\cos(x\cdot \xi)-\cos(y\cdot \xi)|\le |x\cdot \xi -y\cdot \xi|\le |x-y||\xi|$, which implies that, by \eqref{e:phifinte},
$$
I_1\le |x-y|\int_{\R^d} |\xi|e^{-t\psi_Z(\xi)}\, d\xi\le |x-y|
\int_{\R^d} |\xi|e^{-\eps \wh \gamma_0^{-1}\Phi(|\xi|)}\, d\xi 
=c_1(\wh \gamma_0, \eps)|x-y|\, .
$$
In order to estimate $I_2$, without loss of generality we assume that $s\le t$. Then by the mean value theorem we have that
$$
\left| e^{-t\psi_Z(\xi)}-e^{-s\psi_Z(\xi)}\right| \le |t-s| \psi_Z(\xi)e^{-s\psi_Z(\xi)} 
\le \wh \gamma_0 |t-s| \Phi(|\xi|) e^{-\eps \wh \gamma_0^{-1}\Phi(|\xi|)}\, .
$$
Therefore, by \eqref{e:phifinte},
$$
I_2\le 
\wh \gamma_0 |t-s| \int_{\R^d}|\xi|^2
e^{-\eps \wh \gamma_0^{-1}\Phi(|\xi|)}\, d\xi =c_2(\wh \gamma_0, \eps) |t-s|\, .
$$
The claim follows by taking $c=c_1\vee c_2$.
\qed

Define the Pruitt function $\sP$ by
\begin{align}
\label{e:Pruitt}
\sP (r)= \int_{\R^d} \left(\frac{|x|^2}{r^2} \wedge 1\right) j(x)dx.
\end{align}
By \cite[(6) and Lemma 1]{BGR}, 
\begin{align}\label{e:Ppsi}
\frac{1}{2\wh \gamma_0}\psi(r^{-1}) \le \frac{1}{2}\Phi(r^{-1}) \le \sP (r) \le \frac{d \pi^2}{2}\Phi(r^{-1})\le  \frac{\wh \gamma_0 d \pi^2}{2}\psi(r^{-1}).
\end{align}
In this paper we will use \eqref{e:Ppsi} several times.

We next discuss the upper estimate of $p(t,x)$ and its derivatives for 
$0<t \le T$
 and all $x\in \R^d$ using \cite[Theorem 3]{KS}.

\begin{prop}\label{p:upperestonp}
For each $T \ge 1$ and 
$k\in \bZ_+$, there
is a constant 
$c=c(k, T, \wh \gamma_0, d, \delta_1, a_1)\ge 1$ such that
$$
|\nabla^k p(t,x)|\,\le \,c\,t \,( \Phi^{-1}(t^{-1}))^k \rho(t,x)\, ,\qquad 
0<t\le T, x\in \R^d,
$$
where $\nabla^k$ stands for the $k$-th order gradient with respect to the spatial variable $x$.
\end{prop}
\pf 
First, we recall that $\int_{\R^d}j_Z(x)dx=\infty$.
Let $f(s):=\frac{\Phi(s^{-1})}{s^d}$. Then by \eqref{e:jZupper} we have 
$j_Z(x) \le C \wh \gamma_0 f(|x|)$. 
Thus for $A \in \sB(\R^d)$,
$$
\int_{A}j_Z(x)dx \le C  \wh \gamma_0  \int_{A} \frac{\Phi(|x|^{-1})}{|x|^d}dx 
\le C \wh \gamma_0  \frac{\Phi(\text{dist}(0, A)^{-1})}{\text{dist}(0, A)^d} |A|
\le C  \wh \gamma_0 f(\text{dist}(0, A))
(\text{diam}(A))^d.
$$
Therefore, \cite[(1)]{KS} holds with $\gamma=d$ and 
$M_1=C  \wh \gamma_0$. 

Since $(s \vee|y|) -(|y|/2) \ge s/2$ for  $s>0$, using \eqref{e:Ppsi}
in the last inequality we have that for  $s, r>0$, 
\begin{align}
\label{e:(2)1}
&\int_{|y|>r} f( (s \vee|y|) -(|y|/2)) j_Z(y)dy  \le 2^d \frac{\Phi((s/2)^{-1})}{s^d}
\int_{|y|>r}   j_Z(y)dy \nn\\
&=  2^d \frac{\Phi((s/2)^{-1})}{s^d}
\int_{|y|>r} (\frac{|y|^2}{r^2} \wedge 1)  j_Z(y)dy \le 2^{d+2} \frac{\Phi(s^{-1})}{s^d}  \sP (r) 
 \le 2^{d+1} \wh \gamma_0 d \pi^2 f(s) \psi(r^{-1}).
\end{align}
Therefore, \cite[(2)]{KS} holds with 
$M_1=2^{d+1} \wh \gamma_0 d \pi^2$. 

Furthermore, by \eqref{e:assumption-j} and \eqref{e:lower-scaling}, for 
$k\in \bZ_+$, 
\begin{align*}
&\int_{\R^d}e^{-t\psi_Z(\xi)}|\xi|^k\, d\xi \le \int_{\R^d}
e^{-t\wh \gamma_0^{-1}\Phi(|\xi|)}|\xi|^k\, d\xi  
= d  |B(0,1)| \int_0^{\infty}r^{d+k-1}
e^{-t\wh \gamma_0^{-1}\Phi(r)}\, dr\\
&= d  |B(0,1)|\int_0^{\infty}(\Phi^{-1}(s/t))^{d+k-1}
e^{-\wh \gamma_0^{-1}s} (\Phi^{-1})'(s/t)t^{-1} \, ds\\
& \le d  |B(0,1)| \int_0^{1}(\Phi^{-1}(s/t))^{d+k-1} (\Phi^{-1})'(s/t)t^{-1} \, ds\\
&\quad +d  |B(0,1)| \sum_{n=1}^\infty e^{-\wh \gamma_0^{-1}2^{n-1}}
\int_{2^{n-1}}^{2^{n}}(\Phi^{-1}(s/t))^{d+k-1} (\Phi^{-1})'(s/t)t^{-1} \, ds\\
& =  \frac{d |B(0,1)|}{d+k} \int_0^{1}((\Phi^{-1}(s/t))^{d+k})'  \, ds+
\frac{d |B(0,1)|}{d+k} 
\sum_{n=1}^\infty e^{-\wh \gamma_0^{-1}2^{n-1}}\int_{2^{n-1}}^{2^{n}}
((\Phi^{-1}(s/t))^{d+k})' \, ds\\
&  \le \frac{d |B(0,1)|}{d+k} \left((\Phi^{-1}(t^{-1}))^{d+k}+ 
\sum_{n=1}^\infty 
e^{-\wh \gamma_0^{-1}2^{n-1}}(\Phi^{-1}(2^n/t))^{d+k} \right).
\end{align*}
Since $t\le T$, by \eqref{e:lsc-inverse} we have
$ \Phi^{-1}(2^n/t) \le c_0 2^{n/\delta_1}\Phi^{-1}(t^{-1})$.
Thus we see that 
\begin{align*}
&\int_{\R^d}e^{-t\psi_Z(\xi)}|\xi|^k\, d\xi \le 
\frac{d |B(0,1)|}{d+k}(\Phi^{-1}(t^{-1}))^{d+k} (1+  c_0
\sum_{n=1}^\infty  2^{n(d+k)/\delta_1} e^{-
\wh \gamma_0^{-1}2^{n-1}}) \\
&\le c_1 \Phi^{-1}(t^{-1})^{d+k}\le c_2 \psi^{-}(t^{-1})^{d+k},
\end{align*}
where $c_2=c_2(k)>0$ and $\psi^-$ is the generalized inverse of $\psi$:
$
\psi^-(s)=\inf\{u\ge 0:\, \psi(u)\ge s\}.$
Therefore, \cite[(8)]{KS} holds with the set $(0,T]$. 

We have checked that the conditions in  \cite[Theorem 3]{KS} hold for all 
$k\in \bZ_+$.
Thus by \cite[Theorem 3]{KS} (with $n=d+2$ in  \cite[Theorem 3]{KS}),  
there exists $c_3(k)>0$ such that 
for $t \le T$,
\begin{align*}
&|\nabla^k p(t,x)| \le c_3\psi^-(t^{-1})^k\left(\psi^-(t^{-1})^d \wedge \left(\frac{t\Phi(|x|^{-1})}{|x|^d}
+\frac{ \psi^-(t^{-1})^d }{(1+ |x| \psi^-(t^{-1}))^{d+2}} \right)  \right) \\
&\le c_4\Phi^{-1}(t^{-1})^k\left(\Phi^{-1}(t^{-1})^d \wedge \left(\frac{t\Phi(|x|^{-1})}{|x|^d}+\frac{ \Phi^{-1}(t^{-1})^d}{  (1+ |x| \Phi(t^{-1}))^{d+2}}\right)  \right).
\end{align*}
When  $|x| \Phi^{-1}(t^{-1}) \ge 1$ (so that $t\Phi(|x|^{-1}) \le 1$),
\begin{align*}
\frac{ \Phi^{-1}(t^{-1})^d}{  (1+ |x| \Phi(t^{-1}))^{d+2}} \le \frac{ \Phi^{-1}(t^{-1})^d}{ ( |x|\Phi^{-1}(t^{-1}))^{d+2}}
=   |x|^{-d} \left( \frac{\Phi^{-1}( \Phi(|x|^{-1}))}{\Phi^{-1}(
\frac{\Phi(|x|^{-1})}{t\Phi(|x|^{-1})})}\right)^2 \le  |x|^{-d} (t\Phi(|x|^{-1})).
  \end{align*}
  In the  last inequality we have used \eqref{e:psi-inverse-sbm}.
  Therefore using Proposition \ref{p:p-q} we conclude that for all $0<t\le T$ and $x\in \R^d$, 
  $$
|\nabla^k p(t,x)|\le c_{4}\Phi^{-1}(t^{-1})^k\left(\Phi^{-1}(t^{-1})^d \wedge \frac{t\Phi(|x|^{-1})}{|x|^d}\right)  \le c_{4}2^{d+2} t\Phi^{-1}(t^{-1})^k \rho(t,x)\,.
$$
\qed

\subsection{Further properties of $p(t, x)$}
We will need the following simple inequality, cf.~\cite[(2.9)]{CZ}: Let $a>0$ and $x\in \R^d$. For every $z\in \R^d$ such that $|z|\le (2a)\vee (|x|/2)$, we have
\begin{equation}\label{e:2.9}
(a+|x+z|)^{-1}\le 4(a+|x|)^{-1}\, .
\end{equation}
Indeed, if $|z|\le 2a$, then $a+|x|\le a+|x+z|+|z|\le a+|x+z|+2a \le 4(a+|x+z|)$. If $|z|\le |x|/2$, then $4(a+|x+z|)\ge 4a +4|x|-4|z|\ge 4a +4|x|-2|x|\ge a+|x|$.

For a function $f:\R_+\times \R^d\to \R$, we define
\begin{equation}\label{e:delta-f-def}
\delta_f(t,x;z):=f(t,x+z)+f(t,x-z)-2f(t,x)\, .
\end{equation}
Also, $f(x\pm z)$ is an abbreviation for $f(x+z)+f(x-z)$.

The following result is the counterpart of \cite[Lemma 2.3]{CZ}.

\begin{prop}\label{p:gradient-est-kappa}
For every $T\ge 1$, there exists a constant 
$c=c(T, d, \wh \gamma_0, d, \delta_1, a_1)>0$ 
such that for every $t\in (0,T]$ and $x,x',z\in\R^d$,
\begin{equation}\label{e:difference-p-kappa}
\left|p(t,x)-p(t,x')\right|\le c\left((\Phi^{-1}(t^{-1})â|x-x'|)\wedge 1\right)t\left(\rho(t,x)+\rho(t,x')\right)\, ,
\end{equation}
\begin{equation}\label{e:delta-p-kappa}
\left|\delta_{p}(t,x;z)\right| \le c\left((\Phi^{-1}(t^{-1})|z|)^2\wedge 1\right)t\left(\rho(t,x\pm z)+\rho(t,x)\right),
\end{equation}
and
\begin{eqnarray}\label{e:delta-difference-kappa}
\lefteqn{|\delta_{p}(t,x;z)-\delta_{p}(t,x';z)|
\le c
\left((\Phi^{-1}(t^{-1})|x-x'|)\wedge 1\right)\left((\Phi^{-1}(t^{-1})|z|)^2\wedge 1\right)} \nonumber \\
& & \qquad \qquad \qquad \times t \left(\rho(t,x\pm z)+\rho(t,x)+\rho(t,x'\pm z)+\rho(t,x')\right)\, .
\end{eqnarray}
\end{prop}
\pf 
(1) 
Note that, by Proposition \ref{p:upperestonp} with $k=0$, \eqref{e:difference-p-kappa} is clearly true
if $\Phi^{-1}(t^{-1})|x-y|\ge 1$. Thus we assume that $\Phi^{-1}(t^{-1})|x-y|\le 1$. We use Proposition \ref{p:upperestonp} 
for $k=1$ and
\begin{equation}\label{e:difference-p-grad}
p(t,x)-p(t,y)=(x-y)\cdot \int_0^1\nabla p(t, x+\theta (y-x))\, d\theta
\end{equation}
to estimate
$
|p(t,x)-p(t,y)|\le  c_1 t\Phi^{-1}(t^{-1})  |x-y|\int_0^1  \rho(t,x+\theta (y-x))
d \theta\, .$
Since $\theta|y-x|\le 1/\Phi^{-1}(t^{-1})$, we get from \eqref{e:2.9} that
$$
\left(\frac{1}{\Phi^{-1}(t^{-1})}+|x+\theta(y-x)|\right)^{-1}\le 4 \left(\frac{1}{\Phi^{-1}(t^{-1})}+|x|\right)^{-1}\, .
$$
Therefore  using \eqref{e:Berall}  we have
$
|p(t,x)-p(t,y)|\le  c_2|x-y|\Phi^{-1}(t^{-1}) t\rho(t,x)$,  $t\in (0,T]$.

\noindent
(2) Note that \eqref{e:delta-p-kappa} is clearly true if $\Phi^{-1}(t^{-1})|z|\ge 1$.  In order to prove \eqref{e:delta-p-kappa} when $\Phi^{-1}(t^{-1})|z|\le 1$ we use \eqref{e:difference-p-grad} twice to obtain
\begin{eqnarray}\label{e:difference-p-grad-2}
\delta_{p}(t,x;z)
&=&z\cdot \int_0^1 \left(\nabla p(t,x+\theta z)-\nabla p(t,x-\theta z)\right)\, d\theta \nonumber\\
&=&2(z\otimes z)\cdot \int_0^1 \int_0^1 \theta \nabla^2 p(t, x+(1-2\theta')\theta z)\, d\theta'\, d\theta\, .
\end{eqnarray}
Note that $|(1-2\theta')\theta z|\le |z|\le \frac{1}{\Phi^{-1}(t^{-1})}$. Hence, by Proposition \ref{p:upperestonp} and \eqref{e:2.9} we get the estimate
$$
\left|\theta \nabla^2 p(t, x+(1-2\theta')\theta z)\right| \le c_3 \left(\Phi^{-1}(t^{-1})\right)^2 t\rho(t,x)\, .
$$
Therefore,
$
\delta_{p}(t,x;z)\le  c_4 \left(\Phi^{-1}(t^{-1})|z|\right)^2 t\rho(t,x)$, $t\in (0,T]$.

\noindent
(3)
It follows from \eqref{e:delta-p-kappa} that it suffices to prove
\eqref{e:delta-difference-kappa} in the case when $\Phi^{-1}(t^{-1})|x-y|\le 1$. To do this, 
we start with the subcase when $\Phi^{-1}(t^{-1})|z|\le 1$ and $\Phi^{-1}(t^{-1})|x-y|\le 1$. Then by \eqref{e:difference-p-grad-2},
\begin{eqnarray*}
\lefteqn{|\delta_{p}(t,x;z)-\delta_{p}(t,y;z)|}\nonumber \\
&\le & c_5 |x-y|\cdot |z|^2 \int_0^1  \int_0^1  \int_0^1 |\nabla^3 p(t,x+(1-2\theta')\theta z+\theta''(y-x))|\, d\theta d\theta' d\theta'' . \nonumber
\end{eqnarray*}
Note that $|(1-2\theta')\theta z+\theta''(y-x))|\le \frac{2}{\Phi^{-1}(t^{-1})}$. Hence, by Proposition \ref{p:upperestonp} and \eqref{e:2.9} we get 
\begin{eqnarray*}
|\delta_{p}(t,x;z)-\delta_{p}(t,y;z)|
& \le & c_6
 \Phi^{-1}(t^{-1})|x-y|(\Phi^{-1}(t^{-1})|z|)^2 t\rho(t,x)\, .
\end{eqnarray*}
If $\Phi^{-1}(t^{-1})|z|\ge 1$ and $\Phi^{-1}(t^{-1})|x-y|\le 1$, then again by Proposition \ref{p:upperestonp} and \eqref{e:2.9},
\begin{eqnarray*}
\lefteqn{|\delta_{p}(t,x;z)-\delta_{p}(t,y;z)|}\\
&\le & c_7\left( |x-y|\int_0^1 |\nabla p(t, x\pm z+\theta (y-x))|\, d\theta+|x-y|\int_0^1 |\nabla p(t, x+\theta(y-x))|\, d\theta\right)\\
& \le & c_8 \Phi^{-1}(t^{-1}) |x-y|\left(t\rho(t,x\pm z)+t\rho(t,x)\right)\,, \quad t\in (0,T].
\end{eqnarray*}
\qed

The following result is the counterpart of \cite[Theorem 2.4]{CZ}.

\begin{thm}\label{t:fract-der-est} For every $T\ge 1$, there exists a constant 
$c=c(T, d, \wh \gamma_0, d, \delta_1, a_1)>0$ 
such that for all $t\in(0,T]$ and all $x,x'\in\R^d$,
\begin{equation}\label{e:fract-der-est1}
\int_{\R^d} \left|\delta_{p}(t,x;z)\right| j(|z|)\, dz \le c \rho(t,x)
\end{equation}
and
\begin{align}\label{e:fract-der-est2}
\int_{\R^d}\left|\delta_{p}(t,x;z)-\delta_{p}(t,x';z)\right| j(|z|)\, dz
\le  c\left((\Phi^{-1}(t^{-1})|x-x'|)\wedge 1\right)(\rho(t,x)+\rho(t,x')) \, .
\end{align}
\end{thm}

\pf
By \eqref{e:delta-p-kappa} we have
\begin{eqnarray}
\lefteqn{\int_{\R^d} \left|\delta_{p}(t,x;z)\right| j(|z|)\, dz} \nonumber \\
&\le &c_0 \int_{\R^d}\left((\Phi^{-1}(t^{-1})|z|)^2\wedge 1\right)t\left(\rho(t,x\pm z)+\rho(t,x)\right)j(|z|)\, dz \label{e:fract-der-est3}\\
&=&c_0 \left(\int_{\R^d}\left((\Phi^{-1}(t^{-1})|z|)^2\wedge 1\right)t\rho(t,x\pm z)j(|z|)\, dz + t\rho(t,x)\sP (1/\Phi^{-1}(t^{-1}))
\right)\nonumber \\
&=:& c_0\, (I_1+I_2)\, .\nonumber
\end{eqnarray}
Clearly by \eqref{e:Ppsi},
$
I_2 \le   c_1 t \rho(t,x) \Phi (\Phi^{-1}(t^{-1})) = c_1 \rho(t,x)\, .
$
Next,\begin{eqnarray*}
I_1&=&\Phi^{-1}(t^{-1})^2\int_{\Phi^{-1}(t^{-1})|z|\le 1}|z|^2 t \rho(t,x\pm z)j(|z|)\, dz + \int_{\Phi^{-1}(t^{-1})|z|> 1}t \rho(t,x\pm z)j(|z|)\, dz\\
&=:&I_{11}+I_{12}\, .
\end{eqnarray*}
By using \eqref{e:2.9} in the first inequality below and \eqref{e:Ppsi} in the third, we further have
\begin{eqnarray*}
I_{11}&\le & 4^{d+1}  t \rho(t,x) \int_{|z|\le \frac{1}{\Phi^{-1}(t^{-1})}}
((\Phi^{-1}(t^{-1})  |z|)^2 \wedge 1)j(|z|)\, dz\\
& \le  &  4^{d+1}   t \rho(t,x)  \sP (1/\Phi^{-1}(t^{-1})) \le c_2 \rho(t,x)\, .
\end{eqnarray*}
Next, we have
\begin{eqnarray*}
I_{12} &\le&  t\int_{|z|> \frac{1}{\Phi^{-1}(t^{-1})}}{\Phi\left(\left(\frac{1}{\Phi^{-1}(t^{-1})}\right)^{-1}\right)}{\left(\frac{1}{\Phi\-(t^{-1})}\right)^{-d}}j(|z|)\, dz\\
&=&  \Phi^{-1}(t^{-1})^d \int_{|z|> \frac{1}{\Phi^{-1}(t^{-1})}} ((\Phi^{-1}(t^{-1})  |z|)^2 \wedge 1) j(|z|)\, dz\\
&\le &  \Phi^{-1}(t^{-1})^d \sP (1/\Phi^{-1}(t^{-1}))
\le  c_3 \Phi^{-1}(t^{-1})^d \Phi (\Phi^{-1}(t^{-1}))
= c_3 \Phi^{-1}(t^{-1})^d t^{-1}\, ,
\end{eqnarray*}
where in the last line we used \eqref{e:Ppsi}.
If $|x|\le 2/\Phi^{-1}(t^{-1})$, we have that 
$$
 \rho(t,x) \ge  
{\Phi\left(\left(\frac{3}{\Phi^{-1}(t^{-1})}\right)^{-1}\right)}
{\left(\frac{3}{\Phi\-(t^{-1})}\right)^{-d}}
\ge  c_4 t^{-1} \Phi^{-1}(t^{-1})^{d},
$$
implying that $I_{12}\le c_5 \rho(t,x)$.

If $|x|> 2/\Phi^{-1}(t^{-1})$, then by  \eqref{e:Ppsi},
\begin{eqnarray*}
\lefteqn{I_{12}=\left(\int_{\frac{|x|}{2}\ge |z|>\frac{1}{\Phi^{-1}(t^{-1})}}+\int_{|z|>\frac{|x|}{2}}\right)t \rho(t,x\pm z)j(|z|)\, dz}\\
&\le & c_{6}\left(t \rho(t,x)\int_{\frac{|x|}{2}\ge |z|>\frac{1}{\Phi^{-1}(t^{-1})}}j(|z|)\, dz +j(|x|/2) \int_{|z|>\frac{|x|}{2}}t \rho(t,x\pm z)\, dz\right)\\
&\le &c_{7}\left(t \rho(t,x)\int_{|z|>\frac{1}{\Phi^{-1}(t^{-1})}}j(|z|)\, dz + \frac{\Phi(2 |x|^{-1})}{|x|^d}\int_{\R^d}t \rho(t,x\pm z)\, dz \right)\\
&\le &c_{7}\left(t \rho(t,x)\sP (1/\Phi^{-1}(t^{-1})) +
 \frac{\Phi(|x|^{-1})}{|x|^d} \right)\\
&\le & c_{8}\left(\rho(t,x)+ \frac{\Phi(|x|^{-1})}{|x|^d}\right) \le c_{9} \rho(t,x)\, ,
\end{eqnarray*}
where in the last line the second term is estimated by a constant times the first term in view of the assumption that $|x|> 2/\Phi^{-1}(t^{-1})$.
This finishes the proof of \eqref{e:fract-der-est1}.

Next, by \eqref{e:delta-difference-kappa} we have 
\begin{eqnarray*}
\lefteqn{\int_{\R^d}\left|\delta_{p}(t,x;z)-\delta_{p}(t,x';z)\right| j(|z|)\, dz \le c_{10}\left((\Phi^{-1}(t^{-1})|x-x'|)\wedge 1\right)}\\
& &\times\left\{\int_{\R^d}\left((\Phi^{-1}(t^{-1})|z|)^2\wedge 1\right)
\big(t \rho(t,x\pm z)+t \rho(t,x'\pm z)\big)j(|z|)\, dz\right.\\
& & \left. \quad + \big(t \rho(t,x)+t \rho(t,x')\big)\int_{\R^d}\left((\Phi^{-1}(t^{-1})|z|)^2\wedge 1\right)j(|z|)\, dz\right\}\\
&\le &c_{11} \left((\Phi^{-1}(t^{-1})|x-x'|)\wedge 1\right) t^{-1}\big(t \rho(t,x)+t \rho(t,x')\big)\, ,
\end{eqnarray*}
where the last line follows by using the estimates of the integrals $I_1$ and $I_2$ from the first part of the proof. 
\qed

\subsection{Continuous dependence of heat kernels with respect to $\mathfrak{K}$}

Recall that $J: \R^d \to (0, \infty)$ is a  symmetric function satisfying \eqref{e:psi1}.
We now specify the jumping kernel $j_Z$. Let $\mathfrak{K}:\R^d\to (0,\infty)$ be a symmetric function, that is, $\mathfrak{K}(z)=\mathfrak{K}(-z)$. Assume  that there are $0<\kappa_0\le \kappa_1<\infty$ such that
\begin{equation}\label{e:bounds-for-kappa}
\kappa_0\le \mathfrak{K}(z)\le \kappa_1\, ,\qquad \text{for all }z\in \R^d\, .
\end{equation}
Let $j^{\mathfrak{K}}(z):=\mathfrak{K}(z)J(z)$, $z\in \R^d$. Then $j^{\mathfrak{K}}$ satisfies \eqref{e:assumption-j0} 
with ${\wh \gamma_0}=\gamma_0(\kappa_1 \vee \kappa_0^{-1})$. The infinitesimal generator of the corresponding symmetric L\'evy process $Z^{\mathfrak{K}}$ is given by
\begin{eqnarray}
{\mathcal L}^{\mathfrak{K}}f(x)&=&\mathrm{p.v.} \int_{\R^d}(f(x+z)-f(x))\mathfrak{K}(z)J(z)\, dz \nonumber\\
&=&\frac12 \mathrm{p.v.}\int_{\R^d} \delta_f(x;z)\mathfrak{K}(z)J(z)\, dz\, .\label{e:representaion-L}
\end{eqnarray}
We note in passing that, when $f\in C^2_b(\R^d)$, it is not necessary to take the
principal value in the last line above.
The transition density of $Z^{\mathfrak{K}}$ (i.e., the heat kernel of ${\mathcal L}^{\mathfrak{K}}$) will be denoted by $p^{\mathfrak{K}}(t,x)$. 
Then by Lemma \ref{l:partial-time},
\begin{eqnarray}\label{e:prop-p-L}
\frac{ \partial p^{\mathfrak{K}}(t,x)}{\partial t}={\LL}^{\mathfrak{K}}p^{\mathfrak{K}}(t,x)\, ,\qquad \lim_{t\to 0}p^{\mathfrak{K}}(t,x)=\delta_0(x)\, . 
\end{eqnarray}
 
 We will need the following observation for the next result. 
The inequality 
\eqref{e:lsc-inverse}
implies that there exists a constant $c(\kappa_0)\ge 1$ such that 
$$ 
\Phi^{-1}((\kappa_0t/2)^{-1})\le 
a_1^{-1/\delta_1} \Phi^{-1}(T^{-1})^{-1} 
(1 \vee (\kappa_0/2))^{1/\delta_1} 
\Phi^{-1}(t^{-1})\quad \text{ for all } t\in (0,T].
$$ Consequently, for all $z\in \R^d$ and $t\in (0,T]$,
\begin{equation}\label{e:comparison-kappa0}
\left(\Phi^{-1}((\kappa_0 t/2)^{-1})|z|\right)\wedge 1 \le 
a_1^{-1\delta_1} \Phi^{-1}(T^{-1})^{-1} 
(1 \vee (\kappa_0/2))^{1/\delta_1} 
\left(\left(\Phi^{-1}(t^{-1})|z|\right)\wedge 1 \right).
\end{equation}

The following result is the counterpart of \cite[Theorem 2.5]{CZ}, and in its proof
we follow the proof of \cite[Theorem 2.5]{CZ} with some modifications.

\begin{thm}
For every $T\ge 1$, there exists a constant 
$c>0$ depending on $T, d, \kappa_0, \kappa_1, \gamma_0, a_1$ and $\delta_1$ 
such that
for any  two symmetric functions $\mathfrak{K}_1$ and ${{\mathfrak{K}}}_2$ in $\R^d$ satisfying \eqref{e:bounds-for-kappa}, every $t\in(0,T]$ and $x\in \R^d$, we have 
\begin{eqnarray}
\left|p^{{\mathfrak{K}}_1}(t,x)-p^{{\mathfrak{K}}_2}(t,x)\right| &\le  & c \|{\mathfrak{K}}_1-{\mathfrak{K}}_2\|_{\infty}\ t \rho(t,x)\, ,\label{e:difference-kappa}\\
\left|\nabla p^{{\mathfrak{K}}_1}(t,x) -\nabla p^{{\mathfrak{K}}_2}(t,x)\right| &\le  & c  \|{\mathfrak{K}}_1-{\mathfrak{K}}_2\|_{\infty}\Phi^{-1}(t^{-1}) t\rho(t,x)\label{e:difference-nkappa}
\end{eqnarray}
and 
\begin{equation}\label{e:delta-difference-abs}
\int_{\R^d}\left|\delta_{p^{{\mathfrak{K}}_1}}(t,x;z)-\delta_{p^{{\mathfrak{K}}_2}}(t,x;z)\right| j(|z|)\, dz \le c \|{\mathfrak{K}}_1-{\mathfrak{K}}_2\|_{\infty} \rho(t,x)\, .
\end{equation}
\end{thm}

\pf (i) 
Using \eqref{e:prop-p-L} in  the second line,
the fact $\LL^{{\mathfrak{K}}_1}$ is self-adjoint in the third and fourth lines,  we have
\begin{eqnarray*}
\lefteqn{p^{{\mathfrak{K}}_1}(t,x)-p^{{\mathfrak{K}}_2}(t,x)=\int_0^t \frac{d}{ds}\left(\int_{\R^d}p^{{\mathfrak{K}}_1}(s,y)p^{{\mathfrak{K}}_2}(t-s,y-x)\, dy\right) ds}\\
&=&\int_0^t \left(\int_{\R^d}\left(\LL^{{\mathfrak{K}}_1}p^{{\mathfrak{K}}_1}(s,\cdot) (y) p^{{\mathfrak{K}}_2}(t-s,y-x)-p^{{\mathfrak{K}}_1}(s,y)\LL^{{\mathfrak{K}}_2}p^{{{\mathfrak{K}}_2}}(t-s,\cdot) (y-x)\right)dy\right)ds\\
&=&\int_0^{t/2} \left(\int_{\R^d} p^{{\mathfrak{K}}_1}(s,y)\left(\LL^{{\mathfrak{K}}_1}-\LL^{{\mathfrak{K}}_2}\right) p^{{\mathfrak{K}}_2}(t-s,\cdot) (y-x)dy\right)ds\\
&&+\int_{t/2}^t \left(\int_{\R^d} \left(\LL^{{\mathfrak{K}}_1}-\LL^{{\mathfrak{K}}_2}\right) p^{{\mathfrak{K}}_1}(s,\cdot) (y) p^{{\mathfrak{K}}_2}(t-s,y-x)dy\right)ds\\
&=&\frac12 \int_0^{t/2} \left(\int_{\R^d} p^{{\mathfrak{K}}_1}(s,y)\left(\int_{\R^d}\delta_{p^{{\mathfrak{K}}_2}}(t-s,x-y; z)({\mathfrak{K}}_1(z)-{\mathfrak{K}}_2(z))J(z)dz\right)dy\right)ds\\
&&+\frac12 \int_{t/2}^t\left(\int_{\R^d} p^{{\mathfrak{K}}_2}(t-s,x-y)\left(\int_{\R^d}\delta_{p^{{\mathfrak{K}}_1}}(s,y; z)({\mathfrak{K}}_1(z)-
{\mathfrak{K}}_2(z))J(z)dz\right)dy\right)ds.
\end{eqnarray*}
By using \eqref{e:fract-der-est1}, Proposition \ref{p:upperestonp}
and the convolution inequality \eqref{e:convolution-3}, we have
\begin{align*}
&\int_0^{t/2} \left(\int_{\R^d} p^{{\mathfrak{K}}_1}(s,y)\left(\int_{\R^d}\delta_{p^{{\mathfrak{K}}_2}}(t-s,x-y; z)({\mathfrak{K}}_1(z)-
{\mathfrak{K}}_2(z))J(z)dz\right)dy\right)ds\\
&+ \int_0^{t/2} \left(\int_{\R^d} p^{{\mathfrak{K}}_2}(s,x-y)\left(\int_{\R^d}\delta_{p^{{\mathfrak{K}}_1}}(t-s,y; z)({\mathfrak{K}}_1(z)-
{\mathfrak{K}}_2(z))J(z)dz\right)dy\right)ds\\
  \le & \wh{\gamma}_0 \|{\mathfrak{K}}_1-{\mathfrak{K}}_2\|_{\infty}
  \Bigg( \int_0^{t/2}\left(\int_{\R^d}p^{{\mathfrak{K}}_1}(s,y)\left(\int_{\R^d}\left|\delta_{p^{{\mathfrak{K}}_2}}(t-s,x-y; z)\right|j(|z|)dz\right)dy\right)ds\\
  &\quad +
    \int_0^{t/2}\left(\int_{\R^d}p^{{\mathfrak{K}}_2}(s,x-y)\left(\int_{\R^d}\left|\delta_{p^{{\mathfrak{K}}_1}}(t-s,y; z)\right|j(|z|)dz\right)dy\right)ds \Bigg)
  \\
\le & c_1 \|{\mathfrak{K}}_1-{\mathfrak{K}}_2\|_{\infty} \int_0^{t/2} \int_{\R^d} s \left(\rho(s,y)\rho(t-s,x-y)
+  \rho(s,x-y)\rho(t-s,y) \right)dy\, ds
\\
\le & 2 c_1\|{\mathfrak{K}}_1-{\mathfrak{K}}_2\| _{\infty} t^{-1}
\int_0^t \int_{\R^d}  s(t-s) \left(\rho(s,y) \rho(t-s,x-y) +\rho(s,x-y) \rho(t-s,y)\right) dy ds \\
\le & c_2\|{\mathfrak{K}}_1-{\mathfrak{K}}_2\| _{\infty}\ t \rho(t,x),  \qquad \text{for all } t \in (0, T],  \, x \in \R^d \, .
\end{align*}
\noindent 
(ii) 
Set
$
\wh{\mathfrak{K}}_i(z):=\mathfrak{K}_i(z)-{\kappa_0}/{2}, \quad i=1,2.
$
It is straightforward to see that $p^{\kappa_0/2}(t,x)=p^1(\kappa_0 t/2, x)$. Thus,  by the construction of 
the L\'evy process 
we have that for $i=1,2$, 
\begin{equation}\label{e:convolution-kappa}
p^{\mathfrak{K}_i}(t,x)=\int_{\R^d}p^{\kappa_0/2}(t, x-y)p^{{\wh{\mathfrak{K}}_i}}(t,y)\, dy=\int_{\R^d}p^1(\kappa_0 t/2, x-y)p^{{\wh{\mathfrak{K}}_i}}(t,y)\, dy.
\end{equation}

By \eqref{e:convolution-kappa}, Proposition \ref{p:upperestonp}, \eqref{e:difference-kappa}, \eqref{e:convolution-2} in the penultimate line (with $t,2t$ instead of $s,t$), and 
Lemma \ref{l:psi-and-f}(b)  in the last line, we have  that for all  $t \in (0, T]$ and   $x \in \R^d$,
\begin{align*}
&\left|\nabla p^{{\mathfrak{K}}_1}(t,x) -\nabla p^{{\mathfrak{K}}_2}(t,x)\right| =  
 \left|\int_{\R^d}\nabla p^1\left({\kappa_0 t}/{2}, x-y\right)
\left(p^{\wh{\mathfrak{K}}_1}(t,y)-p^{\wh{{\mathfrak{K}}}_2}(t,y)\right)dy\right|\\
 \le &  c_1\|{\mathfrak{K}}_1-{\mathfrak{K}}_2\|_{\infty} \Phi^{-1}(t^{-1})t^2 \int_{\R^d}\rho(t,x-y)
\rho(t,y)
\, dy\\
\le & c_2 \|{\mathfrak{K}}_1-{\mathfrak{K}}_2\|_{\infty} \Phi^{-1}(t^{-1})t
\rho(t,y)\,.
\end{align*}

\noindent 
(iii) By using \eqref{e:convolution-kappa}, \eqref{e:delta-p-kappa}, Lemma \ref{l:convolution}(b) and \eqref{e:difference-kappa}, we have
\begin{align*}
&\left|\delta_{p^{{\mathfrak{K}}_1}}(t,x;z)-\delta_{p^{{\mathfrak{K}}_2}}(t,x;z)\right|\\
 =&\left|\int_{\R^d}\delta_{p^1} \left({\kappa_0 t}/{2}, x-y; z\right) \left(p^{\wh{{\mathfrak{K}}}_1}(t,y)-p^{\wh{{\mathfrak{K}}}_2}(t,y)\right)dy\right| \\
\le &  c_1 \|{\mathfrak{K}}_1-{\mathfrak{K}}_2\|_{\infty} \left((\Phi^{-1}(t^{-1})|z|)^2\wedge 1\right) t^2
 \int_{\R^d}(\rho(t,x-y\pm z)+\rho(t,x-y))\rho(t,y)dy\\
\le & c_2 \|{\mathfrak{K}}_1-{\mathfrak{K}}_2\|_{\infty} 
\left((\Phi^{-1}(t^{-1})|z|)^2\wedge 1\right) t 
\left(\rho(t, x\pm z)+\rho(t, x) \right)\, .
\end{align*}
Now we have
\begin{align*}
&\int_{\R^d}\left|\delta_{p^{{\mathfrak{K}}_1}}(t,x;z)-\delta_{p^{{\mathfrak{K}}_2}}(t,x;z)\right| j(|z|)\, dz\, 
\le\,  c_2 \|{\mathfrak{K}}_1-{\mathfrak{K}}_2\|_{\infty} \\
&\quad \times \int_{\R^d}\left((\Phi^{-1}(t^{-1})|z|)^2\wedge 1\right)t\left(\rho(t, x\pm z)+\rho(t, x)    \right)\,  j(|z|)\,dz\\
= & c_2 \|{\mathfrak{K}}_1-{\mathfrak{K}}_2\|_{\infty} \int_{\R^d}\left((\Phi^{-1}(t^{-1})|z|)^2\wedge 1\right)t\left(\rho(t, x\pm z)+\rho(t, x)   \right)\,  j(|z|)\,dz, 
\end{align*}
which is the same as \eqref{e:fract-der-est3} and was estimated in the proof of Theorem \ref{t:fract-der-est} by $c_3\rho(t,x)$.
 This finishes the proof.
\qed


\section{Levi's construction of heat kernels}

For the remainder of this paper, we always assume that $\kappa:\R^d\times \R^d\to (0,\infty)$ is a Borel  function 
satisfying \eqref{e:intro-kappa} and \eqref{e:intro-kappa-holder}, that $\Phi$ satisfies \eqref{e:intro-wsc} and \eqref{e:intro-psibound} and 
that $J$ satisfies \eqref{e:psi1}.
Throughout the remaining part of this paper, 
$\beta$ is the constant  in \eqref{e:intro-kappa-holder}.

For a fixed $y\in \R^d$, let ${\mathfrak K}_y(z)=\kappa(y,z)$ and let $\LL^{{\mathfrak K}_y}$ be the freezing operator
\begin{align}
\label{e:Lkappapv}
\LL^{{\mathfrak K}_y}f(x)=
\LL^{{\mathfrak K}_y, 0}f(x)=
\lim_{\eps \downarrow 0}\LL^{{\mathfrak K}_y,  \eps}f(x), \,\,\,
\text{where }\LL^{{\mathfrak K}_y,  \eps}f(x)=\int_{|z| > \eps} \delta_f(x;z)  \kappa(y,z)J(z)dz.
\end{align}
Let $p_y(t,x)=p^{{\mathfrak K}_y}(t,x)$ be the heat kernel of the operator $\LL^{{\mathfrak K}_y}$.  
Note that $x\mapsto p_y(t,x)$ is in $C^\infty_0(\R^d)$ and satisfies \eqref{e:prop-p-L}.

\subsection{Estimates on $p_y(t,x-y)$}

The following result is the counterpart of \cite[Lemmas 3.2 and 3.3]{CZ}.

\begin{lemma}\label{l:some-estimates-1}
For every $T\ge 1$ and $\beta_1 \in (0, \delta_1) \cap (0, \beta]$, there exists a constant 
$c=c(T, d$, $\delta_1, \beta_1,\kappa_0$,$\kappa_1, \kappa_2, \gamma_0)>0$
such that 
for all $x\in \R^d$ and $t \in (0,T]$,
\begin{equation}\label{e:some-estimates-2a}
\left|\int_{\R^d}\LL^{{\mathfrak K}_x, \eps} p_y(t,\cdot)(x-y)\, dy\right|\,\le \,c \,t^{-1}\,\Phi^{-1}(t^{-1})^{-{\beta_1}}\, , \text{ for all }\eps \in [0,1], 
\end{equation}
\begin{equation}\label{e:some-estimates-2b}
\left|\int_{\R^d}\partial_t p_y(t,x-y)\, dy\right|\,\le \,c\, t^{-1}\,\Phi^{-1}(t^{-1})^{-{\beta_1}},
\end{equation}
\begin{equation}\label{e:reinstated-estimate}
\left|\int_{\R^d}\nabla p_y(t, \cdot)(x-y)\, dy\right|\le c\,  \Phi^{-1}(t^{-1})^{1-{\beta_1}}\, .
\end{equation}
Furthermore, we have
\begin{equation}\label{e:some-estimates-2c}
\lim_{t\downarrow 0}\sup_{x\in \R^d} \left|\int_{\R^d} p_y(t,x-y)\, dy-1\right|=0\, .
\end{equation}
\end{lemma}
\pf 
Choose  $\gamma \in (0, \delta_1-{\beta_1})\cap (0, 1]$.
Since $\int_{\R^d}p_z(t,\xi-y)dy=1$ for every $\xi, z \in \R^d$, by the definition of $\delta_{p_x}$ we have $\int_{\R^d}\delta_{p_x}(t,x-y;w)dy=0$. Therefore, 
using this, \eqref{e:intro-kappa}, \eqref{e:psi1} and \eqref{e:delta-difference-abs},
for $ \eps \in [0,1]$ and $t\in (0, T]$,
\begin{eqnarray*}
\lefteqn{\left|\int_{\R^d}\LL^{{\mathfrak K}_x, \eps} p_y(t,\cdot)(x-y)\, dy\right|}\\
&=&\left|\int_{\R^d}\left(\int_{|w|>\eps}\left(\delta_{p_y}(t,x-y;w)-\delta_{p_x}(t,x-y;w)\right)\kappa(x,w)J(w)\, dw\right)dy\right|\\
&\le &\kappa_1 \gamma_0\int_{\R^d}\left(\int_{|w|>\eps}\left|\delta_{p_y}(t,x-y;w)-\delta_{p_x}(t,x-y;w)\right|j(|w|)\, dw\right)dy\\
&\le & c_1 \int_{\R^d} \|\kappa(y,\cdot)-\kappa(x, \cdot)\|_{\infty}
\rho(t,x-y)\, dy\\
&\le & c_1 \kappa_2 \int_{\R^d} \left(|x-y|^{{\beta_1}}\wedge 1\right)
 \rho(t,x-y)\, dy 
\le c_2 t^{-1}
\Phi^{-1}(t^{-1})^{-{\beta_1}}.
\end{eqnarray*}
Here the last line follows from \eqref{e:intro-kappa-holder} and  \eqref{e:convolution-integrability} since  ${\beta_1}+\gamma \in (0, \delta_1)$.

For \eqref{e:some-estimates-2b}, by using 
\eqref{e:fract-der-est1} and
\eqref{e:some-estimates-2a} in the third line, we get,  for $t \in (0,T]$,
\begin{eqnarray*}
\lefteqn{\left|\int_{\R^d}\partial_t p_y(t,x-y)\, dy\right|=\left|\int_{\R^d}\LL^{{\mathfrak K}_y}p_y(t, \cdot)(x-y)\, dy\right|}\\
&\le &\left|\int_{\R^d}\left(\LL^{{\mathfrak K}_x}-\LL^{{\mathfrak K}_y}\right)p_y(t,\cdot)(x-y)\, dy\right| +\left|\int_{\R^d}\LL^{{\mathfrak K}_x}p_y(t, \cdot)(x-y)\, dy\right| \\
&\le & c_3 \int_{\R^d}\rho_0^{{\beta_1}}(t,x-y)\, dy +c_2t^{-1}\Phi^{-1}(t^{-1})^{-{\beta_1}}\,\le \, c_4 t^{-1}\Phi^{-1}(t^{-1})^{-{\beta_1}}\, .
\end{eqnarray*}
Here we have used \eqref{e:convolution-integrability} in the last inequality.

For \eqref{e:reinstated-estimate}, by \eqref{e:difference-nkappa} we have
\begin{align*}\
&\left|\int_{\R^d}\nabla p_y(t, \cdot)(x-y)\, dy\right|=\left|\int_{\R^d}\left(\nabla p_y(t, \cdot)-\nabla p_x(t,\cdot)\right)(x-y)\, dy\right|\\
& \le c_5  \int_{\R^d}\|\kappa(x,\cdot)-\kappa(y,\cdot)\|_{\infty}t\Phi^{-1}(t^{-1}) \rho(t,x-y)\, dy\\
&\le c_6  \int_{\R^d} \left(|x-y|^{{\beta_1}}\wedge 1\right)t\Phi^{-1}(t^{-1}) \rho(t,x-y)\, dy\\
&=t\Phi^{-1}(t^{-1})\int_{\R^d}  \rho_0^{\beta_1}(t,x-y)\, dy\\
&\le c_7  t\Phi^{-1}(t^{-1}) t^{-1}\Phi^{-1}(t^{-1})^{-{\beta_1}}= \Phi^{-1}(t^{-1})^{1-{\beta_1}}\, .
\end{align*}
In the last inequality we used Lemma \ref{l:convolution}(a) which requires that ${\beta_1}+\gamma \in (0, \delta_1)$.

Finally, by using \eqref{e:difference-kappa} in the second line and \eqref{e:convolution-integrability} in the last inequality, we get
\begin{align*}
&\sup_{x\in \R^d}\left|\int_{\R^d}p_y(t,x-y)\, dy-1\right| \le \sup_{x\in \R^d}\int_{\R^d}\left|p_y(t,x-y)-p_x(t,x-y)\right|dy \\
&\le c_8  \sup_{x\in \R^d}\int_{\R^d} \|\kappa(y, \cdot)-\kappa(x,\cdot)\|_{\infty} t
 \rho(t,x-y)\, dy\\
&\le c_9t \sup_{x\in \R^d}\int_{\R^d} 
 \rho_0^{\beta_1}(t,x-y)\, dy 
\le c_{10}  \Phi^{-1}(t^{-1})^{-{\beta_1}},  \quad t \in (0,T]\, .
\end{align*}
\qed

\begin{lemma}\label{l:jcontoffzkernel}
The function $p_y(t, x)$ is jointly continuous in $(t, x, y)$.
\end{lemma}

\pf
By the triangle inequality, we have
$$
|p_{y_1}(t_1,x_1)-p_{y_2}(t_2,x_2)|\le |p_{y_1}(t_1,x_1)-p_{y_2}(t_1,x_1)|+|p_{y_2}(t_1,x_1)-p_{y_2}(t_2,x_2)|.
$$
Applying \eqref{e:difference-kappa} and \eqref{e:intro-kappa-holder} to the first term on the right hand side and 
Lemma \ref{l:partial-time}(b) to the second term on the right hand side, we immediately get the desired joint continuity.
\qed

\subsection{Construction of $q(t,x,y)$}
For $(t,x,y)\in (0,\infty)\times \R^d\times \R^d$ define
\begin{align}
q_0(t,x,y)&:=\frac12 \int_{\R^d}\delta_{p_y}(t,x-y;z)
\left(\kappa(x,z)-\kappa(y,z)\right)J(z)\, dz\nn\\
&= \big(\LL^{{\mathfrak K}_x}-\LL^{{\mathfrak K}_y}\big)p_y(t, \cdot)(x-y)\, .
\label{e:q0-definition}
\end{align}
In the next lemma we collect several estimates on $q_0$ that will be needed later on.

\begin{lemma}\label{l:estimates-q0} For every $T\ge 1$ and $\beta_0 \in  (0, \beta]$, there is a constant 
$C_1\ge 1$ depending on $d,\delta_1, \kappa_0,\kappa_1, \kappa_2, \gamma, T$ and $\Phi^{-1}(T^{-1})$
 such that
for $t \in (0, T]$ and 
 $x, x', y, y'\in \R^d$,
\begin{eqnarray}
|q_0(t,x,y)| \le  C_1 (|x-y|^{\beta_0}\wedge 1)\rho(t,x-y)=C_1\rho_0^{\beta_0}(t,x-y),  \label{e:q0-estimate}
\end{eqnarray}
and for all $\gamma\in (0,\beta_0)$, 
\begin{eqnarray}
&&|q_0(t,x,y)-q_0(t,x',y)|\nn\\
&&\le C_1 \left(|x-x'|^{\beta_0-\gamma}\wedge 1\right)\left\{\left(\rho_{\gamma}^0+\rho_{\gamma-\beta_0}^{\beta_0}\right)(t,x-y)
+\left(\rho_{\gamma}^0+\rho_{\gamma-\beta_0}^{\beta_0}\right)(t,x'-y)\right\}\label{e:estimate-step3}
\end{eqnarray}
and 
\begin{equation}
|q_0(t,x,y)-q_0(t,x,y')| \le C_1 \Phi^{-1}(t^{-1})^{\beta_0} \left(|y-y'|^{\beta_0}\wedge 1\right)
\left(\rho(t,x-y)+ \rho(t,x-y')\right).
\label{e:estimate-q0-2}
\end{equation}
\end{lemma}
\pf
(a) \eqref{e:q0-estimate} follows from  \eqref{e:fract-der-est1} and \eqref{e:intro-kappa-holder}.

\noindent (b) By \eqref{e:q0-estimate} and \eqref{e:nonincrease-gamma}, we have that for $t \in (0, T]$,
$$
|q_0(t,x,y)| \le  c_0\rho_0^{\beta_0}(t,x-y)\le c_0 \Phi^{-1}(T^{-1})^{\gamma-\beta_0} \rho_{\gamma-\beta_0}^{\beta_0}(t,x-y), 
$$
which proves \eqref{e:estimate-step3} for $|x-x'| \ge 1 $. 
Now suppose that $1\ge |x-x'|\ge \Phi^{-1}(t^{-1})^{-1}$. Then, by \eqref{e:q0-estimate}, for $t \in (0, T]$,
$$
|q_0(t,x,y)|\le c_1 \left(\Phi^{-1}(t^{-1})\right)^{-(\beta_0-\gamma)}\rho_{\gamma-\beta_0}^{\beta_0}(t,x-y)\le c_1 |x-x'|^{\beta_0-\gamma}\rho_{\gamma-\beta_0}^{\beta_0}(t,x-y),
$$
and the same estimate is valid for $|q_0(t,x',y)|$. By adding we get \eqref{e:estimate-step3} for this case. Finally, assume that $|x-x'|\le  1 \wedge \Phi^{-1}(t^{-1})^{-1}$. Then, by \eqref{e:psi1},  \eqref{e:intro-kappa-holder} and  \eqref{e:fract-der-est2}, for $t \in (0, T]$,
\begin{align*}
\lefteqn{|q_0(t,x,y)-q_0(t,x',y)|=\left|\int_{\R^d}\delta_{p_y}(t,x-y;z)(\kappa(x,z)-\kappa(y,z))J(z)\, dz\right.}\\
&  \quad -\left .\int_{\R^d}\delta_{p_y}(t,x'-y;z)(\kappa(x',z)-\kappa(y,z))J(z)\, dz\right|\\
\le & \gamma_0\int_{\R^d}|\delta_{p_y}(t,x-y;z)-\delta_{p_y}(t,x'-y;z)|\ |\kappa(x,z)-\kappa(y,z)|j(|z|)\, dz\\
 & \quad + \gamma_0\int_{\R^d}|\delta_{p_y}(t,x'-y;z)||\kappa(x,z)-\kappa(x',z)|j(|z|)\, dz\\
\le & \gamma_0 \kappa_2 \left(|x-y|^{\beta_0}\wedge 1\right)\int_{\R^d}|\delta_{p_y}(t,x-y;z)-\delta_{p_y}(t,x'-y;z)| j(|z|)\, dz\\
 & \quad  + \gamma_0\kappa_2\left(|x-x'|^{\beta_0}\wedge 1\right)\int_{\R^d}|\delta_{p_y}(t,x'-y;z) |j(|z|)\, dz\\
\le & c_2 \left(|x-y|^{\beta_0}\wedge 1\right) \big(\rho(t,x-y)+\rho(t,x'-y)\big)+ c_2 |x-x'|^{\beta_0} \rho(t,x'-y).
\end{align*}
By using the definition of $\rho(t,x'-y)$, the obvious equality $x'-y=(x-y)+(x'-x)$, the assumption that $|x-x'|\le \Phi^{-1}(t^{-1})^{-1}$ and \eqref{e:2.9}, we conclude that $\rho^\beta_0(t,x'-y)\le 4\rho^\beta_0(t,x-y)$. Thus, it follows that for $t \in (0, T]$,
\begin{align*}
|q_0(t,x,y)-q_0(t,x',y)|
\,\le \,&\,5\,c_2\, \rho_0^{\beta_0}(t,x-y)+c_2 |x-x'|^{\beta_0} \rho(t,x'-y)\\
\,\le \,&\,5\,c_2 \, |x-x'|^{\beta_0-\gamma}\rho_{\gamma-\beta_0}^{\beta_0}(t,x-y)+c_2 |x-x'|^{\beta_0-\gamma}\rho_{\gamma}^{0}(t,x'-y)\, .
\end{align*}

\noindent (c) First note that
\begin{align*}
&q_0(t,x,y)-q_0(t,x,y')\\
=&\frac12\int_{\R^d} \delta_{p^y}(t,x-y;z)\left(\kappa(y',z)-\kappa(y,z)\right)J(z)\, dz \\
& \ \ \ +\frac12\int_{\R^d}\left(\delta_{p^y}(t,x-y;z)-\delta_{p^y}(t,x-y';z)\right)\left(\kappa(x,z)-\kappa(y',z)\right)J(z)\, dz\\
&\ \ \ +\frac12\int_{\R^d}\left(\delta_{p^y}(t,x-y';z)-\delta_{p^{y'}}(t,x-y';z)\right)\left(\kappa(x,z)-\kappa(y',z)\right)J(z)\, dz\\
=:& I_1+I_2+I_3\, .
\end{align*}
It follows from \eqref{e:intro-kappa-holder}, \eqref{e:psi1} and \eqref{e:fract-der-est1} that for $t\in (0,T]$,
\begin{align*}
| I_1| \le c_1\left( |y-y'|^{\beta_0}\wedge 1\right) \int_{R^d} \left|\delta_{p^y}(t,x-y;z)\right| j(|z|)\, dz \le c_2 \left( |y-y'|^{\beta_0}\wedge 1\right)\rho(t,x-y)\, ,
\end{align*}
which is smaller than or equal to the right-hand side in \eqref{e:estimate-q0-2} since $\Phi^{-1}(t^{-1})\ge \Phi^{-1}(T^{-1})$. By \eqref{e:intro-kappa}, \eqref{e:psi1} and \eqref{e:fract-der-est2} we get that
\begin{align*}
| I_2| \,\le\, & \,c_1\,\int_{\R^d}\left|\delta_{p^y}(t,x-y;z)-\delta_{p^y}(t,x-y';z)\right| j(|z|)\, dz\\
\,\le\, & \,c_2 \,  \left((\Phi^{-1}(t^{-1})|y-y'|)\wedge 1\right) \left(\rho(t,x-y)+\rho(t,x-y')\right)\\
\,\le\, & \,c_2 \, \Phi^{-1}(T^{-1})^{-\beta_0} \Phi^{-1}(t^{-1})^{\beta_0}\left( |y-y'|^{\beta_0}\wedge 1\right)  \left(\rho(t,x-y)+\rho(t,x-y')\right).
\end{align*}
Finally, by \eqref{e:intro-kappa}, \eqref{e:intro-kappa-holder}, \eqref{e:psi1} and \eqref{e:delta-difference-abs},  for $t\in (0,T]$,
\begin{align*}
| I_3| \, \le\, & c_1\,\int_{\R^d}\left|\delta_{p^y}(t,x-y';z)-\delta_{p^{y'}}(t,x-y';z)\right| j(|z|)\, dz\\
\,\le\,& c_3\, \|\kappa(y, \cdot)-\kappa(y', \cdot)\|_{\infty} 
 \rho(t,x-y')
\le \, c_4  \left( 
|y-y'|^{\beta_0}\wedge 1\right)   \rho(t,x-y')\, .
\end{align*}
\qed

\begin{lemma}\label{l:joint-cont-q0}
The function $q_0(t,x,y)$ is jointly continuous in $(t,x,y)$.
\end{lemma}
\pf 
 It follows from Lemma \ref{l:jcontoffzkernel} that $(t,x,y)\mapsto p_y(t,x-y)$ is jointly continuous and hence also that $\delta_{p_y}(t,x-y; z)$ is jointly continuous in $(t,x,y)$.  To prove the joint continuity of $q_0(t,x,y)$,  let $(t_n,x_n,y_n) \to (t,x,y)\in (0,T]\times \R^d\times \R^d$ and assume that $0<\eps \le t_n\le T$. The integrands will converge because of the joint continuity of $\delta_{p_y}$ and continuity of $\kappa$ in the first variable. Moreover, by \eqref{e:delta-p-kappa},
\begin{align*}
\lefteqn{\left|\delta_{p_{y_n}}(t_n,x_n-y_n;z)\right| |\kappa(x_n,z)-\kappa(y_n,z)| j(|z|) }\\
&\le  c_1\left( (\Phi^{-1}(t_n^{-1})|z|^2)\wedge 1\right) T\left(\rho(t_n,x_n-y_n\pm z)+\rho(t_n,x_n, y_n)\right) j(|z|)\\
&\le c_2 \rho(\eps, 0) \left((\Phi^{-1}(\eps^{-1})|z|^2)\wedge 1\right) j(|z|).
\end{align*}
Since the right-hand side is integrable on $\R^d$, the joint continuity follows by use of the dominated convergence theorem.
\qed

For $n\in \N$, we inductively define
\begin{equation}\label{e:qn-definition}
q_n(t,x,y):=\int_0^t \int_{\R^d}q_0(t-s,x,z)q_{n-1}(s,z,y)\, dz\, ds, \quad (t,x,y)\in (0, \infty)\times \R^d \times \R^d\, .
\end{equation}
The following result is the counterpart of \cite[Theorem 3.1]{CZ}.

\begin{thm}\label{t:definition-of-q}
The series $q(t,x,y):=\sum_{n=0}^{\infty}q_n(t,x,y)$ is absolutely and locally uniformly convergent on $(0, \infty)\times \R^d \times \R^d$ and solves the integral equation
\begin{equation}\label{e:integral-equation}
q(t,x,y)=q_0(t,x,y)+\int_0^t \int_{\R^d}q_0(t-s,x,z)q(s,z,y)\, dz\, ds\, .
\end{equation}
Moreover, $q(t,x,y)$ is jointly continuous in $(t,x,y)\in (0, \infty)\times \R^d \times \R^d$ and has the following estimates: for every $T\ge 1$ and $\beta_2\in  (0, \beta] \cap (0, \delta_1/2)$  there is a constant 
$C_2=C_2(T, d,\delta_1, \kappa_0, \kappa_1, \kappa_2, \beta_2, \gamma_0)>0$
such that
\begin{equation}\label{e:q-estimate}
|q(t,x,y)|\le C_2 \big(\rho_0^{{\beta_2}}+\rho_{{\beta_2}}^0\big)(t,x-y), \quad (t,x,y)\in (0, T]\times \R^d \times \R^d\, ,
\end{equation}
and for any $\gamma\in (0,{\beta_2})$ and $T \ge 1$ there is a constant $C_3=C_3(T, d,\delta_1, \gamma,\kappa_0, \kappa_1, \kappa_2, \gamma_0, {\beta_2})>0$ such that
for all $(0, T]\times \R^d \times \R^d$,
\begin{eqnarray}
&&|q(t,x,y)-q(t,x',y)|\nn\\
&&\le C_3\left(|x-x'|^{{\beta_2}-\gamma}\wedge 1\right)
\left(\big(\rho_{\gamma}^0+\rho_{\gamma-{\beta_2}}^{{\beta_2}}\big)(t,x-y)+\big(\rho_{\gamma}^0+\rho_{\gamma-{\beta_2}}^{{\beta_2}}\big)(t,x'-y)\right)\, .
\label{e:difference-q-estimate}
\end{eqnarray}
\end{thm}

\pf 
This proof follows  the main idea of the proof of  \cite[Theorem 3.1]{CZ}, except that
we give a full proof of the joint continuity 
in Step 2. 
We give the details for the readers' convenience. 
In this proof, $T\ge 1$ is arbitrary.

\noindent {\it Step 1:} 
By \eqref{e:q0-estimate}, 
\eqref{e:convolution-3} and \eqref{e:convolution-4}, 
we have that
\begin{align*}
|q_1(t,x,y)|\,\le\, &C_1^2  \int_0^t \int_{\R^d} \rho_0^{{\beta_2}}(t-s,x-y-u)\rho_0^{{\beta_2}}(s,u)\, 
du \, ds\\
\,\le\,  &8C_0C_1^2B\left({{\beta_2}}/{2}, {{\beta_2}}/{2}\right)\big(\rho_{2{\beta_2}}^0+\rho_{{\beta_2}}^{{\beta_2}}\big)(t,x-y), 
 \quad t \le T\, .
\end{align*}
Let $C=2^4C_0C_1^2  $ and we claim that for $n\ge 1$ and $t \le T$, 
\begin{equation}\label{e:bdonq_n}
|q_n(t,x,y)|\le \gamma_n \big(\rho_{(n+1){\beta_2}}^0+\rho_{n{\beta_2}}^{{\beta_2}}\big)(t,x-y)
\end{equation}
with 
$$
\gamma_n= C^{n+1} \prod_{j=1}^n  B\left({{\beta_2}}/{2}, {j{\beta_2}}/{2}\right)\, .
$$
We have seen that \eqref{e:bdonq_n} is valid for $n=1$. Suppose that it is valid for $n$.
Then by using \eqref{e:convolution-3}, \eqref{e:convolution-4}, 
\eqref{e:nonincrease-gamma} and \eqref{e:nonincrease-beta},  we have
that for $t \le T$,
\begin{align*}
|q_{n+1}(t,x,y)|\le & \int_0^t \int_{\R^d}|q_0(t-s,x,z)|\ |q_n(s,z,y)|\, dz\, ds\\
\le & C_1\gamma_n \int_0^t \int_{\R^d}\rho_0^{{\beta_2}}(t-s,x-z)\big(\rho_{(n+1){\beta_2}}^0+\rho_{n{\beta_2}}^{{\beta_2}}\big)(s,z-y)\, dz\, ds\\
\le & 2^4 C_0 C_1\gamma_nB\left(\frac{{\beta_2}}{2}, \frac{(n+1){\beta_2}}{2}\right)\Big(\rho_{(n+2){\beta_2}}^0 +\rho_{(n+1){\beta_2}}^{{\beta_2}}\Big)(t,x-y)\\
 \le &\gamma_{n+1}\Big(\rho_{(n+2){\beta_2}}^0 +\rho_{(n+1){\beta_2}}^{{\beta_2}}\Big)(t,x-y)\, .
\end{align*}
Thus \eqref{e:bdonq_n} is valid.
Since
$$
\sum_{n=0}^{\infty}\gamma_n \Phi^{-1}(T^{-1})^{-(n+1)
{\beta_2}}= 
\sum_{n=0}^{\infty}\frac{\left(\Phi^{-1}(T^{-1})^{-{\beta_2}} 
C\Gamma\left(\frac{{\beta_2}}{2}\right)\right)^{n+1}}{\Gamma\left(\frac{(n+1){\beta_2}}{2}\right)}=:C_2<\infty \,,
$$
by using \eqref{e:nonincrease-gamma} and \eqref{e:nonincrease-beta} in the second line, it follows that for $t \le T$,
\begin{align*}
&\sum_{n=0}^{\infty}|q_n(t,x,y)|\,\le\,  \sum_{n=0}^{\infty}\gamma_n \big(\rho_{(n+1){\beta_2}}^0+
\rho^{\beta_2}_{n{\beta_2}}\big)(t,x-y)\\
&\le   \,
\sum_{n=0}^{\infty}\gamma_n \Phi^{-1}(T^{-1})^{-(n+1){\beta_2}} \big(\rho_{{\beta_2}}^0+\rho_0^{{\beta_2}}\big)(t,x-y)
=C_2 \big(\rho_{{\beta_2}}^0+\rho_0^{{\beta_2}}\big)(t,x-y)\, .
\end{align*}
This proves that $\sum_{n=0}^{\infty}q_n(t,x,y)$ is
 absolutely and uniformly convergent on $[\eps,T]\times \R^d \times \R^d$ for all $\eps \in (0,1)$ and $T \ge 1$, hence $q(t,x,y)$ is well defined. Further, by \eqref{e:qn-definition},
$$
\sum_{n=0}^{m+1}q_n(t,x,y)=q_0(t,x,y)+\int_0^t \int_{\R^d} q_0(t-s,x,z)\sum_{n=0}^m q_n(s,z,y)\, dz\, ds\, ,
$$
and \eqref{e:integral-equation} follows by taking the limit of both sides as $m\to \infty$.

\noindent {\it Step 2:}  The joint continuity of $q_0(t,x,y)$ was shown in Lemma \ref{l:joint-cont-q0}.
We now prove the joint continuity of $q_1(t, x, y)$.  For any $x, y\in \R^d$ and $t, h>0$, we have
\begin{eqnarray}
&&q_1(t+h, x, y)-q_1(t, x, y)\nn\\
&&=\int^{t+h}_t\int_{\R^d}q_0(t+h-s, x, z)q_0(s, z, y)dzds
\nn\\
&&\ \ \ +\int^t_0\int_{\R^d}\left(q_0(t+h-s, x, z)-q_0(t-s, x, z) \right)q_0(s, z, y)dzds.
\label{e:rscontofq-1int}
\end{eqnarray}
It follows from \eqref{e:q0-estimate} that, there exists $c_1=c_1(T)>0$ such that,  for  $0<h\le t/4$ and $t+h \le  T$,
\begin{eqnarray*}
&&\sup_{x,y \in \R^d}|\int^{t+h}_t\int_{\R^d}q_0(t+h-s, x, z)q_0(s, z, y)dzds|\\
&&\le c_1\sup_{x,y \in \R^d}\int^{t+h}_t\int_{\R^d}\rho^{\beta_2}_0(t+h-s, x-z)\rho^{\beta_2}_0(s, z-y)dzds\\
&&=c_1\sup_{x,y \in \R^d}\int^h_0\int_{\R^d}\rho^{\beta_2}_0(r, x-z)\rho^{\beta_2}_0(t+h-r, z-y)dzdr\\
&&\le  c_1 \int^h_0\sup_{x,y \in \R^d}\int_{\R^d}\rho^{\beta_2}_0(r, x-z)\rho^{\beta_2}_0(t-r, z-y)dzdr.
\end{eqnarray*}
Now applying Lemma \ref{l:convolution}(b),  we get
\begin{eqnarray*}
&&\sup_{x,y \in \R^d}\int_{\R^d}\rho^{\beta_2}_0(r, x-z)\rho^{\beta_2}_0(t-r, z-y)dz\\
&&\le c_2((t-r)^{-1}\Phi^{-1}((t-r)^{-1})^{-{\beta_2}}+r^{-1}\Phi^{-1}(r^{-1})^{-{\beta_2}})\rho(t, 0).
\end{eqnarray*}
It follows from Lemma \ref{l:convoluton-inequality} that the 
right-hand 
side of the display above is integrable in $r\in (0, t)$, so by the dominated convergence theorem,
we get
\begin{equation}\label{e:rscontofq-1int-2}
\lim_{h\downarrow0}\sup_{x, y\in \R^d}|\int^{t+h}_t\int_{\R^d}q_0(t+h-s, x, z)q_0(s, z, y)dzds|=0.
\end{equation}
Using \eqref{e:q0-estimate} again, we get that for $s\in (0, t]$,
\begin{eqnarray*}
&&|\left(q_0(t+h-s, x, z)-q_0(t-s, x, z) \right)q_0(s, z, y)|\\
&&\le c_3\left(\rho^{\beta_2}_0(t+h-s, x-z)+\rho^{\beta_2}_0(t-s, x-z) \right)\rho^{\beta_2}_0(s, z-y)\\
&&\le c_4 \rho^{\beta_2}_0(t-s, x-z)\rho^{\beta_2}_0(s, z-y).
\end{eqnarray*}
It follows from Lemma \ref{l:convolution}(c) that
$$
\int^t_0\int_{\R^d}\rho^{\beta_2}_0(t-s, x-z)\rho^{\beta_2}_0(s, z-y)dzds  \le 
c_5 (\rho^{2{\beta_2}}_0 (t, 0) +\rho^{\beta_2}_{\beta_2} (t, 0))<\infty,
$$
thus we can use the dominated convergence theorem to get that, by the continuity of $q_0$,
\begin{equation}\label{e:rscontofq-1int-3}
\lim_{h\downarrow0}\int^t_0\int_{\R^d}\left(q_0(t+h-s, x, z)-q_0(t-s, x, z) \right)q_0(s, z, y)dzds=0.
\end{equation}

It follows from \eqref{e:estimate-q0-2} that for $s\in (0, T]$,
\begin{align*}
&|q_0(s, z, y)-q_0(s, z, y')|\\
\le& c_6\left((\Phi^{-1}(s^{-1})|y-y'|)^{\beta_2}\wedge1 \right)\left(\rho(s, z-y)+
\rho(s, z-y')\right).
\end{align*}
Now we fix $0<t_1\le t_2\le T$.
Then for any $\eps\in (0, t_1/4)$, $t\in [t_1, t_2]$ and $s\in [\eps, t]$,
\begin{eqnarray*}
&&|q_0(t-s, x, z)\left( q_0(s, z, y)-q_0(s, z, y')\right)|\\
&&\le \ \ c_7\left((\Phi^{-1}(\eps^{-1})|y-y'|)^{\beta_2}\wedge1 \right)\rho^{\beta_2}_0(t-s, x, z)\left(\rho(s, z-y)+
\rho(s, z-y')\right).
\end{eqnarray*}
By Lemma \ref{l:convolution}(c), we have
$$
\sup_{x, y, y'\in \R^d, t\in [t_1, t_2]}\int^t_0\int_{\R^d}\rho^{\beta_2}_0(t-s, x, z)\left(\rho(s, z-y)+
\rho(s, z-y')\right)dzds<\infty.
$$
Thus
$$
\lim_{y'\to y}\sup_{x\in \R^d, t\in [t_1, t_2]}\int^t_\eps \int_{\R^d}|q_0(t-s, x, z)\left( q_0(s, z, y)-q_0(s, z, y')\right)|dzds=0.
$$
Consequently, for  each 
$0 <t_1 < t_2 \le T$ and $\eps\in (0, t_1/4)$,
the family of functions 
$$
\left\{\int^t_\eps \int_{\R^d}q_0(t-s, x, z) q_0(s, z, \cdot)dzds: x\in \R^d, t\in [t_1, t_2]\right\}
$$
is equi-continuous. 
By combining \eqref{e:q0-estimate} and
Lemma \ref{l:rsnewlemma}, we get that
\begin{equation}\label{e:rsnew}
\lim_{\eps\to 0}\sup_{x, y\in \R^d, t\in [t_1, t_2]}
\left(\int^\eps_0+\int^t_{t-\eps}\right)\int_{\R^d}q_0(t-s, x, z)q_0(s, z, y)dzds=0.
\end{equation}
Therefore the family
\begin{equation}\label{e:family1}
\left\{\int^t_0 \int_{\R^d}q_0(t-s, x, z) q_0(s, z, \cdot)dzds: x\in \R^d, t\in [t_1, t_2]\right\}
\end{equation}
is equi-continuous.

Similarly, by using \eqref{e:estimate-step3}, we can show that, 
for each $0 <t_1 < t_2 \le T$ and $\eps\in (0, t_1/4)$,
the family of functions 
$$
\left\{\int^{t-\eps}_0
\int_{\R^d} q_0(t-s, \cdot, z) q_0(s, z, y)dzds: y\in \R^d, t\in [t_1, t_2]\right\}
$$
is equi-continuous. Combining this with \eqref{e:rsnew}, we get
the family of functions 
\begin{equation}\label{e:family2}
\left\{\int^{t}_0\int_{\R^d} q_0(t-s, \cdot, z) q_0(s, z, y)dzds: y\in \R^d, t\in [t_1, t_2]\right\}
\end{equation}
is equi-continuous. 

Now combining the continuity of $t \to q_1(t, x,y)$ (by \eqref{e:rscontofq-1int-2} and \eqref{e:rscontofq-1int-3}) and the equi-continuities of the families \eqref{e:family1} and \eqref{e:family2}, we immediately
get the joint continuity of $q_1$.

The joint continuity of $q_n(t, x, y)$ can be proved by induction
by using the estimate \eqref{e:bdonq_n} of $q_n$
and Lemma \ref{l:rsnewlemma}.
The joint continuity
of $q(t, x, y)$ follows immediately.

\noindent 
{\it Step 3:} 
By replacing $\alpha$ by $2$ and $\beta$ by $\beta_2$, this step is exactly the same as Step 4 in \cite{CZ}.
\qed

\subsection{Properties of $\phi_y(t,x)$}
Let
\begin{equation}\label{e:phi-y-def}
\phi_y(t,x,s):=\int_{\R^d}p_z(t-s,x-z)q(s,z,y)\, dz, \quad x \in \R^d, \, 
0< s<t
\end{equation}
and
\begin{equation}\label{e:def-phi-y-2}
\phi_y(t,x):=\int_0^t \phi_y(t,x,s)\, ds =\int_0^t \int_{\R^d}p_z(t-s,x-z)q(s,z,y)\, dz\, ds\, .
\end{equation}

The following result is the counterpart of \cite[Lemma 3.5]{CZ}.

\begin{lemma}\label{l:phi-y-abs-cont}
For all $x,y\in \R^d$, $x\neq y$, the mapping $t\mapsto \phi_y(t,x)$ is absolutely continuous  on $(0, \infty)$ and
\begin{equation}\label{e:phi-y-partial}
\partial_t \phi_y(t,x)=q(t,x,y)+\int_0^t \int_{\R^d} \LL^{\mathfrak{K}_z}p_z(t-s, \cdot)(x-z)q(s,z,y)\, dz\, ds, \quad t \in (0, \infty)\, .
\end{equation}
\end{lemma}
\pf {\it Step 1:} Here we prove that for 
any $T\ge 1$, $t \in (0, T]$ and $s\in (0,t)$,
\begin{equation}\label{e:phi-y-partial-2}
\partial_t \phi_y(t,x,s)=\int_{\R^d}\partial_t p_z(t-s,x-z)q(s,z,y)\, dz\, .
\end{equation}
Let $|\eps|<(t-s)/2$. We have that
\begin{eqnarray*}
\frac{\phi_y(t+\eps,x,s)-\phi_y(t,x,s)}{\eps}
=\int_{\R^d}\left(\int_0^1 \partial_t p_z(t+\theta \eps -s,x,z)\, d\theta\right) q(s,z,y)\, dz\, .
\end{eqnarray*}
By using \eqref{e:psi1}, \eqref{e:prop-p-L}, \eqref{e:fract-der-est1} and \eqref{e:representaion-L}, we have, 
\begin{eqnarray*}
\lefteqn{\left|\partial_t p_z(t+\theta \eps-s,x-z)\right|=\left|\LL^{\mathfrak{K}_z}p_z(t+\theta \eps-s, \cdot)(x-z)\right|} \\
& \le &\frac12 \gamma_0 \int_{\R^d}|\delta_{p_z}(t+\theta \eps-s, x-z;w)|\kappa(z,w)j(|w|)\, dw\\
&\le & c_1 \rho(t+\theta \eps -s, x-z) \, \le \,  c_2 \rho(t-s, x-z)\, .
\end{eqnarray*}
In the last inequality
we used that $|\eps|<(t-s)/2$ and applied 
Lemma \ref{l:psi-and-f}(b). Together with \eqref{e:q-estimate} this gives that for any $\beta_2\in (0, \beta)\cap (0, \delta_1/2)$ and $t\in (0, T]$
$$
\left|\partial_t p_z(t+\theta \eps -s, x-z)q(s,z,y)\right|\le c_3(T) \rho(t-s,x-z)\big(\rho_{\beta_2}^0+\rho_0^{\beta_2}\big)(s,z-y)=:g(z)\, .
$$
By \eqref{e:convolution-2}, we see that $\int_{\R^d}g(z)\, dz<\infty$. Thus, by the dominated convergence theorem,
$$
\lim_{\eps\to 0} \frac{\phi_y(t+\eps,x,s)-\phi_y(t,x,s)}{\eps}=\int_{\R^d}\partial_t p_z(t-s,x-z)q(s,z,y)\, dz\, ,
$$
proving \eqref{e:phi-y-partial-2}.

\noindent 
{\it Step 2:} Here we prove that for all $x\neq y$ and $t \in (0, T]$,
$T\ge 1$,
\begin{equation}\label{e:phi-y-partial-3}
\int_0^t \int_0^r \left|\partial_r \phi_y(r,x,s)\right|\, ds \, dr  \, \le \,  c_1(T) \,  t \frac{ \Phi(|x-y|^{-1})}{|x-y|^d} <+\infty\, .
\end{equation}
By \eqref{e:phi-y-partial-2} we have
\begin{eqnarray*}
\left|\partial_r \phi_y(r,x,s)\right|&\le &\int_{\R^d}\left| \partial_r p_z(r-s,x-z)\right|\, \left|q(s,z,y)-q(s,x,y)\right|\, dz\\
& &+\left|q(s,x,y)\right| \left|\int_{\R^d}\partial_r p_z(r-s,x-z)\, dz\right|
=:Q_y^{(1)}(r,x,s)+ Q_y^{(2)}(r,x,s)\, .
\end{eqnarray*}
For  $Q_y^{(1)}(r,x,s)$, by \eqref{e:difference-q-estimate}, \eqref{e:representaion-L},  \eqref{e:fract-der-est1} and Lemma \ref{l:convolution}(a) and (c), for $\beta_2 \in (0, \delta_1/2) \cap (0, \beta]$ and 
$\gamma\in ((2-\delta_1)\beta_2/2, \beta_2)$,
\begin{eqnarray*}
&&\int_0^t \int_0^r Q_y^{(1)}(r,x,s)\, ds\, dr \\
&&\le\ c_2 \int_0^t \int_0^r \int_{\R^d}\left|\LL^{\mathfrak{K}_z}p_z(r-s, x-z)\right| \left(|x-z|^{{\beta_2}-\gamma}\wedge 1\right)\\
& & \quad \times \left\{\left( \rho_{\gamma}^0 +\rho_{\gamma-{\beta_2}}^{{\beta_2}}\right)(s,x-y)+\left( \rho_{\gamma}^0 +\rho_{\gamma-{\beta_2}}^{{\beta_2}}\right)(s,z-y)\right\}dz\,  ds\,  dr\\
&&\le\ c_3 \int_0^t \int_0^r \left(\int_{\R^d} \rho_0^{{\beta_2}-\gamma}(r-s,x-z)dz \right) \Big(\rho_{\gamma}^0+\rho_{\gamma-{\beta_2}}^{{\beta_2}}\Big)(s,x-y)\, ds\, dr\\
& & \quad+c_3\int_0^t \int_0^r \int_{\R^d} \rho_0^{{\beta_2}-\gamma}(r-s,x-z)\Big(\rho_{\gamma}^0+\rho_{\gamma-{\beta_2}}^{{\beta_2}}\Big)(s,z-y)\, dz\, ds\, dr\\
& &\le\ c_4 \int_0^t \int_0^r (r-s)^{-1}\Phi^{-1}((r-s)^{-1})^{\gamma-{\beta_2}}\Big(\rho_{\gamma}^0+\rho_{\gamma-{\beta_2}}^{{\beta_2}}\Big)(s,x-y)\, ds\, dr\\
& &\quad +c_4 \int_0^t \Big(\rho_{{\beta_2}}^0 +\rho_0^{{\beta_2}}+\rho_{\gamma}^{{\beta_2}-\gamma}\Big)(r,x-y)\, dr\\
&&\le  c_4 \frac{\Phi(|x-y|^{-1})}{|x-y|^d} \int_0^t \int_0^r (r-s)^{-1}\Phi^{-1}((r-s)^{-1})^{\gamma-{\beta_2}}\Big(\Phi^{-1}(s^{-1})^{-\gamma}+\Phi^{-1}(s^{-1})^{{\beta_2}-\gamma}\Big)\, ds\, dr\\
& & \quad+c_4\frac{\Phi(|x-y|^{-1})}{|x-y|^d} \int_0^t  \big(\Phi^{-1}(r^{-1})^{-{\beta_2}}+1+\Phi^{-1}(r^{-1})^{-\gamma}\big)dr\\
& &\le c_5 \frac{\Phi(|x-y|^{-1})}{|x-y|^d} \int_0^t  \big(\Phi^{-1}(r^{-1})^{-{\beta_2}}+1+\Phi^{-1}(r^{-1})^{-\gamma}\Big)dr 
\le c_6  t \frac{\Phi(|x-y|^{-1})}{|x-y|^d} 
< +\infty\, .
\end{eqnarray*}
The second to
last inequality follows from Lemma \ref{l:convoluton-inequality}.

For $Q_y^{(2)}$, by \eqref{e:some-estimates-2b}, \eqref{e:q-estimate} and Lemma \ref{l:convoluton-inequality} we have
\begin{align*}
&\int_0^t \int_0^r Q_y^{(2)}(r,x,s)\, dr\, ds \le 
c_7 \int_0^t \int_0^r \Big(\rho_{{\beta_2}}^0+\rho_0^{{\beta_2}}\Big)(s,x-y) (r-s)^{-1}\Phi^{-1}((r-s)^{-1})^{-{\beta_2}}\, ds \,dr\\
&\le 2 c_7  \frac{\Phi(|x-y|^{-1})}{|x-y|^d}   \int_0^t \left(\int_0^r  (r-s)^{-1}\Phi^{-1}((r-s)^{-1})^{-{\beta_2}}\, ds\right) \,dr  \le c_8 t \frac{\Phi(|x-y|^{-1})}{|x-y|^d} 
 <+\infty\, .
\end{align*}

\noindent
{\it Step 3:} 
We claim that for fixed $s >0$ and $x, y \in \R^d$,
\begin{equation}\label{e:phi-y-partial-4}
\lim_{t\downarrow s}\phi_y(t,x,s)=q(s,x,y)\, .
\end{equation}
Assume $t \le T$, $T\ge 1$. 
For any $\delta >0$ we have
\begin{align*}
\lefteqn{\left|\int_{\R^d} p_z(t-s,x-z)\left(q(s,z,y)-q(s,x,y)\right)\, dz\right|}\\
&\le \int_{|x-z|\le \delta}  p_z(t-s,x-z)\left|q(s,z,y)-q(s,x,y)\right|\, dz \\
 & +\int_{|x-z|> \delta}  p_z(t-s,x-z)\left(|q(s,z,y)|+|q(s,x,y)|\right)\, dz=: J_1(\delta,t,s)+J_2(\delta,t,s)\, .
\end{align*}
By \eqref{e:difference-q-estimate}, for any $\eps >0$ there exists $\delta=\delta(s,x,y, T)>0$ such that if $|z-x|\le \delta$, then
$
|q(s,z,y)-q(s,x,y)|\le \eps\, .
$
Therefore, by Proposition \ref{p:upperestonp} and Lemma \ref{l:convolution}(a),
\begin{eqnarray*}
J_1(\delta,t,s)\le  \eps \int_{\R^d}p_z(t-s,x-z)\, dz\le  \eps (t-s) \int_{\R^d}\rho(t-s,z)\, dz  \le c_1 \eps \, .
\end{eqnarray*}
For $J_2(\delta,t,s )$, since $p_z(t-s,x-z)\le c_2 (t-s)\rho(t-s, x-z)\le c_2(t-s)\rho(0,x-z)$, by \eqref{e:q-estimate} we have
\begin{eqnarray*}
J_2(\delta,t,s)
\le  c_3(t-s)\left(\frac{\Phi(\delta^{-1})}{\delta^d}\int_{\R^d} 
\rho(s,z-y)\, dz +\rho(s,x-y)
\int_{|x-z|>\delta} \frac{\Phi(|x-z|^{-1})}{|x-z|^d}\, dz\right)
\end{eqnarray*}
where $c_3=c_3(T)>0$ is independent of $t$. 
By \eqref{e:convolution-integrability}, the term in parenthesis is finite. Hence, the last line converges to 0 as $t\downarrow s$. This and \eqref{e:some-estimates-2c} prove \eqref{e:phi-y-partial-4}.

\noindent
{\it Step 4:} By  \eqref{e:phi-y-partial-4}, we have that 
$$
\phi_y(t,x,s)-q(s,x,y)=\int_s^t \partial_r \phi_y(r,x,s)\, dr\, .
$$
Integrating both sides with respect to $s$ from $0$ to $t$, using first \eqref{e:phi-y-partial-3} and Fubini's theorem, and then \eqref{e:phi-y-partial-2} and \eqref{e:prop-p-L}, we get
\begin{eqnarray*}
\phi_y(t,x)-\int_0^t q(s,x,y)\, ds &=&\int_0^t \int_s^t \partial_r \phi_y(r,x,s)\, dr\, ds=\int_0^t \int_0^r  \partial_r \phi_y(r,x,s)\, ds\, dr\\
&=& \int_0^t \int_0^r \int_{\R^d} \LL^{\mathfrak{K}_z} p_z(r-s,\cdot)(x-z) q(s,z,y)\, dz\, ds\, dr\, .
\end{eqnarray*}
This proves that $t\mapsto \phi_y(t,x)$ is absolutely continuous and gives its Radon-Nykodim derivative 
\eqref{e:phi-y-partial}.  \qed

The following result is the counterpart of \cite[Lemma 3.6]{CZ}.

\begin{lemma}\label{l:L-on-phi-y}
For all  $t >0$, $x\neq y$ and $\eps \in [0,1]$, we have
\begin{equation}\label{e:L-on-phi-y}
\LL^{\mathfrak{K}_x, \eps} \phi_y(t,x)=\int_0^t \int_{\R^d} \LL^{\mathfrak{K}_x, \eps} p_z(t-s,\cdot) (x-z) q(s,z,y)\, dz\, ds\,
\end{equation}
and 
\begin{equation}\label{e:L-on-phi-y-cont}
t\mapsto \LL^{\mathfrak{K}_x}p_y(t,x-y) \text{ and } t\mapsto \LL^{\mathfrak{K}_x} \phi_y(t,x) \text{ are  continuous on }
(0,\infty)\, .
\end{equation}
Furthermore, if $\beta+\delta_1 >1$ and $\delta_1 \in (2/3,2)$ we also have
\begin{equation}\label{e:gradient-phi-y}
\nabla_x\phi_y(t,x)=\int_0^t \int_{\R^d} \nabla p_z(t-s,\cdot)(x-z) q(s,z,y)\, dz\, ds.
\end{equation}
\end{lemma}

\pf Fix $x\neq y$ and $T \ge 1$. 
In this proof we assume $0<t <T$  and all the constants will depend on $T$, but independent of $s$ and $t$.

\noindent
 (a) 
 By \eqref{e:psi1}, \eqref{e:intro-kappa},   \eqref{e:fract-der-est1}, \eqref{e:q-estimate} and  Lemma \ref{l:convolution}(b), for each 
 $s \in (0,t)$,
\begin{align}\label{e:Fubini1}
&\int_{\R^d}\int_{\R^d} |\delta_{p_z}(t-s, x-z;w)|\kappa(x,w)J(w)dw |q(s,z, y)|
dz\nn\\ 
& \le c_1\int_{\R^d}\rho(t-s, x-z)
\rho(s,z-y)
dz< \infty.
\end{align}
Thus we can use Fubini's theorem so that 
from \eqref{e:phi-y-def} we have  that for  $s \in (0,t)$,
\begin{equation}\label{e:L-on-phi-y2}
\LL^{\mathfrak{K}_x, \eps}\phi_y(t, \cdot ,s)(x)=\int_{\R^d}\LL^{\mathfrak{K}_x, \eps }p_z(t-s, \cdot)(x-z) q(s,z,y)\, dz, \quad  \eps \in [0,1]\, .
\end{equation}
Let $\beta_2 \in (0, \delta_1/2) \cap (0, \beta]$ and $\gamma \in (0, \beta_2)$.
By the definition of $\phi_y$, \eqref{e:phi-y-def}, and Fubini's theorem, 
using the notation \eqref{e:delta-f-def} 
we have for $\eps\in (0,1]$ and  $s \in (0,t)$,
\begin{eqnarray}
&&\left| \LL^{\mathfrak{K}_x, \eps}\phi_y(t, \cdot ,s)(x)\right| \nn\\
&=&\frac12\left|\int_{|w|>\eps}\left(\int_{\R^d}\delta_{p_z}(t-s,x-z;w)q(s,z,y)\, dz\right)\kappa(x,w)J(w)\, dw\right|\nn\\
&=&\frac12\left|\int_{\R^d}\left(\int_{|w|>\eps}\delta_{p_z}(t-s,x-z;w)\kappa(x,w)J(w)\, dw \right) q(s,z,y)\, dz\right|\nn\\
&\le &\frac12\int_{\R^d}\left(\int_{|w|>\eps}|\delta_{p_z}(t-s,x-z;w)|\kappa(x,w)J(w)\, dw \right) |q(s,z,y)-q(s,x,y)|\, dz\nn\\
& &+ \frac12\left|\int_{\R^d}\left(\int_{|w|>\eps}\delta_{p_z}(t-s,x-z;w)\kappa(x,w)J(w)\, dw \right) \, dz\right|\, |q(s,x,y)|\, .\nn
\end{eqnarray} 
By using \eqref{e:psi1},  \eqref{e:fract-der-est1}, \eqref{e:some-estimates-2a}, \eqref{e:q-estimate} and \eqref{e:difference-q-estimate} first and then using Lemma \ref{l:convolution}(a)--(b), we have that  for $\eps\in (0,1]$ and  $s \in (0,t)$,
\begin{eqnarray}
\lefteqn{\left| \LL^{\mathfrak{K}_x, \eps}\phi_y(t, \cdot ,s)(x)\right|} \nn\\
&\le& c_2 \int_{\R^d}\rho_0^{{\beta_2}-\gamma}(t-s,x-z)\left(\rho_{\gamma}^0+\rho_{\gamma-{\beta_2}}^{{\beta_2}}\right)(s,z-y)\, dz \nn\\
& & +c_2\left(\int_{\R^d}\rho_0^{{\beta_2}-\gamma}(t-s,x-z)\, dz\right)\left(\rho_{\gamma}^0+\rho_{\gamma-{\beta_2}}^{{\beta_2}}\right)(s,x-y) \nn\\
& & + c_2(t-s)^{-1}\Phi^{-1}((t-s)^{-1})^{-{\beta_2}}\left(\rho_0^{{\beta_2}}(s,x-y)+\rho_{{\beta_2}}^0(s,x-y)\right) \nn\\
&\le & c_2 \int_{\R^d}\rho_0^{{\beta_2}-\gamma}(t-s,x-z)\rho_{\gamma}^0(s,z-y)\, dz +c_2 \int_{\R^d}\rho_0^{{\beta_2}-\gamma}(t-s,x-z)\rho_{\gamma-{\beta_2}}^{{\beta_2}}(s,z-y)\, dz \nn\\
& & +c_3(t-s)^{-1}\Phi^{-1}((t-s)^{-1})^{\gamma-{\beta_2}}\left(\rho_{\gamma}^0+\rho_{\gamma-{\beta_2}}^{{\beta_2}}\right)(s,x-y) \nn\\
& & +c_3 (t-s)^{-1}\Phi^{-1}((t-s)^{-1})^{-{\beta_2}}\left(\rho_0^{{\beta_2}}(s,x-y)+\rho_{{\beta_2}}^0(s,x-y)\right) \nn\\
&\le& c_4\big((t-s)^{-1}\Phi^{-1}((t-s)^{-1})^{\gamma-2{\beta_2}}\Phi^{-1}(s^{-1})^{{\beta_2}-\gamma}\nn\\
&&+ (t-s)^{-1}\Phi^{-1}((t-s)^{-1})^{\gamma-{\beta_2}}\Phi^{-1}(s^{-1})^{{\beta_2}-\gamma} \nn\\
& &\quad+ (t-s)^{-1}\Phi^{-1}((t-s)^{-1})^{\gamma-{\beta_2}}\Phi^{-1}(s^{-1})^{-\gamma}+ 
(t-s)^{-1}\Phi^{-1}((t-s)^{-1})^{-{\beta_2}} \nn\\
& &\quad+s^{-1}\Phi^{-1}(s^{-1})^{-{\beta_2}}+ s^{-1}\Phi^{-1}(s^{-1})^{-\gamma}\big)\rho(0, x-y) \nn\\
& \le& c_5 (t-s)^{-1}\Phi^{-1}((t-s)^{-1})^{\gamma-{\beta_2}}s^{-1}\Phi^{-1}(s^{-1})^{-\gamma}\rho(0, x-y). \label{e:Iepsilon}
\end{eqnarray}
In the last inequality above we have used the inequality 
$$\Phi^{-1}(s^{-1})^{{\beta_2}} \le a_1^{-{\beta_2}/\delta_1} \Phi^{-1}(T^{-1})^{-{\beta_2}}  s^{-{\beta_2}/\delta_1} \le  a_1^{-{\beta_2}/\delta_1} \Phi^{-1}(T^{-1})^{-{\beta_2}}  T^{1-{\beta_2}/\delta_1}  s^{-1}.$$
Using the fact that $x\neq y$ and Lemma \ref{l:convoluton-inequality} we see that the 
term on the right hand side of \eqref{e:Iepsilon} is  integrable in $s\in (0, t)$.
Moreover, by \eqref{e:intro-kappa}, \eqref{e:psi1}, \eqref{e:q-estimate}   and Proposition \ref{p:upperestonp},
\begin{align}
&\int_{|w|>\eps} \int_0^t |\delta_{\phi_y}(t,x,s;w)|\kappa(x,w)J(w)\, ds\, dw\nn\\ 
 \le &2\kappa_1\gamma_0 C_2
\int_{|w|>\eps} \int_0^t\int_{\R^d}p_z(t-s,x-z) (\rho_0^{\beta_2}(s,z-y)+\rho^0_{\beta_2}(s,z-y))dz
j(|w|)\, ds\, dw \nn\\
&+ \kappa_1\gamma_0 C_2
\int_{|w|>\eps} \int_0^t\int_{\R^d}p_z(t-s,x\pm w -z) (\rho_0^{\beta_2}(s,z-y)+\rho^0_{\beta_2}(s,z-y))dz
j(|w|)\, ds\, dw\nn\\
 \le &c_6\int_{|w|>\eps} j(|w|)dw
 \int_0^t(t-s)  \left(\int_{\R^d}\rho(t-s,x-z) (\rho_0^{\beta_2}(s,z-y)+\rho^0_{\beta_2}(s,z-y))dz \right)
\, ds\, \nn\\
&+ c_6 j(\eps)
 \int_0^t\int_{\R^d}
(t-s) \left(\int_{\R^d} \rho(t-s,x\pm w -z)dw \right)
 (\rho_0^{\beta_2}(s,z-y)+\rho^0_{\beta_2}(s,z-y))dz
\, ds\,, \label{e:FubiniJ1}
\end{align}
which is, by Lemma \ref{l:convolution}(a)--(b), less than or equal to 
\begin{align}
&c_7(\eps)  \left(
 \int_0^t  s^{-1}   \Phi^{-1}(s^{-1})^{-{\beta_2}} \rho(t, x-y) \, ds+
 \int_0^t\int_{\R^d}
 (\rho_0^{\beta_2}(s,z-y)+\rho^0_{\beta_2}(s,z-y))dz
\, ds\, \right) \nn\\
& \le 
c_8(\eps) \left(
 \int_0^t  s^{-1}   \Phi^{-1}(s^{-1})^{-{\beta_2}} ds \rho(t, x-y)+
 \int_0^t     s^{-1}   \Phi^{-1}(s^{-1})^{-{\beta_2}} ds\, \right)  < \infty. \label{e:FubiniJ2}
\end{align}
Thus we can apply Fubini's theorem to see that, by \eqref{e:L-on-phi-y2},
\eqref{e:L-on-phi-y} holds for $ \eps \in (0, 1]$. Moreover, by Fubini's theorem 
 and the dominated convergence theorem in the first equality and the second equality below respectively:
$$
{\LL}^{{\mathfrak K}_x}\phi_y(t,x)=\lim_{\eps \downarrow 0}
\int_0^t \LL^{\mathfrak{K}_x, \eps}\phi_y(t, \cdot ,s)(x)\, ds
=\int_0^t \lim_{\eps\downarrow 0}  \LL^{\mathfrak{K}_x, \eps}\phi_y(t, \cdot ,s)(x)\, ds\, ,
$$
which together with \eqref{e:L-on-phi-y2} yields \eqref{e:L-on-phi-y} for $\eps=0$.

\noindent
(b) Now we prove \eqref{e:L-on-phi-y-cont}.  Note that, by Lemma \ref{l:partial-time}(b),  $t\mapsto \delta_{p_y}(t,x-y;z)=p_y(t,x-y+z)+p_y(t,x-y-z)-2p_y(t,x-y)$ is continuous. Let $\eps\in (0,t)$. By \eqref{e:delta-p-kappa}, 
\begin{eqnarray*}
|\delta_{p_y}(t,x-y;z)| &\le & c_{11} \left(\Phi^{-1}(t^{-1})|z|^2\wedge 1\right)t \left(\rho(t,x-y\pm z)+\rho(t,x-y)\right)\\
&\le &c_{12}\frac{t}{\eps}\left(\Phi^{-1}(\eps^{-1})|z|^2\wedge 1\right)\eps \left(\rho(\eps,x-y\pm z)+\rho(\eps,x-y)\right)\, .
\end{eqnarray*}
By \eqref{e:psi1} and  the proof of \eqref{e:fract-der-est1} we see that the right-hand side multiplied by $\kappa(x,z)J(z)$ is integrable with respect to $dz$. This shows that the family $\{\delta_{p_y}(t,x-y;z)\kappa(x,z)J(z):\, t\in (\eps,T)\}$ is dominated by an integrable function. Now 
by the dominated convergence theorem we see that 
 $t\mapsto \LL^{\mathfrak{K}_x}p_y(t,x-y)$ is continuous on $(0,T]$.

Let $\beta_2 \in (0,\delta_1/2)\cap (0,\beta]$ and $\gamma \in (0,\beta_2)$.
By \eqref{e:Iepsilon},
\begin{align}
\left|\LL^{\mathfrak{K}_x}\phi_y(t,x,s)\right|
\le c_5 (t-s)^{-1}\Phi^{-1}((t-s)^{-1})^{\gamma-\beta_2}s^{-1}
\Phi^{-1}(s^{-1})^{-\gamma}\rho(0, x-y)
\, .\label{e:ub4I-2}
\end{align}
Note that for $0<t\le t+h\le T$,
\begin{eqnarray}
&&{\LL}^{{\mathfrak K}_x}\phi_y(t+h,x)-{\LL}^{{\mathfrak K}_x}\phi_y(t,x)\nn\\
&&=\int^{t+h}_t
\LL^{\mathfrak{K}_x}\phi_y(t+h,x,s)ds +\int^t_0\left( \LL^{\mathfrak{K}_x}\phi_y(t+h,x,s)-\LL^{\mathfrak{K}_x}\phi_y(t,x,s)\right)ds.
\end{eqnarray}
When $h\le t/2$, by \eqref{e:psi-inverse-sbm} and \eqref{e:lsc-inverse}, we have
\begin{eqnarray*}
&&\int^{t+h}_t(t+h-s)^{-1}
\Phi^{-1}((t+h-s)^{-1})^{\gamma-\beta_2}s^{-1}
\Phi^{-1}(s^{-1})^{-\gamma}ds\nn\\
&&=\int^h_0r^{-1}\Phi^{-1}(r^{-1})^{\gamma-\beta_2}(t+h-r)^{-1}
\Phi^{-1}((t+h-r)^{-1})^{-\gamma}dr\nn\\
&&\le c_{13}\int^h_0r^{-1}\Phi^{-1}(r^{-1})^{\gamma-\beta_2}(t-r)^{-1}
\Phi^{-1}((t-r)^{-1})^{-\gamma}dr\,,
\end{eqnarray*}
and so by Lemma \ref{l:rsconvolution-inequality} and \eqref{e:ub4I-2} we get
\begin{equation}\label{e:contoflphi1}
\lim_{h\to0}\int^{t+h}_t\LL^{\mathfrak{K}_x}\phi_y(t+h,x,s)ds=0.
\end{equation}
Note that,  by \eqref{e:Fubini1} we can apply  
the
dominated convergence theorem  and use the continuity of 
$t\mapsto \LL^{\mathfrak{K}_x}p_y(t,x-y)$ so that  for each $s \in (0, t)$,
\begin{align}\label{e:contoflphi22}
&\lim_{h\to0}(\LL^{\mathfrak{K}_x}\phi_y(t+h,x,s)-\LL^{\mathfrak{K}_x}\phi_y(t,x,s))\nn\\
=& 
\int_{\R^d}\lim_{h\to0} (\LL^{\mathfrak{K}_x}p_z(t+h-s, \cdot)(x-z)-\LL^{\mathfrak{K}_x}p_z(t-s, \cdot)(x-z) )q(s,z,y)\, dz
=0.
\end{align}
By Lemma \ref{l:convoluton-inequality}, 
$s\mapsto  (t-s)^{-1}\Phi^{-1}((t-s)^{-1})^{\gamma-\beta_2}s^{-1}
\Phi^{-1}(s^{-1})^{-\gamma}$
is integrable in $(0, t)$, so using \eqref{e:ub4I-2}, we can apply  
the
dominated convergence theorem and use \eqref{e:contoflphi22} to get that
\begin{align}\label{e:contoflphi2}
&\lim_{h\to0}\int^{t}_0(\LL^{\mathfrak{K}_x}\phi_y(t+h,x,s)-\LL^{\mathfrak{K}_x}\phi_y(t,x,s))ds
=0.
\end{align}
Combining \eqref{e:contoflphi1}--\eqref{e:contoflphi2} we get the desired continuity.

\noindent
(c)
Finally we show \eqref{e:gradient-phi-y}. 
Since $\beta+\delta_1 >1$ and $\delta_1 \in (2/3,2)$, we can and will choose 
$\beta_2 \in ( 0 \vee(1-\delta_1)  ,  \delta_1/2) \cap (0, \beta]$ and 
$\gamma\in (0,\beta_2 \wedge (\beta_2+\delta_1-1) \wedge (\delta_1-2\beta_2) )$.
For example, one can take $\beta_2=\beta \wedge (1/3)$.

For each fixed $0<s<t$ and 
$he_i=(0, \dots, h, \dots, 0) \in \R^d$ with $|h| \le 1/(2\Phi^{-1}((t-s)^{-1}))$, by \eqref{e:difference-p-kappa},  \eqref{e:2.9}, \eqref{e:Berall}  and \eqref{e:q-estimate} we have
\begin{align}
\label{e:works1}
&\frac{1}{h}\big|p_z(t-s,x-z+he_i)-p_z(t-s,x-z)\big||q(s,z,y)| \nn \\
&\le c\frac{1}{h}\left((\Phi^{-1}((t-s)^{-1})|h|)\wedge 1\right)(t-s)\left(\rho(t-s,x-z+he_i)+\rho(t-s,x-z)\right)|q(s,z,y)|\nn \\
&\le 2^{d+2}c (t-s) \Phi^{-1}((t-s)^{-1}) \rho(t-s,x-z)(\rho_0^{\beta_2} +\rho^0_{\beta_2})(s,z-y)
\end{align}
which is integrable in $z \in \R^d$ by Lemma \ref{l:convolution}(b).
Thus we can use  the dominated convergence theorem 
and  \eqref{e:phi-y-def} to get
that for  $s \in (0,t)$,
\begin{equation}\label{e:L-on-phi-y21}
\partial_i \phi_y(t, \cdot ,s)(x)=\int_{\R^d}\partial_i p_z(t-s, \cdot)(x-z) q(s,z,y)\, dz\, .
\end{equation}

Let
\begin{eqnarray}
\partial_i \phi_y(t, \cdot ,s)(w) 
&=&\int_{\R^d}\partial_i p_z(t-s, \cdot) ( w-z) q(s, z, y)\, dz\nn\\
&=&\ind_{[t/2,t)}(s)\int_{\R^d} \partial_i p_z(t-s, \cdot) ( w-z)\, (q(s,z,y)-q(s,w,y)) dz\nn\\
& &+\ \ind_{[t/2,t)}(s)\int_{\R^d}\partial_i p_z(t-s, \cdot) ( w-z)\, q(s,w,y) dz\nn\\
& &+\ \ind_{(0,t/2)}(s)\int_{\R^d} \partial_i p_z(t-s, \cdot) ( w-z)\, q(s,z,y) dz\nn\\
&=:&\ind_{[t/2,t)}(s)R_1(t, s, w, y)+\ind_{[t/2,t)}(s)R_2(t, s, w, y)\nn\\
& & +\ \ind_{(0,t/2)}(s)R_3(t, s, w, y)\, . \label{e:comp1}
\end{eqnarray}

Let $x' \in B( x, |x-y|/4)$. 
Then 
it follows from Proposition \ref{p:upperestonp} and \eqref{e:difference-q-estimate} that for 
$s\in [t/2,t)$, 
\begin{align}
&\big |R_1(t, s, x', y)\big |\nn\\
& \le
  \int_{\R^d}| \partial_i p_z(t-s, \cdot) ( x'-z)||q(s,z,y)-q(s,x',y)| dz  \nn \\
& \le \int_{\R^d} \left( (t-s)\Phi^{-1}((t-s)^{-1})\rho(t-s,x'-z)\left(|x'-z|^{{\beta_2}-\gamma}\wedge 1\right) \left(\rho_{\gamma}^0+\rho_{\gamma-{\beta_2}}^{{\beta_2}}\right)(s,x'-y)\right.\nn\\
&   \left.  \ \ +(t-s)\Phi^{-1}((t-s)^{-1})\rho(t-s,x'-z)\left(|x'-z|^{{\beta_2}-\gamma}\wedge 1\right)\left(\rho_{\gamma}^0+\rho_{\gamma-{\beta_2}}^{{\beta_2}}\right)(s,z-y)\right)\, \, dz\nn \\
& =(t-s) \left( \int_{\R^d} \rho_{-1}^{{\beta_2}-\gamma}(t-s,x'-z) dz  \right)\left(\rho_{\gamma}^0+\rho_{\gamma-{\beta_2}}^{{\beta_2}}\right)(s,x'-y). \nn\\
&  \ \   +(t-s) \int_{\R^d} \rho_{-1}^{{\beta_2}-\gamma}(t-s,x'-z)\rho_{\gamma}^0(s,z-y) dz\nn\\&  \ \   +(t-s)\int_{\R^d} \rho_{-1}^{{\beta_2}-\gamma}(t-s,x'-z)\rho_{\gamma-{\beta_2}}^{{\beta_2}}(s,z-y) dz\nn\\
& \le c_9  \left(\Phi^{-1}((t-s)^{-1})^{1-{\beta_2}+\gamma}\left(\rho_{\gamma}^0+\rho_{\gamma-{\beta_2}}^{{\beta_2}}\right)(s,x'-y) \right.\nn \\
& \ \  + \left(\Phi^{-1}((t-s)^{-1})^{1-2{\beta_2}+\gamma}\Phi^{-1}(s^{-1})^{-\gamma+{\beta_2}}
+\Phi^{-1}((t-s)^{-1})^{1-{\beta_2}+\gamma}\Phi^{-1}(s^{-1})^{-\gamma}
\right.\nn\\
&\ \ \left. \left. \left.+(t-s)s^{-1}\Phi^{-1}(s^{-1}) (
\Phi^{-1}(s^{-1})^{-\gamma} +\Phi^{-1}(s^{-1})^{-{\beta_2}})\right)  \rho(t, x'-y)\right) \right)\nn \\
&\le c_{10}\left(\Phi^{-1}((t-s)^{-1})^{1-{\beta_2}+\gamma} 
\Phi^{-1}(s^{-1})^{-\gamma+{\beta_2}}
  \right. \nn\\
 &+ \Phi^{-1}((t-s)^{-1})^{1-2{\beta_2}+\gamma}\Phi^{-1}(s^{-1})^{-\gamma+{\beta_2}}
+\Phi^{-1}((t-s)^{-1})^{1-{\beta_2}+\gamma}\Phi^{-1}(s^{-1})^{-\gamma}
\nn\\
&\ \ \left. +(t-s)s^{-1}\Phi^{-1}((t-s)^{-1}) 
\Phi^{-1}(s^{-1})^{-\gamma} \right)  \rho(t, (x-y)/2). \label{e:comp2}
\end{align}
Here the third inequality follows from Lemma \ref{l:convolution}(a)--(b). 
Since $\delta_1>2/3>1/2$ and $\gamma< \delta_1+\beta_2-1$, 
using Lemma \ref{l:convoluton-inequality} (so that 
$ \int_{t/2}^t \Phi^{-1}((t-s)^{-1})^{1-\beta_2+\gamma}ds$ 
and $ \int_{t/2}^t (t-s) \Phi^{-1}((t-s)^{-1})ds$ are finite)
 it is straightforward
 to see that the function on the right-hand side above is integrable in $s$ over $[t/2,t)$.

Next, for $s\in [t/2,t)$,  using \eqref{e:q-estimate} in the second and 
\eqref{e:reinstated-estimate} in the third line below,
\begin{align}
&\big |R_2(t, s, x', y)\big | = \left|\int_{\R^d} \partial_i p_z(t-s, \cdot) ( x'-z)\, dz\right| q(s,x',y)\nn\\
&\le  \left|\int_{\R^d} \partial_i p_z(t-s, \cdot) ( x'-z)\, dz\right|  
\left(\rho_0^{\beta_2}+\rho_{\beta_2}^0\right)(s,x',y)\nn\\
& \le c \Phi^{-1}((t-s)^{-1})^{1-\beta_2} \rho(t,x'-y) \nn\\
& \le c  \Phi^{-1}((t-s)^{-1})^{1-\beta_2} \rho(t,(x-y)/2).\label{e:comp3}
\end{align}
Since $\int_{t/2}^t \Phi^{-1}((t-s)^{-1})^{1-\beta_2}\, ds  <\infty$ because $\beta_2+\delta_1>1$, 
the right-hand side above is integrable in $s$ over $[t/2,t)$.

Finally for $s\in (0,t/2]$, since $\beta_2 < \delta_1/2$, 
\begin{align}
&\big |R_3(t, s, x', y)\big |  \le
\int_{\R^d}| \partial_i p_z(t-s, \cdot) ( x'-z)|q(s,z,y)| dz \nn \\
& \le c \int_{\R^d}  (t-s)\Phi^{-1}((t-s)^{-1})\rho(t-s,x'-z)
\left(
\rho_{\beta_2}^0+\rho^{\beta_2}_0   \right)(s,z-y) dz\nn\\
& = c  (t-s)\int_{\R^d} \rho_{-1} (t-s,x'-z)\left(
\rho_{\beta_2}^0+\rho^{\beta_2}_0   \right)(s,z-y)dz\nn\\
&\le c  (t-s)\Big( (t-s)^{-1}\Phi^{-1}((t-s)^{-1})^{1-\beta_2} + (t-s)^{-1}\Phi^{-1}((t-s)^{-1})
\nn\\
&\quad
+ (t-s)^{-1}
 \Phi^{-1}((t-s)^{-1})
\Phi^{-1}(s^{-1})^{-\beta_2}  +  \Phi^{-1}((t-s)^{-1}) s^{-1}
\Phi^{-1}(s^{-1})^{-\beta_2}
   \Big)\rho(t, x'-y)\nn\\
   &\le c  \Big( \Phi^{-1}((t-s)^{-1})+ 
 \Phi^{-1}((t-s)^{-1})
\Phi^{-1}(s^{-1})^{-\beta_2}\nn\\
&\quad +(t-s) \Phi^{-1}((t-s)^{-1}) s^{-1}
\Phi^{-1}(s^{-1})^{-\beta_2}
   \Big)\rho(t, x'-y),   \label{e:comp4}
\end{align}
which is integrable using Lemma \ref{l:convoluton-inequality}.

 Hence we can use the dominated convergence theorem and \eqref{e:L-on-phi-y21} to conclude that
\begin{align*}
&\lim_{h\to 0}\frac{1}{h}\big(\phi_y(t,x+w)-\phi_y(t,x) \big)
=\lim_{h\to 0}\int_0^t \int_0^1 \partial_i \phi_y(t, \cdot ,s) ( x+ \theta w)\, d\theta ds ds\nn\\
&=\int_0^t  
\partial_i \phi_y(t, \cdot ,s) ( x) ds=\int_0^t 
\int_{\R^d}\partial_i p_z(t-s, \cdot)(x-z) q(s,z,y)\, dz ds,
\end{align*}
which gives \eqref{e:gradient-phi-y}.
 \qed

\subsection{Estimates and Smoothness of $p^\kappa(t, x, y)$}

Now we define and study the function 
\begin{equation}\label{e:p-kappa}
p^{\kappa}(t,x,y):=p_y(t,x-y)+\phi_y(t,x)=p_y(t,x-y)+\int_0^t \int_{\R^d}p_z(t-s,x-z)q(s,z,y)\, dz\, ds\, .
\end{equation}

\begin{lemma}\label{l:p-kappa-difference}
(1) For every $T\ge 1$ and $\beta_2\in (0, \beta]\cap(0, \delta_1/2)$, 
there is a constant $c_1=
c_1(T, d,\delta_1,\beta_2,\gamma,\kappa_0,\kappa_1,\kappa_2)>0$ 
so that for all $t\in (0,T]$ and $x,y\in \R^d$,
$
p^{\kappa}(t,x,y)\le c_1 t \rho(t,x-y).
$
(2) For any 
$\gamma\in (0,\delta_1) \cap (0,1]$
 and $T \ge 1$ there exists $c_2=
c_2(T, d,\delta_1,\beta_2,\gamma,\kappa_0,\kappa_1,\kappa_2)>0$
 such that for all $x, x', y\in \R^d$ and $t\in (0, T]$,
$$
\left|p^{\kappa}(t,x,y)-p^{\kappa}(t,x',y)\right| \le c_2 |x-x'|^{\gamma}\ t\Big(\rho_{-\gamma}^0(t,x-y)+\rho_{-\gamma}^0(t,x'-y)\Big).
$$
\end{lemma}
\pf 
Throughout this proof we assume that $x, x', y\in \R^d$ and $t\in (0, T]$.

\noindent
(1) By the estimate of $p_z$ (Proposition \ref{p:upperestonp}), \eqref{e:q-estimate}, Lemma \ref{l:convolution}(c), \eqref{e:nonincrease-gamma} and \eqref{e:nonincrease-beta}, we have
\begin{eqnarray}
\lefteqn{\int_0^t \int_{\R^d}p_z(t-s,x-z) |q(s,z,y)|\, dz\, ds} \nonumber \\
&\le & c_1\int_0^t \int_{\R^d} (t-s)\rho(t-s,x-z)
\left(\rho_{\beta_2}^0+\rho_0^{\beta_2}\right)(s,z-y)\, dz\, ds 
\nonumber \\
&\le & c_2 t\left(\rho_{\beta_2}^0+\rho_0^{\beta_2}\right)(t,x-y)
\label{e:main-proof-i} \\ 
&\le & 2 \Phi^{-1} (T^{-1})^{-\beta_2} 
c_2t\rho(t,x-y), \quad \text{ for all } t\in (0,T]\, .\nonumber
\end{eqnarray}
Therefore, $p^{\kappa}(t,x,y)\le p_y(t,x-y)+|\phi_y(t,x)|\le c_4 t \rho(t,x-y)$.

\noindent
(2)
We have by \eqref{e:difference-p-kappa} and the fact that $\gamma\le 1$,
\begin{eqnarray*}
|p_z(t,x-z)-p_z(t,x'-z)|
&\le & c_1 |x-x'|^{\gamma} t \Phi^{-1}(t^{-1})^{\gamma}\big(\rho(t,x-z)+\rho(t,x'-z)\big)\\
&=&  c_1 |x-x'|^{\gamma} t \big(\rho_{-\gamma}^0(t,x-z)+\rho_{-\gamma}^0(t,x'-z)\big)\, .
\end{eqnarray*}
Thus, by \eqref{e:q-estimate} and a change of the variables, we further have
\begin{align*}
&|\phi_y(t,x)-\phi_y(t,x')|
\le  \int_0^t \int_{\R^d}|p_z(t-s,x-z)-p_z(t-s,x'-z)| \, |q(s,z,y)|\, dz \, ds\\
\le & c_2 |x-x'|^{\gamma}\int_0^t \int_{\R^d}(t-s) 
\Big(\rho_{-\gamma}^0(t-s,x-z)+\rho_{-\gamma}^0(t-s,x'-z)\Big)
\big(\rho_0^{\beta_2}+\rho_{\beta_2}^0\big)(s,z-y)\, dz\, ds\\
\le & c_3 |x-x'|^{\gamma} t
\Big(\rho_{-\gamma+\beta_2}^0(t,x-y)+\rho_{-\gamma}^{\beta_2}(t,x-y)+\rho_{-\gamma+\beta_2}^0(t,x'-y)+\rho_{-\gamma}^{\beta_2}(t,x'-y)\Big)\\
\le & 2c_3 \Phi^{-1}(T^{-1})^{-\beta_2} 
|x-x'|^{\gamma} t \big(\rho_{-\gamma}^0(t,x-y)+\rho_{-\gamma}^0(t,x'-y)\big), \quad \text{ for all } t\in (0,T] \, .
\end{align*}
Since $\gamma\in (0,\delta_1)$, the penultimate line follows from \eqref{e:convolution-3} (with $\theta=0$), 
and the last line by \eqref{e:nonincrease-gamma} and \eqref{e:nonincrease-beta}.
The claim of the lemma follows by combining the two estimates. \qed

The following result is the counterpart  of \cite[Lemma 3.7]{CZ}.

\begin{lemma}\label{l:continuity-p-kappa}
The function $p^{\kappa}(t,x,y)$ defined in $\eqref{e:p-kappa}$ is jointly continuous on $(0, \infty)\times \R^d \times \R^d$.
\end{lemma}
\pf The joint continuity of $p_y(t,x-y)$ was shown in Lemma \ref{l:jcontoffzkernel}.
For $\phi_y(t,x)$ we use \eqref{e:def-phi-y-2} and the joint continuity of $q(s,z,y)$ on $(0, \infty)\times \R^d\times \R^d$ together with the  dominated convergence theorem. This is justified by the estimates $p_z(t-s,x-z)\le c_1 (t-s)\rho(t-s,x-z)$ and \eqref{e:q-estimate} which yield that $|p_z(t-s,x-z)q(s,z,y)|\le 
c_2 (t-s)\rho(t-s)\left(\rho_0^{\beta_2}+\rho_{\beta_2}^0\right)(s,z-y)$ for $\beta_2\in (0, \beta]\cap (0, \delta_1/2)$. 
The latter function is integrable over $(0,t]\times \R^d$ with respect to $ds \, dz$ by Lemma \ref{l:convolution}. \qed

Now we define the operator $\LL^{\kappa}$ 
as in \eqref{e:intro-operator} 
which can be rewritten as 
\begin{align}
\LL^{\kappa}f(x)=
\LL^{\kappa, 0}f(x)=\lim_{\eps \downarrow 0} \LL^{\kappa, \eps}f(x), \quad 
\text{where }\LL^{\kappa, \eps}f(x)=\frac{1}{2}\int_{|z| > \eps} \delta_f(x;z)  \kappa(x,z)J(z)\, dz.\label{e:operator-L-kappa}
\end{align}
Note that for a fixed $x\in \R^d$, it holds that $\LL^{\kappa}f(x)=\LL^{{\mathfrak K}_x}f(x)$. This will be used later on.

The following result is the counterpart  of \cite[Lemma 4.2]{CZ}.

\begin{lemma}\label{l:fract-der-p-kappa}
For every $T\ge 1$, there is a constant 
$c_1=c_1 (T, d,\delta_1, a_1, \beta, C_*,  \gamma_0,  \kappa_0, \kappa_1, \kappa_2) >0$
such that for all $\eps \in [0,1]$, 
\begin{equation}\label{e:fract-der-p-kappa-1}
|\LL^{\kappa, \eps} p^{\kappa}(t, \cdot, y)(x)|\le c_1 \rho(t,x-y), \quad \text{ for all } t\in(0,T] \text{ and } x,y\in \R^d, x\neq y
\end{equation}
and if $\beta+\delta_1 >1$ and $\delta_1 \in (2/3,2)$ we also have
\begin{equation}\label{e:fract-der-p-kappa-2}
\left|\nabla_x p^{\kappa}(t,x,y)\right|\le c_1 t \Phi^{-1}(t^{-1})\rho(t,x-y)\quad \text{ for all } t\in(0,T] \text{ and } x,y\in \R^d, x\neq y \, .
\end{equation}
\end{lemma}
\pf 
By \eqref{e:fract-der-est1} and the fact that for fixed $x$, $\LL^{\kappa, \eps}f(x)=\LL^{\mathfrak{K}_x, \eps}f(x)$ for $\eps \in [0,1]$, we see that
$$
|\LL^{\kappa} p_y(t, \cdot)(x-y)|\le c_1 \rho(t,x-y), \qquad \text{for all } t\in (0,T] \text{ and }\eps \in [0,1].
$$

Let $\eps \in [0,1]$. By recalling the definition \eqref{e:def-phi-y-2} of $\phi_y$ and using \eqref{e:L-on-phi-y}, we have
\begin{eqnarray*}
\LL^{\kappa, \eps}\phi_y(t,x) &=& \int_{t/2}^t \int_{\R^d}\LL^{\mathfrak{K}_x, \eps}p_z(t-s,\cdot)(x-z)\left(q(s,z,y)-q(s,x,y)\right)\, dz\, ds\\
 & & +\int_{t/2}^t \left( \int_{\R^d}\LL^{\mathfrak{K}_x, \eps}p_z(t-s,\cdot)(x-z)\, dz\right) q(s,x,y)\, ds \\
 & & +\int_0^{t/2} \int_{\R^d}\LL^{\mathfrak{K}_x, \eps}p_z(t-s,\cdot)(x-z) q(s,z,y)\, dz\, ds\\
 &=:& Q_1(t,x,y)+Q_2(t,x,y)+Q_3(t,x,y)\, .
\end{eqnarray*}
Let $\beta_2 \in (0, \delta_1/2) \cap (0, \beta]$.
For $Q_1(t,x,y)$ we use \eqref{e:fract-der-est1}, Lemmas  \ref{l:psi-and-f}(b),   \ref{l:convoluton-inequality} and  \ref{l:convolution}(a) and (c) to  get that for any 
$\gamma \in ((2-\delta_1){\beta_2}/2, {\beta_2})$,
\begin{eqnarray*}
|Q_1(t,x,y)|&\le & c_1 \int_{t/2}^t \left( \int_{\R^d} \rho_0^{{\beta_2}-\gamma}(t-s,x-z)\, dz \right)\left(\rho_{\gamma}^0+\rho_{\gamma-{\beta_2}}^{{\beta_2}}\right)(s,x-y)\, ds\\
& & +c_1\int_{t/2}^t \int_{\R^d} \rho_0^{{\beta_2}-\gamma}(t-s,x-z)\left(\rho_{\gamma}^0+\rho_{\gamma-{\beta_2}}^{{\beta_2}}\right)(s,z-y)\, dz\, ds\\
&\le & c_2 \left(\rho_{\gamma}^0+\rho_{\gamma-{\beta_2}}^{{\beta_2}}\right)(t,x-y)    \int_{0}^t \int_{\R^d} \rho_0^{{\beta_2}-\gamma}(t-s,x-z)\, dz  ds\\
& & +c_1\int_{0}^t \int_{\R^d} \rho_0^{{\beta_2}-\gamma}(t-s,x-z)\left(\rho_{\gamma}^0+\rho_{\gamma-{\beta_2}}^{{\beta_2}}\right)(s,z-y)\, dz\, ds\\
&\le & c_3 \rho_{\gamma-{\beta_2}}^{0}(t,x-y) \Phi^{-1}(t^{-1})^{-{\beta_2}-\gamma} +c_3\left(\rho_{{\beta_2}}^0+\rho_{\gamma}^{{\beta_2}-\gamma}+\rho_0^{{\beta_2}}\right)(t,x-y)\\
&\le & c_4 \rho(t,x-y), \quad \text{ for all } t\in (0,T]\, ,
\end{eqnarray*}
where the last two lines follow from \eqref{e:nonincrease-gamma} and \eqref{e:nonincrease-beta}.

For $Q_2(t,x,y)$, by \eqref{e:some-estimates-2a}, \eqref{e:q-estimate}, Lemmas  \ref{l:psi-and-f}(b), \ref{l:convoluton-inequality}, \eqref{e:nonincrease-gamma} and \eqref{e:nonincrease-beta},
\begin{eqnarray*}
&&|Q_2(t,x,y)|\le c_5 \int_{t/2}^t (t-s)^{-1}\Phi^{-1}((t-s)^{-1})^{-{\beta_2}}\left(\rho_{{\beta_2}}^0+\rho_0^{{\beta_2}}\right)(s,x-y)\, ds\\
&&\le c_6 \left(\rho_{{\beta_2}}^0+\rho_0^{{\beta_2}}\right)(t,x-y)   \int_{0}^t (t-s)^{-1}\Phi^{-1}((t-s)^{-1})^{-{\beta_2}}\, ds\\ 
&& \le c_7 \rho(t,x-y) \Phi^{-1}(t^{-1})^{-{\beta_2}}  \le  c_7 \Phi^{-1}(T^{-1})^{-{\beta_2}} \rho(t,x-y), \quad \text{ for all } t\in (0,T]\, .
\end{eqnarray*}

For $Q_3(t,x,y)$, by \eqref{e:fract-der-est1}, \eqref{e:q-estimate}, Lemma \ref{l:convolution}(c), \eqref{e:nonincrease-gamma} and \eqref{e:nonincrease-beta},
\begin{eqnarray*}
|Q_3(t,x,y)|&\le & c_7\int_0^{t/2}\int_{\R^d} \rho(t-s,x-z)\left(\rho_{{\beta_2}}^0+\rho_0^{{\beta_2}}\right)(s,z-y)\, dz\, ds \\
&\le & 2 \frac{c_7}{t}\int_0^{t}\int_{\R^d} (t-s)\rho(t-s,x-z)\left(\rho_{{\beta_2}}^0+\rho_0^{{\beta_2}}\right)(s,z-y)\, dz\, ds\\
&\le & c_8 \left(\rho_{{\beta_2}}^0+\rho_0^{{\beta_2}}\right)(t,x-y)\le 2c_8  \Phi^{-1}(T^{-1})^{-\beta} \rho(t,x-y)\, .
\end{eqnarray*}
Combining the above calculations and \eqref{e:p-kappa} we obtain \eqref{e:fract-der-p-kappa-1}.

\noindent
(ii) 
Since $\beta+\delta_1 >1$ and $\delta_1 \in (2/3,2)$, we can and will choose 
$\beta_2 \in ( 0 \vee(1-\delta_1)  ,  \delta_1/2) \cap (0, \beta]$ and 
$\gamma\in (0,\beta_2 \wedge (\beta_2+\delta_1-1) \wedge (\delta_1-2\beta_2) )$.
By \eqref{e:gradient-phi-y} and \eqref{e:comp1}--\eqref{e:comp4} we have
\begin{align}
&|\nabla_x\phi_y(t,x)| \le c_1 \rho(t, x-y)  \Bigg(\int_0^{t/2}   
\Phi^{-1}((t-s)^{-1})+
 \Phi^{-1}((t-s)^{-1})
\Phi^{-1}(s^{-1})^{-\beta_2}
\nn \\
&\qquad 
+(t-s) \Phi^{-1}((t-s)^{-1}) s^{-1}
\Phi^{-1}(s^{-1})^{-\beta_2}
  ds  \nn\\
 & +
\int_{t/2}^t \Phi^{-1}((t-s)^{-1})^{1-{\beta_2}} 
+ \Phi^{-1}((t-s)^{-1})^{1-{\beta_2}+\gamma} 
\Phi^{-1}(s^{-1})^{-\gamma+{\beta_2}}
   \nn\\
 &\qquad
+\Phi^{-1}((t-s)^{-1})^{1-{\beta_2}+\gamma}\Phi^{-1}(s^{-1})^{-{\beta_2}}
 +(t-s)s^{-1}\Phi^{-1}((t-s)^{-1}) 
\Phi^{-1}(s^{-1})^{-\gamma}  ds \Bigg). \label{e:works7}
\end{align}
Since $\beta+\delta_1>1$, $\delta_1>2/3>1/2$ and $\gamma< \delta_1+\beta_2-1$, 
using Lemma \ref{l:convoluton-inequality} we see that 
$\int_{t/2}^t \Phi^{-1}((t-s)^{-1})^{1-\beta_2}\, ds \le c_2t \Phi^{-1}(t^{-1})^{1-\beta_2}$, 
$ \int_{t/2}^t \Phi^{-1}((t-s)^{-1})^{1-\beta_2+\gamma}ds \le c_3 t \Phi^{-1}(t^{-1})^{1-\beta_2+\gamma}$
and $ \int_0^t (t-s) \Phi^{-1}((t-s)^{-1})ds \le c_4 t^2 \Phi^{-1}(t^{-1})$.
Thus, by Lemma \ref{l:convoluton-inequality}, \eqref{e:works7} is bounded above by $c_5 t\Phi^{-1}(t^{-1}) \rho(t, x-y)$.
Now, 
 \eqref{e:fract-der-p-kappa-2} follows immediately from this, \eqref{e:p-kappa}, \eqref{e:gradient-phi-y} and  Proposition \ref{p:upperestonp}. \qed
 
 We will also need the following 
corollary, 
 which follows from  \eqref{e:L-on-phi-y-cont}.
\begin{corollary}\label{l:L-p-continuity}
For $x\neq y$, the function $t\mapsto \LL^{\kappa}p^{\kappa}(t,x,y)$ is continuous on $(0, \infty)$.
\end{corollary}


\section{Proofs of main results}

\subsection{A nonlocal maximum principle}
We first establish a somewhat different version of \cite[Theorem 4.1]{CZ}.
\begin{thm}
\label{t:nonlocal-max-principle}
Suppose there exists a function $g:\R^d \to (0, \infty)$ such that \eqref{e:nshf} holds.
Let $T>0$ and   $u\in C_b([0,T]\times \R^d)$ be such that
\begin{equation}\label{e:nonlocal-max-principle-1}
\lim_{t\downarrow 0} \sup_{x\in \R^d} |u(t,x)-u(0,x)|=0\,, 
\end{equation}
and  for each $x\in \R^d$,
\begin{equation}\label{e:nonlocal-max-principle-2}
t\mapsto \LL^{\kappa}u(t,x) \ \textrm{is continuous on } (0,T].
\end{equation}
Suppose that  $u(t,x)$ satisfies the following inequality: for all $(t,x)\in (0,T]\times \R^d$,
\begin{equation}\label{e:nonlocal-max-principle-4}
\partial_t u(t,x)\le\LL^{\kappa}u(t,x)\, .
\end{equation}
Then for all $t\in (0,T)$,
\begin{equation}\label{e:nonlocal-max-principle-5}
\sup_{x\in \R^d}u(t,x)\le \sup_{x\in \R^d}u(0,x)\, .
\end{equation}
\end{thm}

\pf 
Choose $a>0$ such that 
\begin{align}
\label{e:Lglg}
\LL^{\kappa}g(x) \le a g(x), \quad \text{for all } x \in \R^d.
\end{align}
Let $\delta, \eps>0$ and $u^{\delta}_\eps (t,x):=u(t,x)-\delta(t-\eps+e^{at} g(x))$. Then by  \eqref{e:nonlocal-max-principle-4} and \eqref{e:Lglg}, for all 
$(t,x)\in (0,T]\times \R^d$, we have
\begin{align}\label{e:nonlocal-max-principle-4e}
&\partial_t u^{\delta}_\eps (t,x) 
=\partial_t u (t,x) -\delta(1+ae^{at} g(x))
\le  \LL^{\kappa} u (t,x) -\delta-\delta a e^{at} g(x)\nn\\
&= \LL^{\kappa}u^{\delta}_\eps (t,x) -\delta +\delta e^{at} (\LL^{\kappa}g(x) -a g(x)) \le \LL^{\kappa}u^{\delta}_\eps (t,x) -\delta. 
\end{align}
Since $u\in C_b([0,T]\times \R^d)$, 
by letting $\delta \to 0$ and $\eps \to 0$, 
it suffices to show that 
\begin{equation}\label{e:nonlocal-max-principle-51}
\sup_{x\in \R^d}u^{\delta}_\eps(t,x)\le \sup_{x\in \R^d}u^{\delta}_\eps (\eps, x), \quad t\in 
(\eps,T] \, .
\end{equation}
Fix $\delta, \eps>0$ and suppose that \eqref{e:nonlocal-max-principle-51} does not hold. 
Then, by the continuity of $u^{\delta}_\eps$ and the fact that $\lim_{x \to \infty}u^{\delta}_\eps (t, x) = -\infty$ 
(which is a consequence of \eqref{e:nshf}), 
there exist $t_0  \in (\eps,T]$  and  $x_0 \in \R^d$ such that 
\begin{equation}\label{e:nonlocal-max-principle-52}
\sup_{t  \in (\eps, T], x\in \R^d} u^{\delta}_\eps(t,x) = u^{\delta}_\eps(t_0,x_0).
\end{equation}
Thus by \eqref{e:nonlocal-max-principle-4e},  for $h \in  (0, t_0-\eps)$,
\begin{align*}
&0 \le \frac{1}{h} (u^{\delta}_\eps(t_0,x_0)-u^{\delta}_\eps(t_0-h,x_0))= \frac{1}{h} \int_{t_0-h}^{t_0}  \partial_t u^{\delta}_\eps(s,x_0)ds
 \le \frac{1}{h} \int_{t_0-h}^{t_0}  \LL^{\kappa}u^{\delta}_\eps(s,x_0)ds -\delta.
\end{align*}
Letting $h \to 0$ and using \eqref{e:nonlocal-max-principle-2} and \eqref{e:nonlocal-max-principle-52} we get
\begin{align*}
0 &\le  \LL^{\kappa}u^{\delta}_\eps(t_0,x_0) -\delta\\
&=
\textrm{p.v.} \int_{\R^d}\left(u^{\delta}_\eps(t_0,x_0+z)-u^{\delta}_\eps(t_0,x_0)\right)\kappa(x_0,z)J(z)\, dz -\delta \le -\delta, 
\end{align*}
which gives a contradiction. Therefore \eqref{e:nonlocal-max-principle-51} holds.
 \qed
 
\begin{remark}\label{r:sufjeps}
{\rm
Suppose that $\int_{|z|>1} |z|^\eps j(|z|)dz < \infty$ for some $\eps>0$. 
Let 
$g(x) = (1+|x|^2)^{\eps/2}$.
Note that 
\begin{align}
\label{e:dif}
|\partial_{i,j} g(x)| \le c_1(1+|x|)^{\eps-2}, \quad i, j=1, \dots, d. 
\end{align}
By \eqref{e:dif} and \eqref{e:Ppsi}, we have that 
for $|x| \le 1$,
\begin{align}
&|\LL^{\kappa}g(x)| \le \gamma_0 \int_{|z| \le 1}|\delta_g(x;z)|  j(|z|)dz +  \gamma_0 g(x) \int_{|z| > 1}  j(|z|)dz+\gamma_0 \int_{|z| > 1} g(x \pm z) j(|z|)dz\nn\\
&\le c_2\left(\int_{|z| \le 1}|z|^2 j(|z|)dz +  \int_{|z| > 1}  j(|z|)dz+ \int_{|z| > 1} |z|^\eps j(|z|)dz \right) \le c_3 \le c_3g(x).
\end{align}
If $|x| >1$,
then by \eqref{e:dif} and \eqref{e:Ppsi},
\begin{align}
&|\LL^{\kappa}g(x)| \le \gamma_0 \int_{|z| \le |x|}|\delta_g(x;z)|  j(|z|)dz +  \gamma_0 g(x) \int_{|z| >  |x|}  j(|z|)dz+\gamma_0 \int_{|z| >  |x|} g(x \pm z) j(|z|)dz\nn\\
&\le c_3\left(\int_{|z| \le |x|} |x|^{\eps-2} |z|^2 j(|z|)dz +  g(x) \int_{|z| > 1}  j(|z|)dz+ \int_{|z| > |x|} |z|^\eps j(|z|)dz \right)\nn\\
&\le c_4\left(|x|^{\eps} \int_{\R^d} ((|z|/|x|)^2 \wedge 1) j(|z|)dz +  g(x)+ 1 \right) \le c_5 g(x).
\end{align}
Therefore $g$ satisfies \eqref{e:nshf}.
}
\end{remark}

\subsection{Properties of the semigroup  $(P^{\kappa}_t)_{t\ge 0}$}

Define
$$
P_t^{\kappa}f(x)=\int_{\R^d}p^\kappa(t,x, y)f(y)dy. $$

\begin{lemma}\label{l:L-int-commute0}
For any bounded  function $f$, we have
\begin{equation}\label{e:L-int-commute-2}
\LL^{\kappa}P_t^{\kappa} f(x)=\int_{\R^d}\LL^{\kappa}p^{\kappa}(t,\cdot, y)(x)f(y)dy\, .
\end{equation}
\end{lemma}
\pf
By the same computation as in the proof of \eqref{e:fract-der-est1} we have that for all 
$t \le T$, $T\ge 1$,  and $\eps>0$,
\begin{align*}
&t \int_{|z|>\eps}\rho(t, x \pm z)j(|z|)\, dz \\
&\le 
\int_{\Phi^{-1}(t^{-1})|z|\le 1, |z|>\eps} t \rho(t,x\pm z)j(|z|)\, dz + \int_{\Phi^{-1}(t^{-1})|z|> 1}t \rho(t,x\pm z)j(|z|)\, dz\\
&\le c_1  4^{d+1}   t \rho(t,x)  \int_{ |z|>\eps}
j(|z|)\, dz + c_1\rho (t.x), 
\end{align*}
thus  by  Lemma \ref{l:p-kappa-difference}(1),
\begin{align*}
&\int_{\R^d}\left(\int_{|w|>\eps}\left|p^\kappa(t,x\pm w, y) -2p^\kappa(t,x, y)\right| 
\kappa(x,w)J(w)\, dw\right)  dy \\
\le&  2 \gamma_0 \kappa_1 \int_{\R^d}\int_{|w|>\eps}|p^\kappa(t,x, y)|
j(|w|)\, dwdy
+  \gamma_0 \kappa_1 \int_{\R^d}\int_{|w|>\eps}|p^\kappa(t,x\pm w, y)| 
j(|w|)\, dw dy\\
\le&  c_2 t \left( \int_{|w|>\eps}  j(|w|)\, dw\right)  \int_{\R^d}  \rho(t, x-y) 
 dy
+ c_2  t\int_{\R^d}\left(\int_{|w|>\eps}\rho(t, x\pm w- y)
j(|w|)\, dw\right)  dy\\
 <&\infty.
\end{align*}
Thus
by Fubini's theorem, for all for bounded  function $f$ and $\eps \in (0,1]$,
\begin{align*}
&\LL^{\kappa, \eps}P_t^{\kappa} f(x)=  
\int_{\R^d}\LL^{\kappa, \eps}p^{\kappa}(t,\cdot, y)(x)f(y)dy.
\end{align*}
Now, \eqref{e:L-int-commute-2} follows from this,  \eqref{e:fract-der-p-kappa-1} and the dominated convergence theorem. \qed

The following result is the counterpart  of \cite[Lemma 4.4]{CZ}.

\begin{lemma}\label{l:continuity-of-LP}
(a) For any $p\in [1,\infty]$, there exists a constant $c=c(p,d,\delta_1,\beta, \kappa_0,\kappa_1, \kappa_2)>0$ such that for all $f\in L^p(\R^d)$ and $t>0$,
\begin{equation}\label{e:LP-p-estimate}
\|\LL^{\kappa} P_t^{\kappa}f\|_p \le ct^{-1} \|f\|_p\, .
\end{equation}
(b)  If $f\in L^{\infty}(\R^d)$, $t\mapsto \LL^{\kappa} P_t^{\kappa}f$ is a continuous function on 
$(0,\infty)$.

\noindent
(c) For any $p\in [1,\infty)$ and $f\in L^p(\R^d)$, $t\mapsto \LL^{\kappa} P_t^{\kappa} f$ is continuous from  
$(0,\infty)$
into $L^p(\R^d)$. \end{lemma}
\pf 
(a) Let $p\in [1,\infty]$. 
By
\eqref{e:L-int-commute-2},  Lemma \ref{l:fract-der-p-kappa}, Young's inequality and Lemma \ref{l:convolution}(a), we have that for all
$f \in L^p(\R^d) \cap L^{\infty}(\R^d)$,
\begin{eqnarray*}
\|\LL^{\kappa} P_t^{\kappa}f\|_p&\le & c_1 \left(\int_{\R^d} \left| \int_{\R^d} \rho(t,x-y) |f(y)|\, dy\right|^p dx\right)^{1/p}\\
&\le & c_1 \|\rho(t, \cdot)\|_1\, \|f\|_p \le c_2 t^{-1}\|f\|_p\, .
\end{eqnarray*}
Inequality \eqref{e:LP-p-estimate}
for $f \in L^p(\R^d)$ now follows by a standard density argument.

\noindent
(b)
For any $\eps\in (0,1)$, by Lemma \ref{l:fract-der-p-kappa} we have for $x\neq y$,
$$
\sup_{t\in (\eps,T)}\left|\LL^{\kappa}p^{\kappa}(t,x,y)\right|\le c \sup_{t\in (\eps,T)}\rho(t,x-y)\le c \rho(\eps, x-y)\, .
$$
Assume that $f$ is bounded and measurable.  By 
Corollary
\ref{l:L-p-continuity}, $t\mapsto \LL^{\kappa}p^{\kappa}(t,x,y)f(y)$ is continuous for $x\neq y$. By the above display, the family $\{\LL^{\kappa}p^{\kappa}(t,x,y)f(y):\, t\in (\eps,1)\}$ is bounded by the integrable function $\rho(\eps, x-y)|f(y)|$. Now it follows from the dominated convergence theorem and \eqref{e:L-int-commute-2} that $t\mapsto \LL^{\kappa}P_t^{\kappa}f(x)$ is continuous. 

\noindent
(c)
Let $p\in [1,\infty)$.  When $f \in L^p(\R^d) \cap L^{\infty}(\R^d)$, the claim
follows similarly as (b)  by using \eqref{e:L-int-commute-2} and the domination by the $L^p$-function $\int_{\R^d}\rho(\eps, x-y)f(y)\, dy$. 
The claim for $f \in L^p(\R^d)$ now follows by standard density argument and \eqref{e:LP-p-estimate}.
\qed

\begin{remark}\label{r:L-p-continuity}
{\rm
Note that Lemma \ref{l:continuity-of-LP} uses only the following properties of $p^{\kappa}(t,x,y)$:  \eqref{e:L-int-commute-2}, $|\LL^{\kappa}p^{\kappa}(t, \cdot,y)(x)|\le c_1(T) \rho(t,x-y)$ for $t \in (0, T]$ and $t\mapsto \LL^{\kappa}p^{\kappa}(t, \cdot,y)(x)$ is continuous on $(0,T]$. Moreover, Lemma \ref{l:L-int-commute0} uses only the following properties of $p^{\kappa}(t,x,y)$: $p^{\kappa}(t, \cdot,y)(x) \le c_2(T) t \rho(t,x-y)$ and  $|\LL^{\kappa, \eps}p^{\kappa}(t, \cdot,y)(x)|\le c_3(T) \rho(t,x-y)$ for $\eps \in [0,1]$ and $t \in (0, T]$.
}
\end{remark}

The following result is the counterpart  of \cite[Lemma 4.3]{CZ}.

\begin{lemma}\label{l:L-int-commute}
For any bounded  H\"older continuous function 
$f \in C^\eta_b(\R^d)$, 
we have
\begin{equation}\label{e:L-int-commute}
\LL^{\kappa}\left(\int_0^t P_s^{\kappa}f(\cdot)ds\right)(x)=\int_0^t \LL^{\kappa} P_s^{\kappa}f(x)ds\, , \quad x\in \R^d\, .
\end{equation}
\end{lemma}
\pf 
Define 
$$
T_tf(x)=\int_{\R^d}p_y(t, x-y)f(y)dy, \quad S_tf(x)=\int_{\R^d}q(t, x, y)f(y)dy
$$
and 
$$
R_tf(x)=\int^t_0T_{t-s}S_sf(x)ds.
$$
Then, by Fubini's theorem and \eqref{e:q-estimate}, 
for all for bounded  function $f$,
\begin{align}
\label{e:PTR}
P^\kappa_tf(x)=T_tf(x)+R_tf(x).
\end{align}

We now assume $\eps  \in (0,1]$ and 
$0<s < t \le T$, $T\ge 1$.
Suppose that $|f(x)-f(y)|\le c_1(|x-y|^{\eta}\wedge 1)$. Without loss of generality we may and will assume that $\eta<\beta$. 
By Fubini's theorem, \eqref{e:psi1}, \eqref{e:intro-kappa} and  \eqref{e:fract-der-est1},
$$
\LL^{\kappa, \eps}T_tf(x)= \int_{\R^d}{\LL}^{\kappa, \eps}p_z(s, \cdot)(x-z)f(z)\, dz.$$
Thus, 
\begin{eqnarray*}
 |{\LL}^{\kappa, \eps} T_sf(x)|
&\le & \int_{\R^d}\left(\int_{|w|>\eps}|\delta_{p_z}(s,x-z;w)|
\kappa(x,w)J(w)\, dw\right)\, |f(z)-f(x)|\, dz\\
& & +\left|\int_{\R^d}\left(\int_{|w|>\eps}\delta_{p_z}(s,x-z;w)
\kappa(x,w)J(w)\, dw\right)\, dz\right| |f(x)|\,.
\end{eqnarray*}
By using \eqref{e:psi1}, \eqref{e:fract-der-est1}, \eqref{e:some-estimates-2a} and  \eqref{e:convolution-integrability}, 
for any $\beta_1\in(0, \delta_1)\cap (0, \beta]$,
$ |{\LL}^{\kappa, \eps} T_sf(x)|$ is bounded by 
\begin{eqnarray*}
&& c_1 \int_{\R^d}\rho(s,x-z)\left(|x-z|^{\eta}\wedge 1\right)\, dz +
c_1\ s^{-1}\Phi^{-1}(s^{-1})^{-\beta_1}\\
&\le & c_2\ s^{-1}\Phi^{-1}(s^{-1})^{-\eta}+
c_1\ s^{-1}\Phi^{-1}(s^{-1})^{-\beta_1}\, ,
\end{eqnarray*}
and the right hand side is integrable by Lemma \ref{l:convoluton-inequality}. 
Thus by the dominated convergence theorem and Fubini's theorem, 
\begin{align}
{\LL}^{\kappa}\int_0^t T_sf(x)\, ds
=\lim_{\eps\downarrow 0} {\LL}^{\kappa, \eps} \int_0^t T_sf(x)\, ds=\int_0^t \lim_{\eps\downarrow 0} {\LL}^{\kappa, \eps} T_sf(x) ds =
\int_0^t {\LL}^{\kappa} T_s f(x)\, ds\, .\label{e:rs3}
\end{align}

It follows from \eqref{e:difference-q-estimate}, \eqref{e:convolution-integrability} and the boundedness of $f$  that
for any $\beta_2\in (0, \beta]\cap(0, \delta_1/2)$ and $\gamma\in (0, \beta_2)$,
we have
\begin{equation}\label{e:rs1}
|S_sf(x)-S_sf(x')|\le c_3 s^{-1}\Phi^{-1}(s^{-1})^{-\gamma}
\left(|x-x'|^{\beta_2-\gamma}\wedge 1\right).
\end{equation}
It follows from \eqref{e:q-estimate}, \eqref{e:convolution-integrability} and the boundedness of $f$
that
\begin{equation}\label{e:rs2}
|S_sf(x)|\le c_4 s^{-1}\Phi^{-1}(s^{-1})^{-\beta_2}.
\end{equation}
We use Lemma \ref{l:p-kappa-difference}(1) and Fubini's theorem  in the first line below, which can be justified by an argument similar to \eqref{e:FubiniJ1} and  \eqref{e:FubiniJ2}:
\begin{eqnarray*}
&&|{\LL}^{\kappa, \eps} R_s f(x)|\\
&\le&\int^s_0\left|\int_{\R^d}\left(\int_{|w|>\eps}\delta_{p_z}(s-r,x-z;w)
\kappa(x,w)J(w)\, dw\right) S_rf(z)\, dz \right| dr \\
&\le & \int^s_0\int_{\R^d}\left(\int_{|w|>\eps}|\delta_{p_z}(s-r,x-z;w)|
\kappa(x,w)J(w)\, dw\right)\, |S_rf(z)-S_rf(x)|\, dzdr\\
& & +\int^s_0\left|\int_{\R^d}\left(\int_{|w|>\eps}\delta_{p_z}(s-r,x-z;w)
\kappa(x,w)J(w)\, dw\right)\, dz\right| |S_rf(x)|dr\,.
\end{eqnarray*}
By using \eqref{e:psi1}, \eqref{e:fract-der-est1}, \eqref{e:some-estimates-2a}, \eqref{e:convolution-integrability}, \eqref{e:rs1}, \eqref{e:rs2} and Lemma \ref{l:convoluton-inequality},
we further have that
\begin{align*}
&|{\LL}^{\kappa, \eps} R_s f(x)|\,\le \, c_5 \int^s_0\int_{\R^d}\rho(s-r,x-z)r^{-1}\Phi^{-1}(r^{-1})^{-\gamma}
\left(|x-z|^{\beta_2-\gamma}\wedge 1\right)\, dzdr\\
& \qquad +c_5\int^s_0\ r^{-1}\Phi^{-1}(r^{-1})^{-\beta_2}dr\\
&\le c_6\int^s_0(s-r)^{-1}\Phi^{-1}((s-r)^{-1})^{-(\beta_2-\gamma)}r^{-1}
\Phi^{-1}(r^{-1})^{-\gamma}dr+
c_5\int^s_0\ r^{-1}\Phi^{-1}(r^{-1})^{-\beta_2}dr\\
&\le c_7s^{-1}\Phi^{-1}(s^{-1})^{-\beta_2}+c_5\Phi^{-1}(s^{-1})^{-\beta_2} 
= 2c_7 s^{-1}\Phi^{-1}(s^{-1})^{-\beta_2}\,,
\end{align*}
and the right hand side is integrable by Lemma \ref{l:convoluton-inequality}. This justifies the use of the dominated convergence theorem in the second line of the following calculation:
\begin{align}
&{\LL}^{\kappa}\int_0^t R_sf(x)\, ds =
\lim_{\eps\downarrow 0}
{\LL}^{\kappa, \eps} \int_0^t R_sf(x)\, ds
=\int_0^t \lim_{\eps\downarrow 0}  {\LL}^{\kappa, \eps} R_s f(x)\, ds  =
\int_0^t {\LL}^{\kappa} R_s f(x)\, ds\, .\label{e:rs4}
\end{align}

Combining \eqref{e:rs4} with \eqref{e:rs3} and \eqref{e:PTR}, we arrive at the conclusion of this lemma.
\qed

\subsection{Proofs of Theorems \ref{t:intro-main}--\ref{t:intro-semigroup}}

\noindent
{\bf Proof of Theorem \ref{t:intro-main}}.
By using Lemma \ref{l:phi-y-abs-cont} in the second equality, \eqref{e:q0-definition} in the third, \eqref{e:integral-equation} in the fourth, \eqref{e:q0-definition} in the fifth, and Lemma \ref{l:L-on-phi-y} in the sixth equality, we have
\begin{eqnarray*}
\partial_t p^{\kappa}(t,x,y)&=&\partial_t p_y(t,x-y)+\partial_t \phi_y(t,x)\\
&=&\LL^{{\mathfrak K}_y}p_y(t,x-y)+\left(q(t,x,y)+\int_0^t \int_{\R^d}\LL^{{\mathfrak K}_z}p_z(t-s,\cdot)(x-z)q(s,z,y)\, dz\, ds\right)\\
&=&\left(\LL^{{\mathfrak K}_x}p_y(t,x-y)-q_0(t,x,y)\right)\\
& &+\left(q(t,x,y)+\int_0^t \int_{\R^d}\LL^{{\mathfrak K}_z}p_z(t-s,\cdot)(x-z)q(s,z,y)\, dz\, ds\right)\\
&=&\LL^{{\mathfrak K}_x}p_y(t,x-y)+\int_0^t \int_{\R^d} q_0(t-s,x-z)q(s,z,y)\, dz\, ds\\
& &+\int_0^t \int_{\R^d}\LL^{{\mathfrak K}_z}p_z(t-s,\cdot)(x-z)q(s,z,y)\, dz\, ds \\
&=&\LL^{{\mathfrak K}_x}p_y(t,x-y)+\int_0^t \int_{\R^d}\LL^{{\mathfrak K}_x}p_z(t-s,\cdot)(x-z)q(s,z,y)\, dz\, ds\\
&=& \LL^{\kappa}p^{\kappa}(t,x,y)\, .
\end{eqnarray*}
Thus \eqref{e:intro-main-1} holds. The joint continuity of  $p^{\kappa}(t,x,y)$ 
is proved in Lemma \ref{l:continuity-p-kappa}.
Further,  if we apply the maximum principle, Theorem \ref{t:nonlocal-max-principle}, to $u_f(t,x):=P_t^{\kappa}f(x)$ with $f\in C_c^{\infty}(\R^d)$ and $f\le 0$, we get $u_f(t,x)\le 0$ for all $t\in(0,T]$ and all  $x\in \R^d$. This implies that $p^{\kappa}(t,x,y)\ge 0$.

\noindent 
(i) \eqref{e:intro-main-2} is proved in 
Lemma \ref{l:p-kappa-difference}(1).

\noindent
(ii)  
The estimate \eqref{e:intro-main-4} is given in \eqref{e:fract-der-p-kappa-1}, 
while continuity of $t\mapsto \LL^{\kappa}p^{\kappa}(t,\cdot,y)(x)$ is proven in 
Corollary \ref{l:L-p-continuity}. 

\noindent
(iii) Let $f$ be a bounded and uniformly continuous function. For any $\eps>0$, there exists $\delta>0$ such that $|f(x)-f(y)|<\eps$ for all $|x-y|<\delta$. By \eqref{e:some-estimates-2c}, \eqref{e:intro-psibound}, \eqref{e:convolution-integrability}  and the estimate for $p_y(t,x-y)$  in Proposition \ref{p:upperestonp} we have
\begin{eqnarray*}
\lefteqn{\lim_{t\downarrow 0}\sup_{x\in \R^d} \left|\int_{\R^d}p_y(t,x-y)f(y)\, dy -f(x)\right|}\\
&=&\lim_{t\downarrow 0}\sup_{x\in \R^d} \left|\int_{\R^d}p_y(t,x-y)f(y)\, dy -\int_{\R^d}p_y(t,x-y)f(x)\, dy\right|\\
&\le & c_1 \lim_{t\downarrow 0}\sup_{x\in \R^d} \int_{\R^d} t\rho(t,x-y)\, |f(y)-f(x)|\, dy\\
&\le &\eps c_1 \lim_{t\downarrow 0}\sup_{x\in \R^d}\int_{|x-y|<\delta}t\rho(t,x-y)dy +2c_1\|f\|_{\infty}\lim_{t\downarrow 0}\sup_{x\in \R^d}\int_{|x-y|\ge\delta}t\rho(t,x-y)dy\\
&\le & c_2 \eps \lim_{t\downarrow 0}\sup_{x\in \R^d} \int_{\R^d}t \rho(t,x-y)dy +2c_1\|f\|_{\infty}\lim_{t\downarrow 0}t \sup_{x\in \R^d}\int_{|x-y|\ge \delta}\frac{\Phi(|x-y|^{-1})}{|x-y|^d}dy\\
&\le & c_2 \eps +2c_1\|f\|_{\infty}\lim_{t\downarrow 0}t \int_{|z|\ge \delta} \frac{\Phi(|z|^{-1})}{|z|^d} dz = c_2 \eps \, .
\end{eqnarray*}
This implies that
\begin{align}\label{e:main-proof-iv}
\lim_{t\downarrow 0}\sup_{x\in \R^d} \left|\int_{\R^d}p_y(t,x-y)f(y)\, dy -f(x)\right|=0\, .
\end{align}
Further, by \eqref{e:main-proof-i} and \eqref{e:convolution-integrability},
for any $\beta_2\in (0, \beta]\cap (0, \delta_1)$, we have
\begin{eqnarray*}
\lefteqn{\left|\int_{\R^d}\int_0^t \int_{\R^d}p_z(t-s,x-z)q(s,z,y)\, dz\, ds\,f(y) dy\right|}\\
&\le & c_3\, \|f\|_{\infty}\,  t\int_{\R^d} 
\left(\rho_0^{\beta_2}+\rho_{\beta_2}^0\right)(t,x-y)\, dy\, 
\le \, c_4\, \Phi^{-1}(t^{-1})^{-\beta_2} \longrightarrow 0\, ,
\quad t\downarrow 0\, .
\end{eqnarray*}
The claim now follows from this, \eqref{e:p-kappa} and \eqref{e:main-proof-iv}.

\medskip
\noindent
{\it Uniqueness of the kernel satisfying \eqref{e:intro-main-1}-\eqref{e:intro-main-5}}:
Let $\wt{p}^{\kappa}(t,x,y)$ be another non-negative jointly continuous kernel satisfying \eqref{e:intro-main-1}--\eqref{e:intro-main-5}.
For any function
$f\in C_c^{\infty}(\R^d)$, define $\wt{u}_f(t,x):= \int_{\R^d}\wt{p}^{\kappa}(t,x,y)f(y)\, dy$. By the joint continuity of $\wt{p}^{\kappa}(t,x,y)$, (i) and (iii) we have that
$$
\wt{u}_f\in C_b([0,T]\times \R^d), \qquad \lim_{t\downarrow 0} \sup_{x\in \R^d} |\wt{u}_f(t,x)-f(x)|=0\, .
$$
 By Lemma \ref{l:L-int-commute0}
 and Remark \ref{r:L-p-continuity}, 
 \begin{align}
 \label{e:L-int-commute-2n}
\LL^{\kappa} \wt{u}_f(t,x)=\int_{\R^d}\LL^{\kappa} \wt{p}^{\kappa}(t,x,y)f(y)\, dy \quad \text{and} \quad \LL^{\kappa} {u}_f(t,x)=\int_{\R^d}\LL^{\kappa}  {p}^{\kappa}(t,x,y)f(y)\, dy.
 \end{align}
 Moreover, by Lemma \ref{l:continuity-of-LP} and Remark \ref{r:L-p-continuity}, $t\mapsto \LL^{\kappa}u_f(t,x)$ and $t\mapsto \LL^{\kappa}\wt{u}_f(t,x)$ are continuous on $(0,T]$. 
Here and in \eqref{e:L-int-commute-2n} we use that $\wt{p}^{\kappa}$ satisfies (i)--(ii).

Let $w(t,x):=u_f(t,x)-\wt{u}_f(t,x)$. Then $w(0,x)=0$, $\lim_{t\downarrow 0}\sup_{x\in \R^d} |w(t,x)-w(0,x)|=0$, 
and $t\mapsto \LL^{\kappa}w(t,x)$ is continuous on $(0,T]$. 
Note that by \eqref{e:intro-main-4} and \eqref{e:intro-main-1},
$$
|\partial_t  p^{\kappa}(t,x,y)|+ |\partial_t  \wt{p}^{\kappa}(t,x,y)| \le c_5 \rho(t,x-y), \quad t \in (0, T]\,.
$$
Thus, by the dominated convergence theorem,
$$
\partial_t \wt{u}_f(t,x)=\int_{\R^d}\partial_t \wt{p}^{\kappa}(t,x,y)f(y)\, dy \quad \text{and} \quad \partial_t  {u}_f(t,x)=\int_{\R^d}\partial_t  {p}^{\kappa}(t,x,y)f(y)\, dy.
$$
By this,  \eqref{e:intro-main-1} and  \eqref{e:L-int-commute-2n}, we have 
 $\partial_t w(t,x)=\LL^{\kappa}w(t,x)$. Hence, 
all the assumptions
 of Theorem \ref{t:nonlocal-max-principle} are satisfied and we can conclude that for every $t\in(0,T]$, $\sup_{x\in \R^d}w(t,x)\le \sup_{x\in \R^d}w(0,x)=0$. By applying the theorem to $-w$ we get that $w(t,x)=0$ for all $t\in(0,T]$ and every $x\in \R^d$. Hence, $u_f=\wt{u}_f$ for every $f\in C_c^{\infty}(\R^d)$, which implies that $\wt{p}^{\kappa}(t,x,y)=p^{\kappa}(t,x,y)$.

The last statement of the theorem about the dependence of constants $c_1$ and $c_2$ has been already proved in the results above.
\qed

\noindent
{\bf Proof of Theorem \ref{t:intro-further-properties}}. (1)
The constant function $u(t,x)=1$ solves $\partial_t u(t,x)=\LL^{\kappa}u(t,x)$, hence applying Theorem \ref{t:nonlocal-max-principle} to $\pm(P_t^{\kappa} 1(x)-1)$ we get that $P_t^{\kappa} 1(x)\equiv 1$ proving \eqref{e:intro-main-6}.

\noindent
(2) Same as the proof of \cite[Theorem 1.1(3)]{CZ}.

\noindent
(3) 
By \eqref{e:intro-main-1} and \eqref{e:intro-main-4} we see that $\left|\partial_t p^{\kappa}(t,x,y)\right| \le c_2 \rho(t,x-y)$ for $t\in(0,T]$ and $x\neq y$. Hence by the mean value theorem, for $0<s\le t\le T$ and $x\neq y$,
\begin{align}
\left| p^{\kappa}(s,x,y)-p^{\kappa}(t,x,y)\right| \le c_2 |t-s|\rho(s,x-y)\, .\label{e:lip-cont-new5}
\end{align}
Let 
$\gamma\in (0,\delta_1) \cap (0,1]$.
By Lemma \ref{l:p-kappa-difference} and by the definition of $\rho_{-1}^0$, we have that for every $t\in(0,T]$ ,
\begin{eqnarray}
|p^{\kappa}(t,x,y)-p^{\kappa}(t,x',y)|&\le &c_1 |x-x'|^\gamma \Phi^{-1}(t^{-1}) t\left(\rho(t,x-y)+\rho(t,x'-y)\right) \nonumber \\
&\le & 2c_1 |x-x'|^\gamma \Phi^{-1}(t^{-1}) t\left(\rho(t,x-y)\vee \rho(t,x'-y)\right)\, .\label{e:lip-cont-new}
\end{eqnarray}
By use of the
triangle inequality, this  together with 
\eqref{e:lip-cont-new5} implies the first claim.

By \eqref{e:intro-main-2}, if $\Phi^{-1}(t^{-1})|x-x'|\ge 1$, 
\begin{align}
&|p^{\kappa}(t,x,y)-p^{\kappa}(t,x',y)|\le 
p^{\kappa}(t,x,y)+p^{\kappa}(t,x',y) \le 
c_1 t\left(\rho(t,x-y)+\rho(t,x'-y)\right) \nonumber \\
&\le  2c_1 |x-x'| \Phi^{-1}(t^{-1}) t\left(\rho(t,x-y)\vee \rho(t,x'-y)\right)\, .\label{e:lip-cont-new2}
\end{align}

Suppose $\Phi^{-1}(t^{-1})|x-x'|\ge 1$, $\beta+\delta_1 >1$ and $\delta_1 \in (2/3,2)$. Then by 
\eqref{e:fract-der-p-kappa-2}
\begin{align}
&|p^{\kappa}(t,x,y)-p^{\kappa}(t,x',y)| \le |x-x'|\cdot \int_0^1|\nabla p(t, x+\theta (x'-x), y)| \, d\theta
\nn\\
&\le c t\Phi^{-1}(t^{-1})  |x-x'|\int_0^1  \rho(t,(x-y)+\theta (x'-x))
d \theta. \label{e:lip-cont-new21}
\end{align}
Since $\theta|x'-x|\le 1/\Phi^{-1}(t^{-1})$, from \eqref{e:lip-cont-new21} we have 
\begin{align}
&|p^{\kappa}(t,x,y)-p^{\kappa}(t,x',y)| \le c t\Phi^{-1}(t^{-1})  |x-x'|\rho(t,x-y)\nn\\
 &\le c t\Phi^{-1}(t^{-1})  |x-x'| \left(\rho(t,x-y)\vee \rho(t,x'-y)\right). \label{e:lip-cont-new22}
\end{align}
\eqref{e:lip-cont-new5}, \eqref{e:lip-cont-new2} and \eqref{e:lip-cont-new22} imply the second claim.

\noindent
(4)
This follows immediately from the second part of Lemma \ref{l:fract-der-p-kappa}.
\qed

\noindent 
{\bf Proof of Theorem \ref{t:intro-semigroup}}. (1)
We first claim that for $f\in C_b^{2,\eps}(\R^d)$, $\LL^{\kappa} f$ is bounded H\"older continuous. We will use results from \cite{B}.
For $f\in C_b^{2,\eps}(\R^d)$ and $x, z \in \R^d$, let
$$
E_zf(x)=f(x+z)-f(x) \quad \text{and} \quad F_zf(x)=f(x+z)-f(x) -\nabla f(x) \cdot z.
$$
Using the assumption that $\kappa (y, z)=\kappa (y,-z)$, we have
\begin{align*}
\LL^{{\mathfrak K}_y} f(x) =\int_{|z|<1}  F_zf(x)  \kappa (y,z) J(z)  dz +\int_{|z| \ge 1}  E_zf(x)  \kappa (y,z) J(z)  dz.
\end{align*}
Thus,   $\LL^{\kappa} f$ is bounded by \eqref{e:psi1} and \eqref{e:intro-kappa}.
Moreover,   using \eqref{e:intro-kappa-holder},  \eqref{e:psi1} and \cite[Theorem 5.1 (b) and (e)]{B} with $\gamma=2+\eps$,
\begin{align*}
&|\LL^{\kappa} f(x)-\LL^{\kappa} f(y) |\\
\le& 
|\int_{\R^d} \delta_f(x ; z) (\kappa (x,z)-\kappa (y,z)) J(z)dz|+
| \LL^{{\mathfrak K}_y} f(x)-\LL^{{\mathfrak K}_y} f(y) |\\
\le&c_1(|x-y|^{\beta}\wedge 1)
\int_{\R^d}(|z|^2 \wedge 1)  j(|z|) dz +c_1\int_{|z|<1}  |F_zf(x)-F_zf(y)|  \kappa (y,z) j(|z|) dz  \\
&+c_1\int_{|z| \ge 1}  |E_zf(x)-E_zf(y)|  \kappa (y,z) j(|z|) dz\\
\le& c_2 |x-y|^{\beta}+c_2\left( \int_{|z|<1}  |z|^2 j(|z|) dz  \right)
|x-y|^{\eps}
+c_2 \left(\int_{|z| \ge 1}  j(|z|) dz\right)  |x-y|. 
\end{align*}
Thus we have proved the claim.

For $f\in C_b^{2,\eps}(\R^d)$,  we define
$
u(t,x):=f(x)+\int_0^tP_s^{\kappa}\LL^{\kappa}f(x)\, ds\, .
$
Note that 
$$
|u(t,x)-u(0,x)|\le 
 \int_0^t |P_s^{\kappa}\LL^{\kappa}f(x)| ds \le t \|\LL^{\kappa}f\|_\infty. 
$$
Thus \eqref{e:nonlocal-max-principle-1} holds.
Since $\LL^\kappa f$  is bounded H\"older continuous, 
we can use  \eqref{e:L-int-commute} (together with \eqref{e:intro-main-4},  \eqref{e:intro-main-1} and  
\eqref{e:L-int-commute-2n})
to get $\LL^{\kappa} P_s^{\kappa}\LL^{\kappa} f(x)=\partial_s \left(P_s^{\kappa}\LL^{\kappa}f\right)(x)$) and 
obtain
\begin{eqnarray*}
\LL^{\kappa}u(t,x)&=&\LL^{\kappa} f(x)+\int_0^t \LL^{\kappa} P_s^{\kappa}\LL^{\kappa} f(x)\, ds\\
&=&\LL^{\kappa} f(x)+\int_0^t \partial_s \left(P_s^{\kappa}\LL^{\kappa}f\right)(x)\, ds= P_t\LL^{\kappa} f(x)=\partial_t u(t,x)\, .
\end{eqnarray*}
Therefore  $u(t,x)$ satisfies the assumptions of Theorem \ref{t:nonlocal-max-principle}. Since $u(0,x)=f(x)$, it follows from the maximum principle that
\begin{equation}\label{e:main-proof-7}
P_t^{\kappa}f(x)=u(t,x)=f(x)+\int_0^tP_s^{\kappa}\LL^{\kappa}f(x)\, ds\, .
\end{equation}
Since $\LL^{\kappa}f$ is bounded and uniformly continuous,  we can use \eqref{e:intro-main-5} to get
$$
\lim_{t\downarrow 0}\frac{1}{t}\left(P_t^{\kappa}f(x)-f(x)\right)=\lim_{t\downarrow 0}\frac{1}{t}\int_0^t P_s^{\kappa}\LL^{\kappa}f(x)ds=\LL^{\kappa}f(x)
$$
and the convergence is  uniform.

\noindent
(2) Using our Theorem \ref{t:intro-main}(iii), Theorem \ref{t:intro-further-properties}(1)
and Lemma \ref{l:continuity-of-LP}, the proof of this part is the same as in \cite{CZ}. 
\qed

\subsection{Lower bound estimate of $p^{\kappa}(t,x,y)$}

By Theorem \ref{t:intro-semigroup}, we have that $(P_t^{\kappa})_{t  \ge 0}$ is a Feller semigroup and there 
exists 
a Feller process $X=(X_t, \P_x)$ corresponding to $(P_t^{\kappa})_{t  \ge 0}$. Moreover, by \eqref{e:main-proof-7}
for $f\in C_b^{2,\eps}(\R^d)$, 
\begin{align}
\label{e:MG}
f(X_t)-f(x)-\int_0^t \LL^{\kappa}f(X_s)\, ds
\end{align}
is a martingale with respect to the filtration $\sigma(X_s, s \le t)$.
 Therefore by the same argument as that in \cite[Section 4.4]{CZ}, we have 
the following L\'evy system formula:
 for every function $f:\R^d \times \R^d \to [0,\infty)$ vanishing on the diagonal and every stopping time $S$,
\begin{align}
\label{e:LSF}
\E_x \sum_{0<s\le S} f(X_{s-}, X_s)=\E_x \int_0^S f(X_s,y)J_X(X_s,dy) ds\, ,
\end{align}
where $J_X(x,y):=\kappa(x,y-x)J(x-y)$.

For $A \in \sB(\R^d)$
we define
$\tau_A:=\inf\{t \ge 0: X_t \notin A\}.
$

The following result is the counterpart  of \cite[Lemma 4.6]{CZ}.

\begin{lemma}\label{l:exit-probability}
For each $\gamma\in (0,1)$ there exists $A=A(\gamma)>0$ such that for  every $r>0$,
\begin{equation}\label{e:exit-probability-1}
\sup_{x \in \R^d} \P_x\left(\tau_{B(x,r)}\le (A\Phi(1/(4r)))^{-1}\right)\le \gamma\, .
\end{equation}
\end{lemma}
\pf Without loss of generality, we take $x=0$. The constant $A$ will be chosen later.  Let $f\in C_b^\infty(\R^d)$ with $f(0)=0$ and $f(y)=1$ for $|y|\ge 1$. For any $r>0$ set $f_r(y)=f(y/r)$. By the definition of $f_r$ and the martingale property in \eqref{e:MG} we have
\begin{eqnarray}
\P_0\left(\tau_{B(0,r)}\le (A\Phi(1/(4r)))^{-1}\right)&\le &\E_0\left[f_{r}\left(X_{\tau_{B(0,r)}\wedge (A\Phi(1/(4r)))^{-1}}\right)\right] \nonumber \\
&=& \E_0\left(\int_0^{\tau_{B(0,r)}\wedge (A\Phi(1/(4r)))^{-1}} \LL^{\kappa} f_{r}(X_s)\, ds\right) \, .\label{e:exit-probability-2}
\end{eqnarray}
By the definition of $\LL^{\kappa}$, \eqref{e:intro-kappa} and  \eqref{e:psi1} we have
\begin{eqnarray*}
\lefteqn{|\LL^{\kappa}f_{r}(y)|=\frac12 \left|\int_{\R^d}\left(f_{r}(y+z)+f_{r}(y-z)-2f_{r}(y)\right)\kappa(y,z)J(z)\, dz\right|}\\
&\le &\frac{\kappa_1 \gamma_0 \|\nabla^2 f_{r}\|_{\infty}}{2}\int_{|z|\le  r}|z|^2 j(|z|)\, dz +2\kappa_1 \gamma_0\|f_{r}\|_{\infty}\int_{|z| >  r}j(|z|)\, dz\\
&\le & c_1\left(\frac{\|\nabla^2 f\|_{\infty}}{r^2} r^2 \sP(r)+\|f\|_{\infty}\sP(r)\right)\, \le \,   c_2 \, \Phi(r^{-1})\, ,
\end{eqnarray*}
where $c_2=c_2(\kappa_1,\gamma_0,f)$. Here the last inequality is a consequence of  \eqref{e:Ppsi}. Substituting in \eqref{e:exit-probability-2} we get that
$$
\P_0\left(\tau_{B(0,r)}\le (A\Phi(1/(4r)))^{-1}\right)\le c_2  \Phi(r^{-1}) (A\Phi(1/(4r)))^{-1} \le 4c_2A^{-1}\, .
$$
With $A=4c_2/ \gamma$ the lemma is proved. \qed

\medskip
\noindent

\noindent
{\bf Proof of Theorem \ref{t:intro-lower-bound}}.
Throughout the proof, 
we fix $T, M\ge 1$ 
and, without loss of generality, we assume that  $\Phi^{-1}(T^{-1})^{-1}=M$.  

By \cite[Theorem 2.4]{CKK} and the same argument as the one in \cite[Proposition 2.2]{CKS}
(see also 
\cite[Proposition 6.4(1)]{GKK} or \cite[Proposition 6.2]{CK}), 
 \eqref{e:intro-wsc}, \eqref{e:intro-wusc}, 
\eqref{e:intro-kappa}
  and \eqref{e:psi1} imply that 
 there exists a constant $c_0>0$ such that 
 \begin{align}
p_y(t,x)\ge c_0 \left(
\Phi^{-1}(t^{-1})^d \wedge t j(|x|) \right)  \quad (t,x,y) \in (0,T] \times B(0,4M) \times \R^d\, .
\end{align}
 Since by \cite[Lemma 3.2(a)]{KSV},
 \begin{align}
\label{e:jlowerb}
j\left( |x|\right) \ge c_1 |x|^{-d} \Phi(|x|^{-1}), \quad |x| \le 4M
 \end{align}
 for some $c_1 \in (0,1)$, by Proposition  \ref{p:p-q}
we have 
\begin{align}
\label{e:lowerbp}
p_y(t,x)\ge c_0 c_1 t\rho(t,x)  \quad (t,x,y) \in (0,T] \times B(0,4M) \times \R^d\, .
\end{align}

\noindent
(1) Let $\lambda=1/A$ where $A$ is  the constant from Lemma \ref{l:exit-probability} for $\gamma=1/2$. Then  for every $t >0$,
\begin{equation}\label{e:exit-probability-3}
\sup_{z \in \R^d} \P_z(\tau_{B(z,2^{-2}\Phi^{-1}(t^{-1})^{-1})}
\le \lambda t)\le \frac12\, .
\end{equation}

Let $t\in(0,T]$ and $|x-y|\le 3\Phi^{-1}(t^{-1})^{-1}$( so that $|x-y| \le 3M$). 
By \eqref{e:main-proof-i}
 we have that there exists a constant $c_2>0$ such that
\begin{eqnarray*}
\lefteqn{\int_0^t\int_{\R^d}p_z(t-s,x-z)q(s,z,y)\, dz\, ds \ge  -c_2 t\left(\rho_{\beta}^0+\rho_0^{\beta}\right)(t,x-y)}\\
&=& -c_2 t \left(\Phi^{-1}(t^{-1})^{-\beta}+|x-y|^{\beta}\wedge 1\right)\rho(t,x-y)\\
&\ge & -c_2 t \left(\Phi^{-1}(t^{-1})^{-\beta}+3^{\beta}\Phi^{-1}(t^{-1})^{-\beta}\right)\rho(t,x-y)\, .
\end{eqnarray*}
We choose $t_0\in (0,1)$ so that for all $t\in (0,t_0)$, $c_2 (1+3^{\beta})\Phi^{-1}(t^{-1})^{-\beta}\le c_1/2$. Together with \eqref{e:lowerbp} and \eqref{e:p-kappa} we 
conclude 
that for all $t\in (0,t_0)$ and all $x,y\in \R^d$ satisfying $|x-y|\le 3\Phi^{-1}(t^{-1})^{-1}$ we have
$$
p^{\kappa}(t,x,y)\ge \frac{c_1}{2} t\rho(t,x-y)\ge c_3 t \frac{\Phi\left(\frac{1}{\Phi^{-1}(t^{-1})}+\frac{3}{\Phi^{-1}(t^{-1})}\right)}{\left(\frac{1}{\Phi^{-1}(t^{-1})}+\frac{3}{\Phi^{-1}(t^{-1})}\right)^d}\ge c_4 \Phi^{-1}(t^{-1})^d\, .
$$
By \eqref{e:intro-main-7} and iterating $\lfloor T/t_0\rfloor +1$ times, we obtain the following near-diagonal lower bound
\begin{equation}\label{e:lower-bound-1}
p^{\kappa}(t,x,y)\ge c_5 \Phi^{-1}(t^{-1})^d \ \ \textrm{for all }t\in(0,T] \textrm{ and }|x-y|\le 3\Phi^{-1}(t^{-1})^{-1}\, .
\end{equation}

Now we assume $|x-y| > 3\Phi^{-1}(t^{-1})^{-1}$ and let $\sigma=\inf\{t\ge 0: X_t\in B(y,2^{-1}\Phi^{-1}(t^{-1})^{-1})\}$. 
By the strong Markov property and \eqref{e:exit-probability-3} we have
\begin{eqnarray}
&&\P_x\left(X_{\lambda t}\in B(y,\Phi^{-1}(t^{-1})^{-1})\right)\ge \P_x\left(\sigma\le \lambda t, \sup_{s\in [\sigma, \sigma+\lambda t]} |X_s-X_{\sigma}|<2^{-1}\Phi^{-1}(t^{-1})^{-1}\right) \nn\\
&=&\E_x\left(\P_{X_{\sigma}}\Big(\sup_{s\in[0,\lambda t]} |X_s-X_0|<2^{-1}\Phi^{-1}(t^{-1})^{-1}\Big); \sigma\le \lambda t\right)\nn\\
&\ge & \inf_{z\in  B(y,2^{-1}\Phi^{-1}(t^{-1})^{-1})} \P_z\big(\tau_{ B(z,2^{-1}\Phi^{-1}(t^{-1})^{-1})}>\lambda t\big) \P_x\big(\sigma \le \lambda t\big)\nn\\
&\ge &\frac12 \P_x\big(\sigma \le \lambda t\big) \, \ge \, \frac12 \P_x\left(X_{\lambda t \wedge \tau_{ B(x,\Phi^{-1}(t^{-1})^{-1})}}\in B(y,2^{-1}\Phi^{-1}(t^{-1})^{-1})\right)\, .\label{e:ineque00}
\end{eqnarray}
Since 
$$
X_s\notin B\left(y,2^{-1}\Phi^{-1}(t^{-1})^{-1}\right)\subset B\left(x,\Phi^{-1}(t^{-1})^{-1}\right)^c\,, \quad s<\lambda t \wedge \tau_{ B(x,\Phi^{-1}(t^{-1})^{-1})},
$$
we have 
$$
{\mathbf 1}_{X_{\lambda t \wedge \tau_{ B(x,\Phi^{-1}(t^{-1})^{-1})}}\in B(y,2^{-1}\Phi^{-1}(t^{-1})^{-1})}=\sum_{s\le \lambda t \wedge \tau_{ B(x,\Phi^{-1}(t^{-1})^{-1})}} {\mathbf 1}_{X_s\in B(y,2^{-1}\Phi^{-1}(t^{-1})^{-1}}\, .
$$
Thus, by the L\'evy system formula in \eqref{e:LSF} we have
\begin{align}
\lefteqn{\P_x\left(X_{\lambda t \wedge \tau_{ B(x,\Phi^{-1}(t^{-1})^{-1})}}\in B(y,2^{-1}\Phi^{-1}(t^{-1})^{-1})\right)} \nn\\
=&\ \E_x\left[\int_0^{\lambda t \wedge \tau_{ B(x,\Phi^{-1}(t^{-1})^{-1})}}\int_{B(y,2^{-1}\Phi^{-1}(t^{-1})^{-1})}J_X(X_s,u)\, du\, ds\right] \nn\\
\ge&\ \E_x\left[\int_0^{\lambda t \wedge \tau_{ B(x,6 \cdot 2^{-4}\Phi^{-1}(t^{-1})^{-1})}}\int_{B(y,2^{-1}\Phi^{-1}(t^{-1})^{-1})}\kappa_0 j(|X_s-u|)
{\bf 1}_{\{u:|X_s-u|<|x-y|\}}
\, du\, ds\right].\label{e:ineque1}
\end{align}

Let  $w$ be the point on the line connecting $x$ and $y$ (i.e., $|x-y|=|x-w|+|w-y|$) such that  
$|w-y|=7 \cdot 2^{-4}  \Phi^{-1}(t^{-1})^{-1}$. Then 
$B(w, 2^{-4}  \Phi^{-1}(t^{-1})^{-1}) \subset B(y, \, 2^{-1} \Phi^{-1}(t^{-1})^{-1})$.
Moreover, for every $(z,u) \in  B(x,6 \cdot 2^{-4} \Phi^{-1}(t^{-1})^{-1}) \times B(w, 2^{-4} \Phi^{-1}(t^{-1})^{-1}) $, we have 
\begin{align*}
|z-u| & \le |z-x|+|w-u|+|x-w| =|z-x|+|w-u|+ |x-y|- |w-y|  \nn\\
&  <  (6 \cdot 2^{-4}+ 2^{-4}) \Phi^{-1}(t^{-1})^{-1} +|x-y|- 7 \cdot 2^{-4}  \Phi^{-1}(t^{-1})^{-1} = |x-y|. 
\end{align*}
Thus 
\begin{align}\label{e:ineque}
B(w, 2^{-4}  \Phi^{-1}(t^{-1})^{-1}) \subset \{u:|z-u|<|x-y|\} \quad\text{for } z \in B(x,6 \cdot 2^{-4} \Phi^{-1}(t^{-1})^{-1}).
\end{align}
\eqref{e:ineque} and \eqref{e:exit-probability-3} imply that 
\begin{align}
&\E_x\left[\int_0^{\lambda t \wedge \tau_{ B(x,6 \cdot 2^{-4}\Phi^{-1}(t^{-1})^{-1})}}\int_{B(y,2^{-1}\Phi^{-1}(t^{-1})^{-1})} j(|X_s-u|){\bf 1}_{\{u:|X_s-u|<|x-y|\}}
\, du\, ds\right] \nn\\
\ge & \E_x\left[\lambda t \wedge \tau_{ B(x,6 \cdot 2^{-4}\Phi^{-1}(t^{-1})^{-1})}\right] \int_{B(w,2^{-4}\Phi^{-1}(t^{-1})^{-1})}   j\left( |x-y|\right)\, du \nn\\
\ge &\lambda t \P_x\left(\tau_{B(x,6 \cdot 2^{-4}\Phi^{-1}(t^{-1})^{-1})}\ge \lambda t\right) \big|B(w,2^{-4}\Phi^{-1}(t^{-1})^{-1})
\big|\,   j\left( |x-y|\right) \nn\\
\ge & c_6 t \Phi^{-1}(t^{-1})^{-d}\, j\left( |x-y|\right).\label{e:ineque2}
\end{align}
By combining \eqref{e:ineque00},  \eqref{e:ineque1} and \eqref{e:ineque2}
 we get that
\begin{equation}\label{e:lower-bound-2}
\P_x\left(X_{\lambda t}\in B(y,\Phi^{-1}(t^{-1})^{-1})\right)\ge  \frac{1}{2}c_6 t \Phi^{-1}(t^{-1})^{-d}\,  j\left( |x-y|\right)
\end{equation}
By \eqref{e:intro-main-7}, \eqref{e:lower-bound-1} and \eqref{e:lower-bound-2} we have
\begin{eqnarray*}
p^{\kappa}(t,x,y)&\ge & \int_{B(y,\Phi^{-1}(t^{-1})^{-1})}p^{\kappa}(\lambda t, x,z)p^{\kappa}((1-\lambda)t,z,y)\, dz\\
&\ge &\inf_{z\in B(y,\Phi^{-1}(t^{-1})^{-1})}p^{\kappa}((1-\lambda)t,z,y) \int_{B(y,\Phi^{-1}(t^{-1})^{-1})}p^{\kappa}(\lambda t, x,z)\, dz\\
&\ge &c_7  \Phi^{-1}(t^{-1})^d  t \Phi^{-1}(t^{-1})^{-d}\,  j\left( |x-y|\right)= c_7 t j\left( |x-y|\right).
\end{eqnarray*}
Combining this estimate with \eqref{e:lower-bound-1} we obtain 
\eqref{e:intro-main-11}.
Inequality \eqref{e:intro-main-111} 
follows from \eqref{e:intro-main-11}, Proposition \ref{p:p-q} and 
\eqref{e:jlowerb}.
\qed

\bigskip
\noindent
{\bf Acknowledgements:} 
We are grateful to Xicheng Zhang for several valuable comments, in particular for suggesting the improvement of the gradient estimate \eqref{e:intro-main-9}.
We also thank Karol Szczypkowski for pointing out some mistakes in an earlier version of this paper and Jaehoon Lee  for reading the manuscript and giving helpful comments.

\end{doublespace}

\bigskip
\noindent

\vspace{.1in}
\begin{singlespace}


\vskip 0.1truein

\parindent=0em

{\bf Panki Kim}

Department of Mathematical Sciences and Research Institute of Mathematics,

Seoul National University, Building 27, 1 Gwanak-ro, Gwanak-gu Seoul 08826, Republic of Korea

E-mail: \texttt{pkim@snu.ac.kr}

\bigskip

{\bf Renming Song}

Department of Mathematics, University of Illinois, Urbana, IL 61801,
USA

E-mail: \texttt{rsong@illinois.edu}

\bigskip

{\bf Zoran Vondra\v{c}ek}

Department of Mathematics, Faculty of Science, University of Zagreb, Zagreb, Croatia, and \\
Department of Mathematics, University of Illinois, Urbana, IL 61801,
USA

E-mail: \texttt{vondra@math.hr}

\end{singlespace}
\end{document}